\documentclass[11pt]{amsart}
\usepackage{amssymb,amsmath,amsthm,amscd}
\usepackage[all]{xy}

\numberwithin{equation}{section}

\newtheoremstyle{fancy1}{10pt}{10pt}{\itshape}{12pt}{\textsc\bgroup}{.\egroup}{8pt}{
}
\newtheoremstyle{fancy2}{10pt}{10pt}{}{12pt}{\itshape}{.}{8pt}{ }

\theoremstyle{fancy1}

\newtheorem{cor}[equation]{Corollary}
\newtheorem{lem}[equation]{Lemma}
\newtheorem{prop}[equation]{Proposition}
\newtheorem{thm}[equation]{Theorem}

\newtheorem{main}{Theorem}
\newtheorem*{main*}{Theorem}

\newtheorem*{cor*}{Corollary}

\renewcommand{\thetable}{\theequation}
\theoremstyle{fancy2}
\newtheorem{definition}[equation]{Definition}
\newtheorem{rem}[equation]{Remark}
\newtheorem*{rem*}{Remark}

\newcommand{\cref}[1]{Corollary~\ref{#1}}

\newcommand{\lref}[1]{Lemma~\ref{#1}}

\newcommand{\gt}{\theta}

\newcommand{\eps}{\varepsilon}

\newcommand{\RP}{\mathbb{R\mkern1mu P}}
\newcommand{\CP}{\mathbb{C\mkern1mu P}}
\newcommand{\HP}{\mathbb{H\mkern1mu P}}
\newcommand{\CaP}{\mathrm{Ca}\mathbb{\mkern1mu P}^2}

\newcommand{\C}{{\mathbb{C}}}
\newcommand{\R}{{\mathbb{R}}}
\newcommand{\Z}{{\mathbb{Z}}}
\newcommand{\QH}{{\mathbb{H}}}

\renewcommand{\H}{\ensuremath{\operatorname{\mathsf{H}}}}
\newcommand{\E}{\ensuremath{\operatorname{\mathsf{E}}}}
\newcommand{\F}{\ensuremath{\operatorname{\mathsf{F}}}}
\newcommand{\G}{\ensuremath{\operatorname{\mathsf{G}}}}
\newcommand{\D}{\ensuremath{\operatorname{\mathsf{D}}}}
\newcommand{\SO}{\ensuremath{\operatorname{\mathsf{SO}}}}
\renewcommand{\O}{\ensuremath{\operatorname{\mathsf{O}}}}
\newcommand{\Sp}{\ensuremath{\operatorname{\mathsf{Sp}}}}
\newcommand{\U}{\ensuremath{\operatorname{\mathsf{U}}}}
\newcommand{\SU}{\ensuremath{\operatorname{\mathsf{SU}}}}
\newcommand{\Spin}{\ensuremath{\operatorname{\mathsf{Spin}}}}
\newcommand{\Pin}{\ensuremath{\operatorname{\mathsf{Pin}}}}

\newcommand{\T}{\ensuremath{\operatorname{\mathsf{T}}}}
\renewcommand{\S}{\ensuremath{\operatorname{\mathsf{S}}}}

\newcommand{\fg}{{\mathfrak{g}}}

\newcommand{\fm}{{\mathfrak{m}}}
\newcommand{\fn}{{\mathfrak{n}}}

\newcommand{\fso}{{\mathfrak{so}}}
\newcommand{\fsu}{{\mathfrak{su}}}

\def\con#1=#2(#3){#1 \equiv #2 \bmod{#3}}

\newcommand{\ml}{\langle}                     
\newcommand{\mr}{\rangle}                    

\def\jjoinrel{\mathrel{\mkern-4mu}}
\def\llongrightarrow{\relbar\jjoinrel\longrightarrow}

\def\Lllongrightarrow{\relbar\jjoinrel\relbar\jjoinrel\llongrightarrow}

\newcommand{\diag}{\ensuremath{\operatorname{diag}}}

\newcommand{\rank}{\ensuremath{\operatorname{rk}}}
\newcommand{\corank}{\ensuremath{\operatorname{corank}}}
\renewcommand{\Im}{\ensuremath{\operatorname{Im}}}
\newcommand{\Ad}{\ensuremath{\operatorname{Ad}}}

\renewcommand{\sec}{\ensuremath{\operatorname{sec}}}

\DeclareMathOperator{\Fix}{Fix}

\DeclareMathOperator{\pr}{pr}

\DeclareMathOperator{\id}{id}
\DeclareMathOperator{\Id}{Id}

\newcommand{\no}{\noindent}
\newcommand{\co}{{cohomogeneity }}
\newcommand{\com}{{cohomogeneity one manifold }}

\newcommand{\Kpmo}{\mathsf{K}_{\scriptscriptstyle{0}}^{\scriptscriptstyle{\pm}}}
\newcommand{\Kpo}{\mathsf{K}_{\scriptscriptstyle{0}}^{\scriptscriptstyle{+}}}
\newcommand{\Kmo}{\mathsf{K}_{\scriptscriptstyle{0}}^{\scriptscriptstyle{-}}}
\newcommand{\Kpm}{\mathsf{K}^{\scriptscriptstyle{\pm}}}
\newcommand{\Kp}{\mathsf{K}^{\scriptscriptstyle{+}}}
\newcommand{\Km}{\mathsf{K}^{\scriptscriptstyle{-}}}
\newcommand{\Ko}{\mathsf{K}_{\scriptscriptstyle{0}}}

\newcommand{\Ho}{\mathsf{H}_{\scriptscriptstyle{0}}}

\newcommand{\Kbpm}{\bar{\mathsf{K}}^{\scriptscriptstyle{\pm}}}
\newcommand{\Kbp}{\bar{\mathsf{K}}^{\scriptscriptstyle{+}}}
\newcommand{\Kbm}{\bar{\mathsf{K}}^{\scriptscriptstyle{-}}}
\newcommand{\Ktm}{\tilde{\mathsf{K}}^{\scriptscriptstyle{-}}}
\newcommand{\Ktp}{\tilde{\mathsf{K}}^{\scriptscriptstyle{+}}}
\newcommand{\Hb}{\bar{\mathsf{H}}}
\newcommand{\Ht}{\tilde{\mathsf{H}}}

\newcommand{\Bp}{B^{\scriptscriptstyle{+}}}

\newcommand{\lpm}{l^{\scriptscriptstyle{\pm}}}

\newcommand{\subo}{_{\scriptscriptstyle{0}}}

\newcommand{\Symp}{\mathsf{Sp}}
\newcommand{\Sph}{\mathbb{S}}
\newcommand{\Disc}{\mathbb{D}}
\newcommand{\W}{\mathsf{W}}
\newcommand{\Q}{\mathsf{Q}}
\newcommand{\gS}{\mathsf{S}}
\newcommand{\gC}{\mathsf{C}}
\newcommand{\gH}{\mathsf{H}}\newcommand{\bgH}{\bar{\mathsf{H}}}
\newcommand{\N}{\mathsf{N}}
\newcommand{\K}{\mathsf{K}}
\renewcommand{\L}{\mathsf{L}}

\renewcommand{\F}{\mathsf{F}}

\newcommand{\gT}{\mathsf{T}}

\newcommand{\gW}{\mathsf{W}}

\newcommand{\gK}{\mathsf{K}}

\begin{document}

\title{Positively Curved Cohomogeneity One Manifolds and 3-Sasakian Geometry}

\dedicatory{Dedicated to Wilhelm Klingenberg on his 80th birthday}

\author{Karsten Grove}
\address{University of Maryland\\
        College Park , MD 20742}
\email{kng@math.umd.edu}

\author{Burkhard Wilking}
\address{University of M\"{u}nster\\
        Einsteinstrasse 62\\
        48149 M\"{u}nster, Germany}
\email{wilking@math.uni-muenster.de}

\author{Wolfgang Ziller}
\address{University of Pennsylvania\\
        Philadelphia, PA 19104}
\email{wziller@math.upenn.edu}

\thanks{The first named author was supported in part by the Danish Research
Council and the last named author by the Francis J. Carey Term
Chair. All three authors were  supported by grants from the National
Science Foundation.}

\maketitle



Since the round sphere of constant positive (sectional) curvature is
the simplest and most symmetric topologically non-trivial Riemannian
manifold, it is only natural that manifolds with positive curvature
always will have a special appeal, and play an important role in
Riemannian geometry. Yet, the general knowledge and understanding of
these objects is still rather limited. In particular, although only
a few obstructions are known, examples are notoriously hard to come
by.

The additional structure provided by the presence of a large
isometry group has had a significant impact on the subject (for a
survey see \cite{grove:survey}). Aside from classification and
structure theorems in this context (as in \cite{hsiang-kleiner},
\cite{grove-searle:rank}, \cite{grove-searle:rep}, \cite{grove-kim},
\cite{Wi: torus}, \cite{Wi: sym} and \cite{rong:s5},
\cite{fang-rong:sym}, \cite{fang-rong:dim8and9}), such
investigations also provide a natural framework for a systematic
search for new examples. In retrospect, the classification of simply
connected homogeneous manifolds of positive curvature
(\cite{Be},\cite{Wa},\cite{AW},\cite{Bb}) is a prime example. It is
noteworthy, that in dimensions above 24, only the rank one symmetric
spaces, i.e., spheres and projective spaces appear in this
classification. The only further known examples of positively curved
manifolds are all  biquotients \cite{Es1,Es2,Ba},
 and so far occur only in dimension 13 and below.

A natural measure for the size of a symmetry group is provided by
the so-called cohomogeneity, i.e. the dimension of its orbit space.
It was recently shown in \cite{Wi: sym}, that the lack of positively
curved homogeneous manifolds  in higher dimensions in the following
sense carries over to any cohomogeneity: If a simply connected
positively curved manifold with cohomogeneity $k\ge 1$ has dimension
at least $18(k+1)^2$, then it is  homotopy equivalent to a rank one
symmetric space.

\smallskip

This paper deals with manifolds of cohomogeneity one.  Recall that
in \cite{GZ} a wealth of new nonnegatively curved examples were
found among such manifolds. Our ultimate goal is to classify
positively curved (simply connected) cohomogeneity one manifolds.
The spheres and projective spaces admit an abundance of such actions
(cf.  \cite{HL,straume,iwata:1,iwata:2}, and \cite{uchida}). In
\cite{searle}, however, it was shown that in dimensions at most six,
these are in fact the only ones. In \cite{PV2} it was shown that
this is also true in dimension 7, as long as the symmetry group is
not locally isomorphic to $\S^3\times \S^3$. Recently Verdiani
completed the classification in even dimensions (see
\cite{PV1,verdiani:1, verdiani:2}) :

\begin{main*}[Verdiani]\label{verdiani}
An even dimensional simply connected cohomogeneity one manifold with
an invariant metric of positive sectional curvature is equivariantly
diffeomorphic to a compact rank one symmetric space with a linear
action.
\end{main*}

The same conclusion is false in odd dimensions. There are  three
normal homogeneous manifolds of positive curvature which admit
cohomogeneity one actions: The Berger space $B^7 = \SO(5)/\SO(3)$
with a subaction by $\SO(4)$. The Aloff Wallach manifold $W^7 =
\SU(3)/\diag(z,z,\bar{z}^{2})$ $ = \SU(3)\SO(3)/\U(2)$ with
subactions by $\SU(2)\SO(3)$, denoted by $W^7_{(1)}$, and  by
$\SO(3)\SO(3)$, denoted by $W^7_{(2)}$. And finally the Berger space
$B^{13} = \SU(5)/\Sp(2)\S^1$ with a subaction by $\SU(4)$. It is
perhaps somewhat surprising that none of the remaining homogeneous
manifolds of positive curvature admit \co one actions. More
interestingly, the subfamily $E^7_p = \diag(z,z,z^p) \backslash
\SU(3)/\diag(1,1,\bar{z}^{p+2}) , $ $p\ge 1 $ of inhomogeneous
positively curved Eschenburg spaces  admit \co one actions by
$\SO(3)\SU(2)$ which extends to $\SO(3)\U(2)$. Similarly, the
subfamily of the inhomogeneous positively curved Bazaikin spaces,
$B^{13}_p = \diag(z,z,z,z,z^{2p-1})$ $\backslash\SU(5)/\Sp(2) $ $
\diag(1,1,1,1,\bar{z}^{2p+3}), p \ge 1$
     admit \co one actions by $\SU(4)$, which extend to $\U(4)$. We
point out  that $ E^7_1 =W^7_{(1)}$ with one of
     its \co one actions, and $ B^{13}_1=B^{13} $.

\bigskip

The goal of this paper is to  give an exhaustive description of all
simply connected cohomogeneity one manifolds that can possibly
support an invariant metric with positive curvature.
In addition to the examples already mentioned, it turns out that
only one isolated 7-manifold, $R$ and two  infinite 7-dimensional
families $P_k$ and $Q_k$  potentially admit invariant cohomogeneity
one metrics of positive curvature.

We will also exhibit an intriguing connection between  these
\emph{new candidates for positive curvature} and the cohomogeneity
one self dual Einstein \emph{orbifold metrics} on $\Sph^4$
constructed by Hitchin \cite{Hi2}. As a biproduct,  the manifolds
$P_k$ and $Q_k$ all support  3-Sasakian Riemannian metrics, i.e.,
their Euclidean cones are Hyper-K\"{a}hler (see \cite{BG} for a
survey), and are in particular Einstein manifolds with positive
scalar curvature. In dimension 7 the known examples, due to Boyer,
Galicky, Mann and Rees \cite{BGM,BGMR},  are constructed as
so-called reductions from the 3-Sasakian metric on a round sphere,
and except for $\Sph^7$, have positive second Betti number. They
include the Eschenburg spaces $E_p$ as a special case. The new
3-Sasakian manifolds $P_k$ are particularly interesting since they
are, apart from $\Sph^7=P_1$, the first 2-connected examples  in
dimension 7 (see Theorem C). Both $P_k$ and $Q_k$ are also the first
seven dimensional non-toric 3-Sasakian manifolds, i.e. do not
contain a 3-torus in their isometry group.

\bigskip

To describe the new candidates for positive curvature, recall
that any simply connected cohomogeneity one $\G$-manifold admits a
decomposition $M = \G \times_{\Km} \Disc_- \cup \G \times_ {\Kp}
\Disc_+$ where $\H \subset \{\Km, \Kp\} \subset \G$ are (isotropy)
subgroups of $\G$, and $\Disc_{\pm}$ are Euclidean discs with
$\partial \Disc_{\pm} = \Sph_{\pm} = \Kpm/\H $. Conversely, any
collection of groups $\H \subset \{\Km, \Kp\} \subset \G$ where
$\Kpm/\H$ are spheres, give rise in this fashion to a \co one
manifold.

Using this  notation, we first describe a sequence of 7-dimensional
manifolds $H_k$. They are  given by the groups
$\Z_2\oplus\Z_2\subset\{\Kmo\cdot\H,\Kpo\cdot\H\}\subset_k
\SO(3)\SO(3)$. Furthermore,  the identity
    components $\Kpmo \cong \SO(2)$   depend on
integers $(p,q)$ which
    describe the slope of their embedding into a maximal torus of
$\SO(3)\SO(3)$. They
    are
$(1,1)$ for $\Kmo$ embedded into the lower $2\times 2$ block of
$\SO(3)$, and $(k,k+2)$ for $\Kpo$ embedded into the upper $2\times
2$ block.

    The universal covers of $H_k$ break up into two families,
    $P_k$ being the universal cover of $H_{2k-1}$ with $\G=\SO(4)$ and
    principal isotropy group
$\Z_2\oplus\Z_2$, and $Q_k$ the
    universal cover of $H_{2k}$ with $\G=\SO(3)\SO(3)$ and
      principal isotropy group
$\Z_2$. The additional manifold $R$ is like $Q_k$ but with slopes
$(3,1)$ on the left and $(1,2)$ on the right.

\smallskip

     Our main result can now be formulated as:

\begin{main}

Any odd dimensional simply connected cohomogeneity one manifold $M$
with an invariant metric of positive sectional curvature is equivariantly
diffeomorphic to one of the following:

\begin{itemize}

\item
A Sphere with a linear action,

\item
     One of $E^7_p,  B^{13}_p$ or
$B^7$,

\item

One of the  $7$-manifolds $P_k, Q_k$, or $R$,

\end{itemize}

\no with one of the actions described above.

\end{main}

The first in each sequence $P_k,Q_k$ admit an invariant metric with
positive curvature since $P_1=\Sph^7$ and $Q_1=W^7_{(2)}$. For more
information and further discussion of the non-linear examples we
refer to Section \ref{sec: examples}.

There are numerous 7 dimensional cohomogeneity one manifolds with
singular orbits of codimension two, all of which by {\cite{GZ}} have
invariant metrics with non-negative curvature. Among these, there
are two subfamilies like the above $P_k$ and $Q_k$, but where the
slopes for $\Kpm$  are arbitrary. It is striking that in positive
curvature, with one exception, only the above slopes are allowed.
The exception is given by the positively curved \co one action on
$B^7$, where the isotropy groups are like those for $P_k$ with
slopes $(1,3)$ and $(3,1)$. In some tantalizing sense then, the
exceptional Berger manifold $B^7$ is associated with the $P_k$
family in an analogues way as the exceptional candidate $R$ is
associated with the $Q_k$ family. It is also surprising that all
non-linear actions in Theorem A, apart from the Bazaikin spaces
$B_p^{13}$, are \co one under a group locally isomorphic to
$\S^3\times\S^3$.

\vspace{10pt}

     As already indicated, the  manifolds $H_k$  have another intriguing characterization.
To describe this in more detail, recall that $\Sph^4$ and $\CP^2$ according to
Hitchin  are the only smooth self dual Einstein 4-manifolds.
However, in the more general context of orbifolds, Hitchin
constructed \cite{Hi1}  a sequence of self dual Einstein orbifolds
$O_k$ homeomorphic to $\Sph^4$, one for each integer $k>0$, which
are invariant under a cohomogeneity one $\SO(3)$ action. It has an
orbifold singularity whose angle normal to a smooth $\SO(3)$ orbit
$\RP^2$ is equal to $2\pi/k$. Here $O_{2k}$ can also be interpreted
as an orbifold metric on $\CP^2$ with normal angle $2\pi/k$, and the
cases of  $k=1,2$  correspond to the smooth standard metrics
 on $\Sph^4$ and on
$\CP^2$ respectively.  In general, any self dual Einstein orbifold
gives rise to a 3-Sasakian orbifold metric on the Konishi bundle,
which is the $\SO(3)$ orbifold principal bundle of the vector bundle
of self dual 2-forms. The action of $\SO(3)$ on the base lifts to
form a \co one $\SO(3)\SO(3)$ action on the total space, and we will
prove the following surprising relationship with our positive
curvature candidates:

\begin{main}\label{Hitchin}
For each $k$, the total space  of the  Konishi bundle corresponding
to the selfdual Hitchin orbifold $O_k$ is a smooth 3-Sasakian
manifold, which is
    equivariantly diffeomorphic to  $H_{k}$
    with its  cohomogeneity one $\SO(3)\SO(3)$ action.
\end{main}

In this context we note that the exceptional manifolds $B^7$ and $R$
can be described, up to covers, as the $\SO(3)$ orbifold principal
bundles of the vector bundle of anti-self dual 2-forms over $O_{3}$
and $O_4$ respectively.

\smallskip

It was shown by O.Dearricott in \cite{De1} that Konishi metrics,
scaled down in direction of the principal $\SO(3)$ orbits, have
positive sectional curvature if and only if the self dual Einstein
orbifold base has positive curvature. Unfortunately, the Hitchin
orbifold metrics  do not have positive curvature for $k > 2$, so
this appealing description does not easily yield the desired metrics
of positive curvature on $P_k$ and $Q_k$.

\vspace{10pt}

Our candidates also have interesting
topological properties:

\begin{main}\label{topology}
The manifolds $P_k$ are two-connected with $\pi_3(P_k) =\Z_k$. For
the manifolds $Q_k$ and $R$ we have $H^2(Q_k,\Z)=H^2(R,\Z)=\Z$ and
$H^4(Q_k,\Z)=\Z_{2k+1}$, respectively  $H^4(R,\Z)=\Z_{35}$.
\end{main}

We note that the cohomology rings of $Q_k$ and $R$ occur as the
cohomology rings of one or more of the  seven dimensional positively
curved  Eschenburg biquotients \cite{Es1},\cite{Es2}. In fact,
surprisingly, $Q_k$ has the same cohomology ring as $E_k$. On the
other hand the manifolds $P_k$ have the same cohomology ring as
$\Sph^3$ bundles over $\Sph^4$, and among such manifolds, so far
only $\Sph^7$ and the Berger space $B^7$  (see  \cite{GKS}) are
known to admit metrics of positive curvature. It would be
interesting to know whether there are other cases where a manifold
in the families $P_k,Q_k$  is diffeomorphic to an Eschenburg space
or to an $\Sph^3$ bundles over $\Sph^4$.

The fact that the manifolds $P_k$ are 2-connected is particularly
significant. Recall that by the finiteness theorem of
Petrunin-Tuschmann \cite{PT} and Fang-Rong
\cite{fang-rong:finiteness}, 2-connected manifolds play a special
role in positive curvature since there exist only finitely many
diffeomorphism types of such manifolds, if one specifies the
dimension and the pinching constant, i.e. $\delta$ with $\delta \le
sec \le 1$. Thus, if $P_k$ admit positive curvature metrics, the
pinching constants $\delta_k$ necessarily go to 0 as $k\to \infty$,
and $P_k$ would be the first examples of this type. The existence of
such metrics would provide counter examples to a conjecture by Fang
and Rong in \cite{fang-rong:sym} (cf. also Fukaya \cite{fukaya},
Problem 15.20).

\bigskip

We conclude the introduction by giving a brief discussion of the
proof of our main result and  how we have organized it.

The most  basic \emph{recognition tool} one has is of course the
group diagram itself. However, given just the richness of linear
actions on spheres, one would expect that looking primarily for such
detailed information might actually hinder classification. It is
thus crucial to have other recognitions tools at our disposal, that
do not need the full knowledge of a group diagram. In fact, in our
proof we often either exclude a potential manifold, or determine
what it is before we actually derive a possible group diagram.

    For this we first note that Straume \cite{straume} has provided a
complete classification of all cohomogeneity one actions on homotopy
spheres. Aside from linear actions on the standard sphere, there are
families of non-linear actions, and also actions on exotic Kervaire
spheres. It was observed by Back and Hsiang \cite{BH} (Searle
\cite{searle}  in dimension 5) that only the linear ones support
invariant metrics of positive curvature (in dimensions other than
five they cannot even support invariant metrics of nonnegative
curvature \cite{kervaire}). In particular, for our purposes it
suffices to recognize the underlying manifold as a homotopy sphere,
and we have
    two specific tools for doing so: One of them
is provided by the (equivariant diffeomorphism) classification of
positively curved \emph{fixed point homogeneous} manifolds
\cite{grove-searle:rep} , i.e., manifolds on which a group $\G$ acts
transitively on the normal sphere to a component of its fixed point
set $M^{\G}$. The other is the  \emph{Chain Theorem} of \cite{Wi:
sym}, which classifies 1-connected positively curved manifolds up to
homotopy that support an isometric action by one of the classical
groups, $\SO(n), \SU(n)$ or $\Sp(n)$ so that its principal isotropy
group contains the same type of group as a standard $3\times3$ block
(or $2\times2$ block in case of $\Sp(n)$).

\smallskip

  Our classification of positively curved manifolds with
an isometric  cohomogeneity one $\G$-action is  done by induction on
the dimension of the  manifold $M$. Here the induction step is
typically done via \emph{reductions}, i.e., by analyzing fixed point
sets of subgroups of $\G$ and how they sit inside of $M$. Since such
fixed point sets are totally geodesic, they are themselves
positively curved manifolds of cohomogeneity at most one and hence
in essence known by assumption. In this analysis, the basic
\emph{connectivity lemma} of \cite{Wi: torus} which asserts that the
inclusion map of a totally geodesic codimension $k$ submanifold in
an $n$ dimensional positively curved manifold is $n-2k+1$ connected,
naturally plays an important role.

Another variable in the proof is $\rank \G$, the rank of $\G$. Here
it is a simple but important fact that in positive curvature, the
corank of the principal isotropy group $\H$, i.e., $\corank \H =
\rank \G - \rank \H$ is 1 in even dimensions, and 0 or 2 in odd
dimensions. The equal rank case is fairly simple and induction is
not used here (see Section \ref{sec: equalrank}).

\smallskip

The following brief description of the content of the sections will
hopefully support the overall understanding of the strategy of the
proof just outlined.

In Section \ref{sec: tech} we recall some essential simple curvature
free facts about cohomogeneity one manifolds we will need
throughout. This includes a discussion of the \emph{Weyl groups} and
\emph{reductions}, i.e. fixed point sets of subgroups, including the
\emph{core} of the action.

Sections \ref{sec: obstruction} and \ref{sec: weyl} form the
geometric heart of the paper. It is here we present and derive all
our \emph{obstructions} stemming from having an invariant metric of
positive curvature. Some of these, which have been derived earlier
in more general settings (see \cite{Wi: torus, Wi: sym}), become
particularly powerful in the context of cohomogeneity one manifolds.
Other than the rank restriction,  which enters from the outset, two
key obstructions used throughout are primitivity, and restrictions
imposed on the isotropy representation of the principal isotropy
group. The full strength of primitivity is derived in Section
\ref{sec: weyl} after a classification of all  Weyl groups corresponding to non trivial cores. It is also shown here that all  Weyl groups are finite and strong bounds on their orders are derived.

 Section \ref{sec: examples} we present and discuss some of the
properties of the \co one actions on the known examples of positive
curvature, as well as on the new candidates.

We start the classification in Section \ref{sec: equalrank} with the
equal rank case and in Section \ref{sec: not semi} we deal with  the
case where $\G$ is not semisimple. For semisimple groups $\G$, it
turns out to be useful to prove the theorem for groups of rank 2 or
3 first, and this is done in the Section \ref{sec: rank2} and
\ref{sec: rank3}. In a sense these two sections form the core of the
classification. It is here all non spherical examples emerge. The
case of semisimple groups $\G$ with $\rank \G \ge 4$ is done
separately for non-simple groups in Sections \ref{sec: rank one
factor} and \ref{sec: not simple} and for simple groups in Section
\ref{sec: simple}.

In Section  \ref{sec: Hitchin} we exhibit our new infinite families
of candidates as 3-Sasakian manifolds (Theorem B), and  in Section
\ref{sec: topology} we prove Theorem C. These sections can be read
independently of the rest of the paper.

Since we need the classification in even dimensions, we have added a
relatively short proof as a service to the reader in Appendix I. As
another service to the reader, we have collected the cohomogeneity
one diagrams for the essential  actions on rank one symmetric
spaces, and other known useful classification results  in Appendix
II.


\bigskip

\section{Cohomogeneity one manifolds.}\label{sec: tech}\label{sec: general}

We begin  by discussing a few useful general facts about closed
cohomogeneity one Riemannian $\G$ manifolds $M$ and fix notation we
will use throughout. Readers with good working knowledge of
cohomogeneity one manifolds may want to proceed
 to Section \ref{sec: obstruction}, Section \ref{sec: weyl} and the
classification starting in Section \ref{sec: equalrank} immediately
and refer back to this section whenever needed.

 Our primary
interest is in positively curved, 1- connected $\G$ manifolds $M$
with $\G$ connected. However, since fixed point sets with induced
\co one actions play a significant role in our proof, it is
important to understand the more general case where $\G$ is not
connected, and $M$ is connected with possibly non-trivial finite
fundamental group.
\smallskip

Since $M$ has finite fundamental group, the orbit space  $M/\G$ is an interval
and not circle.
 The end points of the interval correspond
to two non-principal orbits, and all interior points to principal
orbits. By scaling the metric if necessary we may assume that  $M/\G
= [-1,1]$ as a metric space.

\smallskip

Fix a normal geodesic $c: \R \to M$ perpendicular to all orbits (an
infinite horizontal lift of  $M/\G$). The image $C = c(\R)$ is
either an embedded circle, or a 1-1 immersed line (cf.
\cite[Proposition 3.2]{AA}). We denote by $\H$ the principal isotropy
group $\G_{c(0)}$ at $c(0)$, which is equal to the isotropy groups
$\G_{c(t)}$ for all $t\ne 1 \mod 2\Z$, and by $\Kpm$ the isotropy
groups at $p_\pm=c(\pm 1)$. Then $M$ is the union of tubular
neighborhoods of the non-principal orbits $B_{\pm} = \G/\Kpm$ glued
along their common boundary $\G/\H$, i.e., by the slice theorem
\begin{equation}\label{double disc}
M = \G \times_{\Km} \Disc_- \cup  \G \times_{\Kp} \Disc_+,
\end{equation}

\no where $\Disc_{\pm}$ denotes the normal disc to the orbit $\G
p_\pm = B_{\pm} $ at $p_\pm$. Furthermore, $\Kpm/\H=\partial
\Disc_{\pm}= \Sph_{\pm}$ are spheres, whose dimension we denote by
$l_\pm$. It is important to note that  the  diagram of groups

\begin{equation}\label{diagram}
\begin{split}
        \xymatrix{
           & {\G} & \\
           {\Km}\ar@{->}[ur]^{j_-} & & {\Kp}\ar@{->}[ul]_{j_+} \\
           & {\H}\ar@{->}[ul]^{h_-}\ar@{->}[ur]_{h_+} &
           }
\end{split}
\end{equation}
where $j_\pm$ and $h_\pm$ are the natural inclusions,
 which we also record as
\begin{equation}\label{linediag}
\H \subset \{\Km,\Kp\} \subset \G,
\end{equation}

\no determines $M$. Conversely, such a \emph{group diagram} with
$\Kpm/\H = \Sph^{l_{\pm}}$, defines a cohomogeneity one $\G$-manifold.

In section \ref{sec: Hitchin}, we will see that the above
construction, as well as the principal bundle construction for
cohomogeneity one manifolds in [GZ], naturally carries over to a
large class within the more general context of orbifolds.

We point out that the spheres $\Kpm/\H$ are often highly ineffective
and we denote by $\H_\pm$ their ineffective kernel. It will be
convenient to allow the ineffective kernel of  $\G/\H$ to be finite,
i.e., to allow the action to be \emph{almost effective}.

A non-principal orbit $\G/\K$ is called \emph{exceptional} if $\dim
\G/\K = \dim \G/\H$ or equivalently $\K/\H = \Sph^0$. Otherwise
$\G/\K$ is called \emph{singular}. As usual we refer to the collection
$M_{\subo}$ of principal orbits, i.e., $M - (B_- \cup B_+)$ as the
\emph{regular part} of $M$.

\bigskip

\begin{center}
The Cohomogeneity One Weyl Group.
\end{center}

\smallskip

The \emph{Weyl group}, $\W(\G,M) = \W$ of the action, is by
definition the stabilizer of the geodesic $C$ modulo its kernel
$\H$. If $\N(\H)$ is the normalizer of $\H$ in $\G$, it is easy to
see (cf. \cite{AA}) that $\W$ is a dihedral subgroup of $\N(\H)/\H$,
generated by unique involutions $w_{\pm} \in (\N(\H) \cap \Kpm)/\H$,
and that $M/\G = C/\W$. Each of these involutions can also be
described as the unique element $a\in\Kpm$ mod $\H$ such that $a^2$
but not $a$ lies in $\H$.

Note that $\W$ is finite if and only if $C$ is a closed geodesic,
and in that case the order $|\W|$ is the number of minimal geodesic
segments $C - (B_- \cup B_+)$. Note also that any non principal isotropy
group along $c$ is of the form $w \Kpm w$ for some $w \in \N(\H)$
representing an element of $\W$. The isotropy types $\Kpm$ alternate
along $C$ and hence half of them are isomorphic to $\Kp$ and half to
$\Km$, in the case where $\W$ is finite.

\bigskip

\begin{center}
Group Components.
\end{center}

\smallskip

In this section $\G$ is a not necessarily connected Lie group acting
with cohomogeneity one on a connected manifold $M$ with finite fundamental
group. From the description of $M$ as a double disc
bundle~\eqref{double disc},
 we see that
\begin{eqnarray}\label{transversality}
\text{ \quad $\G/\Kpm \cong B_{\pm} \to M$}&&\text{ is \quad $l_\mp$-connected.}\\ \nonumber
\text{ \quad $\G/\gH  \to M$}&&\text{ is \quad
 $\min\{l_-,l_+\}$-connected.}
\end{eqnarray}
Recall that by definition a map $f: X \to Y$ is \emph{$l$-connected} if
the induced map $f_i: \pi_i(X) \to \pi_i(Y)$ between homotopy groups
is an isomorphism for $i < l$ and surjective for $i = l$.

First observe that it is impossible that both $l_\pm=0$.
 Indeed, if both normal bundles to $\G/\Kpm$ are trivial $M$ is a bundle over $\Sph^1$.
If one of the orbits
say $\G/\Kp$ has non-trivial normal bundle
the two fold cover $\G/\H\to\G/\Kp$ gives rise to a two fold cover
$M'$ of $M$ on which $\G$ acts by cohomogeneity one with diagram $\H
\subset \{\Km,w_+ \Km w_+\} \subset \G$. We are now either in the
first situation, or we can repeat the second argument indefinitely,
contradicting that $\pi_1(M)$ is finite.

If both  $l_\pm >0$ , \eqref{transversality} implies that $\G/\H$ is
connected and hence $\G$ and $\G_{\subo}$ have the same orbits, and in particular the same Weyl group.
 If one of $l_\pm $ say $l_- =0$ and $l_+ >0$,
\eqref{transversality} implies that $\G/\Km$ is connected. Since
$\G/\gH$ is a sphere bundle over $\G/\Km$, it follows that $\G/\gH$
has at most two components. This in turn implies that
\begin{equation}\label{bundlecomp}
\text{The Weyl group of  the $\G_{\subo}$ action  has index at most $2$ in
 the Weyl group of $\G$.}
\end{equation}

We now assume that $M$ is simply
connected and $\G$ is connected. The above covering argument then
implies that there cannot be any exceptional orbits. If both
$l_\pm\ge 2$, \eqref{transversality} implies that all orbits are
simply connected and hence all isotropy groups connected. If one of
$\l_\pm$ say $l_-=1$ and $l_+\ge 2$, then $\G/\Km$ is simply
connected and hence $\Km$ connected. Since $\G/\H$ is a circle
bundle over $\G/\Km$ it follows that $\pi_1(\G/\gH)$ and hence
$\H/\Ho\simeq \Kp/\Kpo$ are cyclic.
 In
summary,

\begin{lem}\label{Hcomp}
Assume that $\G$ acts on $M$ by \co one with $M$  simply connected
and $\G$ connected. Then:

\begin{enumerate}
\item[(a)] There are no exceptional orbits, i.e. $l_\pm \ge 1$.
        \item[(b)] If both $l_\pm \ge 2$, then $\Kpm$ and $\H$ are all
        connected.
\item[(c)] If one of $l_\pm$, say $l_-=1$, and  $l_+\ge 2$, then
$\Km=\H\cdot\S^1=\Ho\cdot\S^1$, $\H=\Ho\cdot\Z_k$ and
$\Kp=\Kpo\cdot\Z_k$.
\end{enumerate}
\end{lem}

The situation where both $l_{\pm} = 1$ is analyzed in the presence of an
invariant positively curved metric in \eqref{both-one}.
Finally we observe

\begin{lem}\label{isocomp}
Suppose $ \Kbpm \subset \Kpm$ are subgroups  with $\Kpm/ \Kbpm$
finite, $\Kbpm\not\subset \gH$, and $ \Kbm \cap \H  =  \Kbp \cap \H =:  \Hb$. Then $\Km /
\Kbm \simeq \H/\Hb \simeq \Kp/  \Kbp$ and the cohomogeneity one
manifold $\bar M$ defined by ${\Hb} \subset \{ \Kbm, \Kbp\} \subset
\G$ is an $\H/ \Hb$ cover of $M$.
\end{lem}

In general, a subcover of a compact \co one manifold  with finite
fundamental group and $\G$ connected, is obtained by a combination
of the following three: We can add components to $\Kpm$ and $\H$ as
in (\ref{isocomp}), or we can divide $\G$ by a central subgroup
which does not intersect $\Kpm$. These two yield orbitspace preserving covering maps. We can also create a subcover where one  of the orbits is
exceptional, if $\Kp$ is the $w$ conjugate of $\Km$
 for an order two element  in $\N(\H)/\H$ represented
by   $w\in \N(\H)$.

\smallskip

\bigskip

\begin{center}
Reductions.
\end{center}

\smallskip

Fixed point sets of subgroups $\L \subset \G$ will play a pivotal
role throughout. It is well known that the fixed point set $M^\L$ of
$\L$ consists of a disjoint union of totally geodesic submanifolds.
If $M^\L$ is non empty, $\L$ is of course a subgroup of an isotropy
group, and hence of $\H$ or of $\Kpm$ (up to conjugacy). In general
when $\L \subset \K \subset \G$, it is well known that $\N(\L)$ acts
with finite orbit space on $(\G/\K)^\L$, and transitively when $\L =
\K$, or when $\L$ is a maximal torus of $\K$ (see e.g.
\cite{Br}, Corollary II.5.7).

\smallskip

Suppose first that  $\L\subset \K_-$ is not conjugate to a  subgroup
of $\H$.
 Then no component
of $M^\L$ intersects the regular part $M_{\subo}$ of $M$. In this case,
all components of $M^\L$ are homogeneous, and we usually consider
the component in one of $B_\pm$ say $B_-$ containing $p_-$ which
equals $\N(\L)_{\subo}/\N(\L)_{\subo}\cap\Km$. As a particular application of
this, we point out that a central involution in $\G$ which lies in
one of $\Kpm$ say $\Km$ but not in $\H$, has $\G/\Km$ as its fixed
point set.

\smallskip

If $\L$ is conjugate to a subgroup of $\H$, the components of
$M^{\L}$ which intersect the regular part of $M$ form a
cohomogeneity one manifold under the action of $\N(\L)$ since
$\N(\L)$ acts with finite quotient on $(\G/\H)^{\L}$.  Each
component of $M^\L$ that intersects the regular part is hence a
cohomogeneity one manifold under the action of the subgroup of
$\N(\L)$ stabilizing the component.   Unless otherwise stated, the
reduction we will consider is the component $M^\L_c$ of $M^\L$
containing the geodesic $c$. We will denote it's stabilizer subgroup
of $\N(\L)$ by $\N(\L)_c$ and  refer to $(M^\L_c,\N(\L)_c)$  as
\emph{reductions}  (for general actions see
\cite{grove-searle:core}). In general the length of
 $M^\L_c/\N(\L)_c$ is an integer multiple
of the length of $M/\G$. The orbit spaces coincide if both $\N(\L)
\cap \Kpm$ act nontrivially on the  normal spheres
 of $M^\L_c\cap B_\pm\subset M^{\L}_c$ at $p_\pm$,
which are given by  $\Sph^{\L}_\pm= \N(\L)_c \cap
\Kpm/\N(\L)_c\cap\H $. If this is the case, $\N(\L)_c$ acts ($\L$
ineffectively) by cohomogeneity one on $M^\L_c$ with orbit space
$M/\G$, and diagram $\N(\L)_c \cap \H \subset \{\N(\L)_c \cap \Km,
\N(\L)_c \cap \Kp\} \subset
\N(\L)_c$.

In the main part of the induction proof, it is usually sufficient to
consider the cohomogeneity one action of the connected component
$\N(\L)_{\subo}$ of $\N(\L)_c$ on $M_c^{\L}$ keeping in mind that its Weyl
group need not be that of $M$.

\smallskip

If  $\L$ is a maximal torus of $\Ho$ and $a\in \N(\H)$, then $a \L a^{-1}\subset \Ho$  is also conjugate to
$\L$ by an element in $\Ho$. In particular, one can represent $w_\pm$ by elements
in the normalizer of $\L$. The same holds by definition of the Weyl
group for $\L=\H$, and hence:

\begin{lem}[Reduction Lemma]\label{reduction}
If $\L$ is either equal to $\H$ or given by a maximal torus of
$\Ho$, then $\N(\L)_c/\L$ acts by cohomogeneity one  on $M^\L_c$
and the corresponding Weyl groups coincide.
\end{lem}

In the most reduced case where $\L = \H$, we refer
to $M^{\H}_c$ as the \emph{core} of $M$ and $\N(\H)_c$ as
the \emph{core group}.

Often we consider also the least reduced case, that is we take the
fixed point set of an involution
or of an element $\iota$ whose square, but not $\iota$
itself, lies in the center of $\G$, i.e. $\iota$ is an involution in
some central quotient of $\G$.
  In this case  we can determine $\N(\langle \iota
\rangle)=\N(\iota)$ using the well known fact that $\G/\N(\iota)$ is
a symmetric space with $\rank(\N(\iota))=\rank(\G)$,
 and appeal to their classification,
see Table \ref{symm}, Appendix II.

\smallskip

In general the codimension of a reduction might be odd.
However, if $\L$ is a subgroup of a torus in $\T \subset \G$,
and $M$ is positively curved and odd dimensional,
 then  all components
of $M^\L$ have even codimension.
One can establish this fact by induction on the dimension, where one
uses that odd dimensional positively curved
manifolds are orientable and that the statement holds
for cyclic subgroups $\L\subset \T$.

As a simple consequence of the Rank
Lemma \ref{rank}, we also see that in positive curvature, $M^{\L}$ has even codimension  when $\rank
\N(\L) = \rank \G$ and $\rank \G - \rank \H = 2$.

\bigskip

\begin{center}
Equivalence of diagrams.
\end{center}

Recall that in order to get a
 group diagram we choose an invariant metric on $M$.
Thus it can happen that different metrics on the manifold give
different group diagrams.
Of course, one can conjugate all three groups by the same element in
$\G$, and one can also switch $\Km$ and $\Kp$.

Let us now fix a point $p$ in the regular part of the manifold and
an orientation of the normal bundle $\G p$. For each invariant
metric $g$ on the manifold we consider the minimal horizontal
geodesic  $c_g\colon [-\eps_1(g),\eps_2(g)] \to M$ from the left
singular orbit to the right with $c_g(0)=p$. We reparametrize these
geodesics relative to a
 fixed parametrization of the orbit space $M/\G = [-1,1]$, where the orbit
 through $p$ corresponds to $0$. The resulting curves $\bar{c}_g$
 are fixed pointwise by $\H$.
 Using a (smooth) family
of such reparametrized geodesics in $M^{\H}$ corresponding to convex
combinations
of two invariant metrics $g_1$, $g_2$ and the fact that $\N(\H)$
acts transitively on $(\G/\H)^{\H}$,  we can find a curve $a\colon
[-1,1]\rightarrow \N(\gH)_{\subo}$ such that the curve $\bar{c}_{g_2}$
is given by $a(t) \bar{c}_{g_1}(t)$. This proves
that we can find two elements $a_-,a_+\in \N(\gH)_{\subo}$ such that the
group diagram from the metric $g_2$ is obtained from the group
diagram for the metric $g_1$ by conjugating $\Kpm$ with $a_{\pm}$.
On the other hand it is easy to see that indeed for any $a_-,a_+ \in
\N(\gH)_{\subo}$ one can find a metric for which there is  a horizontal
geodesic from $a_- c(-\eps_1(g))$ to
$a_+c(\eps_2(g))$. In fact this
can be achieved by  changing the metric on the complement of two
small tubular neighborhoods of $B_{\pm}$.

All in all we conclude that two  group diagrams
 $\H\subset \{\Km,\Kp\}\subset \G$ and $\Ht\subset \{\Ktm,\Ktp\}\subset \G$
 yield the same cohomogeneity one manifold up to equivariant
diffeomorphism if and only if after possibly switching the roles of
$\Km$ and $\Kp$,  the following holds: There is a $b\in \G$ and an
$a\in \N(\gH)_{\subo}$ with $\Km= b \Ktm b^{-1}$, $\H= b \Ht b^{-1}$, and
$\Kp= ab \Ktp b^{-1}a^{-1}$ (cf. also \cite{neumann}).


\section{Positive Curvature Obstructions.}\label{sec:
obstruction}

\bigskip

In this section we will discuss a number of severe obstructions on a
cohomogeneity one manifold to have an invariant metric with positive
curvature. We point out that none of our obstructions are caused by
nonnegative curvature only. We also mention that Alexandrov geometry
of orbit spaces, which is used extensively to obtain our two
geometric recognitions tools \eqref{cofix} and \eqref{chain}, enter
only once directly in our proof, namely the rank two case \eqref{s3s3}.

The simplest obstruction is a direct consequence of the well known
fact that an isometric torus action on a positively curved manifold
has fixed points in even dimensions and orbits of dimension at most
one in odd dimensions. Since spheres $\K/\H$ have corank at most
one, this gives:

\begin{lem}[Rank Lemma]\label{rank}
One of $\Kpm$ has corank $0$, when $M$ is even dimensional,
and at most corank $1$, when $M$ is odd dimensional. In particular
$\H$ has corank $1$ if M is even dimensional, and corank $0$
or $2$ when $M$ is odd dimensional.
\end{lem}

\smallskip

A second powerful and much more difficult result expresses in two
ways how the \emph{representation of the triple} $\H \subset \{\Km,
\Kp \}$ in $\G$ is \emph{maximal}. The first of these, which we will
refer to as \emph{linear primitivity}, follows  from \cite[Corollary
10]{wilking:dual}, and has the Weyl group bound below as an
immediate consequence. As we will see in the next section this type
of primitivity implies that the Weyl group is finite  as well (see
\eqref{weyl}).

     To define the second kind of primitivity,
      we say that a $\G$ - manifold is \emph{non-primitive} if there is a
$\G$ equivariant map $M \to \G/\L$ for some subgroup $\L \subset \G$
(see \cite[p.17]{AA}). Otherwise, the action is said to be
\emph{primitive}. For cohomogeneity one manifolds, non-primitivity
is equivalent to the statement that for some representation we have
$\H \subset \{\Km, \Kp\} \subset \L \subset \G$, i.e., for some
invariant metric and some normal geodesic, $\Kpm$ generate a proper
subgroup of $\G$. In terms of the original groups, the action is
hence primitive if $\Km$ and $n \Kp n^{-1}$ generate $\G$, for any
fixed  $ n \in\N(\H)_c$.

In the next section we will show that the core with its core action
is primitive \eqref{weyl-core}. When this is combined with linear
primitivity for $\G$, we will show that the $\G$ action itself is
primitive (see \eqref{prim}):

\begin{lem}[Primitivity Lemma]\label{primitivity} Let $c\colon \R \to M$ be
any horizontal geodesic as above. Then

\begin{enumerate}

\item [(a)]

     \emph{(Linear Primitivity)}\label{lin} The Lie algebras of the isotropy
groups along $c$
     generate $\fg$
as a vectorspace .

\item [(b)]

     \emph{(Lower Weyl Group
Bound)}\label{bound} The Weyl group is finite, and $ |\W|\ge
2\dim(\G/\H)/(l_-+l_+)$.

\item [(c)]

\emph{(Group Primitivity)}\label{pim} Any of the singular isotropy
groups $\Kpm$, together with any conjugate of the other by an
element of the core group, generate $\G$ as a group. In particular
this is true for conjugation by elements of  $\N(\H)_{\subo}$.
\end{enumerate}

\end{lem}

\bigskip

The following obstructions deal with \emph{isotropy
representations}. The first of these is a special case of a more
general result in \cite{Wi: sym}, although in our situation it also
follows from linear primitivity.  The second part of the lemma
follows from the first part and the classification of transitive
actions on spheres, see Table \ref{trans}, Appendix II.

\begin{lem}[Isotropy Lemma]\label{isotropy}
Suppose $\H$ is non trivial. Then

\begin{enumerate}

\item[(a)]
Any  irreducible subrepresentation of the isotropy representation of
$\G/\H$ is equivalent to a subrepresentation of the isotropy
representation of one of $\K/\H$, where $\gK$ is an isotropy group
of some point in $c(\R) - M_{\subo}$.

\item[(b)]
The isotropy representation of $\G/\H_{\subo}$ is spherical, i.e. $\H_{\subo}$ acts
transitively on the unit sphere of any $k$ dimensional irreducible
subrepresentation if $k>1$.

\end{enumerate}
\end{lem}

Notice that part a) implies that any subrepresentation of $\G/\H$,
i.e. the isotropy representation of $\G/\H$, is weakly equivalent to
a subrepresentation of $\Km/\H$ or $\Kp/\H$. Recall that two
representations of $\H$ are weakly equivalent if they are equivalent
modulo an automorphism of $\H$. We thus often say that a particular
representation has to degenerate in $\Kp/\H$ or $\Km/\H$.

The fact that the isotropy representations are spherical is a
particularly powerful tool. In \cite{Wi: sym} one finds an
exhaustive list of such so-called \emph{spherical subgroups} when
$\H$ and $\G$ are simple (apart from the case where $\H$ is a rank
one group in an exceptional Lie group). We reproduce this list in
Table \ref{sph}, since it will be used frequently.

\smallskip

\no Lemma \ref{isotropy} has the following very useful consequence:

\begin{lem}\label{Hrank2}
If $\G$ is simple, $\H$ can have at most one simple normal subgroup
of rank at least two.
\end{lem}

\begin{proof} Assume that  $\L_1$ and $\L_2$ are two simple normal
subgroups of $\H$ with $\rank \L_i\ge 2$. From the classification of
transitive actions on spheres it follows that either $\L_1$ or
$\L_2$ must act trivially on the irreducible subrepresentations of
$\H$ in $\Kpm$. By the Isotropy Lemma the same then holds for each
irreducible subrepresentation of $\H$ in $\G$.

We  decompose $\fg=\fm_1\oplus\fm_2\oplus\fn $ where
  $\L_1$ acts non-trivially on $\fm_1$ and
trivially on $\fm_2$, $\L_2$ acts trivially on $\fm_1$ and
non-trivially on $\fm_2$ and both act trivially on $\fn$.
Note that $[\fm_1,\fm_2]=0$ since both $\L_1$ and $\L_2$ act
non-trivially on any subrepresentation of $\fm_1\otimes \fm_2$.
Similarly $[\fm_1,\fn]\subset \fm_1$ and in summary
$[\fm_1,\fg]\subset \fm_1+[\fm_1,\fm_1]$. Using the Jacobi identity
we see that $[[\fm_1,\fm_1],\fn]\subset \fm_1+[\fm_1,\fm_1]$ and
$[[\fm_1,\fm_1],\fm_2]=0$. Thus $\fm_1+ [\fm_1,\fm_1]$ is an ideal
of $\fg$, a contradiction.
\end{proof}

\smallskip

For the singular orbits there are two relevant representations, the
\emph{isotropy representation} and the \emph{slice representation}. These are
related via equivariance of the \emph{second fundamental form}

\begin{equation}
B^{\pm}: S^2(T_{\pm}) \to T_{\pm}^{\perp}
\end{equation}

\no where $T_{\pm}$ is the tangent space of $B_{\pm} \cong \G/\Kpm$ at $p_{\pm}
$, and $T_{\pm}^{\perp}$ is the normal space.

\smallskip

As an example of an application of this, it sometimes follows that
equivariance forces a singular orbit to be totally geodesic. In
particular, this singular orbit must then be in the short list of
positively curved homogeneous manifolds, see Table \ref{trans} and
\ref{homog} in Appendix II.

\smallskip

The next  result also follows from equivariance of the second
fundamental form applied to a  singular orbit.

\begin{lem}[Product Lemma]\label{product}
Suppose $\G = \L_1 \times \L_2$ is semisimple and that the identity component one of $\Kpm$
 is a product subgroup, say  $\Kmo=\K_1\times \K_2$
 and that one of $\N^{\L_i}(\K_i)/\K_i$ is finite. Then $M$
cannot carry a positively curved $\G$-  invariant metric if it is odd
dimensional.
\end{lem}

\begin{proof}
    The condition on the normalizers implies, by Schur's Lemma, that
every invariant  metric on $\G/\K$ is (locally) a product metric on
$(\L_1 /\K_1)\times (\L_2/\K_2) $.
      Denote by $U_i$  the subspace  tangent
      to the factor $\L_i/\K_i$. Note that $\dim U_i > 1$ since $\G$
      is semisimple.

      From the classification of transitive actions on spheres, see
Table \ref{trans},
       we may assume that one of the factors, say $\K_1$,
      acts transitively on the normal sphere. Since $\K_1$ acts trivially on
$U_2$, no subrepresentation of
      $S^2U_2$ is equivalent to the slice representation,
      and hence $B_{S^2U_2}=0$. Since any plane generated
      by one  vector in $U_1$ and one vector in $U_2$
      has intrinsic curvature $0$, we see from the Gauss equation that
      $B(u_1,u_2)\neq 0$ for all nonzero $u_i\in U_i$.
      Because $B$ is bilinear, this implies
      that $\dim(U_i)\le \dim(T^{\perp})$.

      If there  exists a
      $\K_1$-invariant subspace $U_1'\subset U_1$ such that the induced
      representation in $U_1'$ is not equivalent to the slice
      representation, then  the equivariance
      of $B$ implies that $B_{U_1'\otimes U_2}=0$ contradicting
$B(u_1,u_2)\neq 0$ for all nonzero $u_i\in U_i$.
      Thus, using in addition the above dimension restriction,
       the representation of $\K_1$ on all of  $U_1$ is equivalent to
      the slice representation. In particular, $\K_1$ acts transitively
      on the unit sphere in $U_1$, and hence  $\L_1/\K_1$ is two point
      homogeneous.
      Thus $\L_1/\K_1$ is isometric to a rank one symmetric space.
      From the classification of rank one symmetric spaces as homogeneous spaces
       we  see that the
representation of $\K_1$ is either of real
      or complex type, but not symplectic.

      Since the manifold is odd dimensional and $U_1$ and the slice
      have the same dimension,
       it follows that $U_2$
      is odd dimensional and therefore $\dim(U_2)\ge 3$.
      Because of $\dim(U_2)\ge 2$ there exists a $\K_1$ invariant
       irreducible subspace $U'\neq 0$
      of $U_1\otimes U_2$ contained in the kernel of $B$.

If the representation of $\K_1$ on $U_1$  is of real type,
      we claim
      that $U'$ is necessarily of the form $U_1\otimes U_2'$, where $U_2'$ is
      a one dimensional subspace of $U_2$, which
contradicts $B(u_1,u_2)\neq 0$.
      To see this, choose a basis $e_{\subo},e_1,\cdots , e_k$  of $U_2$. Any $\K_1 $
      invariant subspace of $U_1\otimes U_2$, which we can assume  projects onto
      $U_1\otimes e_{\subo}$, must be of the form $x\otimes e_{\subo} +
L_1(x)\otimes e_1 + \cdots + L_k(x)\otimes e_k$, where $x\in U_1$
and $L_i$   endomorphisms of $U_1$. To be $\K_1$ invariant implies
that $L_i$ commute with the representation of $\K_1$ on $U_1$. Since
it is of real type, this means that $L_i$ are scalar multiplication
with $\lambda_i$, and hence $e_{\subo}+\lambda_1 e_1 + \cdots + \lambda_k
e_k$ spans $U_2'$.

If the representation of $\K_1$ on $U_1$ is of complex type, we can
    repeat the previous argument in the complexifications
$U_i\otimes\C$.
      Since the  kernel of $B_{U_1\otimes U_2}$ contains $\dim(U_2)-1$
      linear independent $\K_1$ invariant irreducible subrepresentations,
       we may view these subrepresentations as a complex hyperplane
      in $U_2\otimes  \C$. Because of $\dim(U_2)\ge 3$, this hyperplane
      intersects $U_2\otimes  \R$, and we get a contradiction as before.
\end{proof}

We stress that in even dimensions, the statement of the product lemma is no longer valid in general. We will determine the exceptions in (\ref{evennotsimple}).

\bigskip

We conclude this section with a discussion of the \emph{recognition
tools} we will apply in this paper. These tools are indispensable for our
proof.

First of all by combining Straume's classification of cohomogeneity
one homotopy spheres \cite{straume} with the work of Back-Hsiang
\cite{BH} (and Searle \cite{searle} in dimension five) we have

\begin{thm}\label{thm: straume}
     Any cohomogeneity one homotopy sphere $\Sigma^n$ with an
invariant metric of positive curvature is
equivariantly diffeomorphic to the standard sphere $\Sph^n$
      with a linear action.
\end{thm}

     The same conclusion is  true for all manifolds whose rational
cohomology ring
      is like that
     of a nonspherical
rank one symmetric space (see \cite{iwata:1,iwata:2} and \cite{uchida}).

\bigskip

     The following very general recognition theorem was proved in
\cite{Wi: sym}:

\begin{thm}[Chain Theorem]\label{chain}

     Suppose $\G_d \in \{\SO(d), \SU(d),\Sp(d)\}$ acts isometrically
      and nontrivially on a positively curved compact simply connected
      manifold $M$.  Suppose also that a principal isotropy group of
      the action contains up to conjugacy a  $k\times k$~block
      with $k\ge 2$ if $\G_d=\Sp(d)$, and $k\ge 3$ otherwise.
Then $M$ is  homotopy equivalent to a rank one symmetric
      space.

\end{thm}

In conjunction with the reduction idea, the following basic
\emph{connectedness lemma} of \cite{Wi: torus}  provides another general  topological tool that will aid us in the recognition process.

\begin{thm}[Connectedness Lemma]\label{thm1}
\label{thm:connected}\label{thm: connected}\label{thm: connect}
Let $M^{n}$ be a compact positively curved Riemannian manifold, and
$N^{n-k} \subset M^{n}$ a compact totally geodesic submanifold. Then
\begin{enumerate}
\item[(a)] The inclusion map $N^{n-k}\rightarrow M^{n}$ is
$(n-2k+1)-$connected.
\item[(b)]If in addition $N^{n-k}$ is fixed pointwise by a compact
group $L$ of isometries of $M$, then
the inclusion map is $\bigl(n-2k+1+\delta(L)\bigr)-$connected, where
$\delta(L)$ is the dimension of a principal orbit of the $L$ action.
     \item[(c)] If also $V^{n-l}\subset M^{n}$ is a compact totally
geodesic submanifold, and $k\le l$, $k+l\le n$. Then the inclusion
map $N^{n-k}\cap V^{n-l}\longrightarrow V^{n-l}$
is $(n-k-l)-$connected.
\end{enumerate}
\end{thm}

As an example of a simple application of this result, combined with
Poincare duality,  we note (cf. \cite{Wi: torus}):

\begin{equation}\label{codim2}
\text{$V^{n-2} \subset M^n$ totally geodesic and $M$ positively curved $\implies
\tilde{M}$ is a homotopy sphere.}
\end{equation}

\smallskip

We finally recall that a $\G$-manifold is fixed point homogeneous if
$M^{\G}$ is non-empty and $\G$ acts transitively on the normal
spheres to a component of the fixed point set, equivalently $\dim
M/\G - \dim M^{\G} = 1$. The classification of fixed point
homogeneous manifolds with positive curvature
\cite{grove-searle:rep} will be used frequently.

\begin{thm}[Fixed Point Homogeneity]\label{thm: gs}\label{cofix}
      Let $M$ be a compact simply connected manifold of positive
      curvature. If $M$ is fixed point homogeneous, then $M$ is equivariantly
      diffeomorphic to a rank one symmetric space endowed with a linear action.
\end{thm}

Consider the special case, where one of $\Kpm$, say $\Km$ contains a
connected normal subgroup $\G' \triangleleft \G$. Let $\G''
\triangleleft \G$ be a normal subgroup with $\G' \cdot \G'' = \G$.
Clearly $\G'$ acts trivially on $\G/\Km$. Thus  if $\G'$ acts
transitively on the normal sphere $\Sph^{l_-}$, $M$ is fixed point
homogeneous. If not, $\G''\cap \Km$ acts transitively on
$\Sph^{l_-}$, and hence $\G''$ has the same orbits as $\G$ does. In
summary:

\begin{lem}\label{lem: normal subgroup} \label{lem: normal}
     If one of $\Kpm$ contains a normal connected subgroup of $\G$, then  either
there is
     a proper normal subgroup of $\G$ acting orbit equivalently, or $M$ is
     fixed point homogeneous.
\end{lem}

\no This motivates the following:

\begin{definition}\label{essential}
An action is called \emph{essential} if no subaction is fixed point
homogeneous, and no
normal subaction is orbit equivalent to it.
\end{definition}

\no Note that the above Lemma asserts in particular that:

\begin{itemize}
\item
\text{For an essential $\G$-action, none of $\Kpm$ contain a
connected normal subgroup of $\G$.}
\end{itemize}

\smallskip

In the proof of Theorem A we restrict ourselves to essential
actions. In the case that the underlying manifold is sphere this
is justified by Theorem~\ref{thm: straume}. If the underlying
is not a sphere then a cohomogeneity one
action has an essential normal subaction, and by
Lemma~\ref{extensions} below this subaction already determines
the action itself.

In the case of linear actions on spheres, a nonessential action
is either a sum action, including certain modified sum actions, or
 a $\U(1)$ extension of an essential action. The principal
isotropy groups of sum actions are transparent (see Appendix II).
The essential actions on spheres with their isotropy groups, which
we use frequently, are collected in Table \ref{odd} (and for the
even dimensional rank 1 projective spaces in Table \ref{even} ). We
include their normal  extensions since, although not essential in
the above sense, they will also be used in our induction steps.


\section{ Weyl Groups.}\label{sec: weyl}

The main objective in this section is to obtain effective upper
bounds on the Weyl groups of positively curved cohomogeneity one
manifolds, and to prove group primitivity of such manifolds.
The main result asserts that except for the cases of
 $\corank(\gH)=0$, and $\H$  finite and non-cyclic, the order of the
Weyl group divides $4\corank(\gH)\le 8$. We will first analyze the situation in the case of a trivial $\gH$
and later on reduce the general case to this one.
\smallskip

We begin with the following crucial observation

\begin{lem}\label{weyl}
      The Weyl group of a positively curved cohomogeneity one manifold is
finite.
\end{lem}

\begin{proof}
Since the Weyl group is a subgroup of $\N(\H)/\H$ our claim is obvious
when
  $\N(\H)/\H$ is finite. When $\dim
(\N(\H)/\H) > 0$ we will use the fact noted earlier, that the Weyl
group of $M$ coincides with the Weyl group of its core \eqref{reduction}. In
particular, it suffices to prove our claim for $\G$-actions with
trivial principal isotropy group. Now suppose $\W = <w_-, w_+>$ is
infinite, i.e., the Weyl group elements $w_+$, $w_-$ are involutions
in $\G$ and $w_+\cdot w_-$ generates an infinite
     cyclic group. Let $\T^h, h \ge 1$ be the identity component of the closure
     of this cyclic group. Choose a positive integer $l$ with
     $(w_+w_-)^l\in\T^h$. Clearly
     $w_-(w_+\cdot w_-)w_-=w_+(w_+\cdot w_-)w_+=(w_+\cdot w_-)^{-1}$
     and similarly
      \[
       w_-(w_+\cdot w_-)^lw_-=w_+(w_+\cdot w_-)^lw_+=
       (w_+\cdot w_-)^{-l}.
      \]
Since the infinite group generated by $(w_+\cdot w_-)^{l}$ is dense
in
     $\T^h$, it follows that the maps $\T^h \rightarrow \T^h,$
     $a\mapsto w_{\pm} aw_{\pm}$ both coincide with the inverse map $\iota: \T
\to \T$ taking $t$ to $t^{-1}$. Thus $\Ad_{w_+}v=\Ad_{w_-}v=-v$ for
all vectors $v$
     in the Lie algebra of $\T^h$.
     On the other hand, since $\Kpm$ can only be $\Z_2$, $\S^1$ or $\S^3$,
$w_{\pm}$ is central in $\Kpm$, and hence
$\Ad_{w_{\pm}}v_{\pm}=v_{\pm}$ for $v_{\pm}$ in the Lie algebra of
$\Kpm$. If we fix a biinvariant metric we deduce that the Lie
algebras of $\Kpm$ are perpendicular to the Lie algebra of $\T^h$.
Applying the same argument again on any of $w \Kpm w^{-1}$, $w\in
\W$, we see that in fact the Lie algebras of $w\Kpm w^{-1}$ for any
$w\in \W$  are perpendicular to the Lie algebra of $\T^h$.
     This contradicts linear primitivity.
\end{proof}

\smallskip

It is now possible to classify all cores with their core actions
(see also \cite{putmann} for the even dimensional case). However,
the following suffices for our
purposes:

\begin{lem}[Core-Weyl Lemma]\label{weyl-core} Suppose a Lie group $\G$
 acts with cohomogeneity one on a positively curved compact
manifold $M$ with finite fundamental group
 and trivial principal isotropy group.
Then
$\G$ has at most two components and the action is primitive.  Moreover,
 \[
\mid \gW  \mid \quad divides \quad  2 \rank(\G)\cdot | \G/\G_{\subo}|\le 8. \]

\no Furthermore
$\G_{\subo}$ is one of the groups $\S^1, \S^3, \T^2, \S^1\times\S^3, \U(2), \S^3\times\S^3, \SO(3)\times\S^3$, or $\SO(4)$, and $M$ is fixed point homogeneous in all cases but $\G_{\subo} = \SO(3)\times\S^3$.
\end{lem}

\begin{proof} First notice that the rank of $\G$ is $1$ or $2$ by
     the rank lemma. Since the group $\H$ is trivial, it follows that
     $\Kpm$ is isomorphic to one of the groups $\Z_2,\S^1$ or $\S^3$. Moreover,
at most one of $\Kpm$ is $\Z_2$ and $\G$ has at most two components
(cf.\eqref{bundlecomp}). Furthermore if $\G$ is not connected then
the Weyl group of the $\G_{\subo}$  action has index 2 in $\gW$, and the
bound follows from the connected case. It is also easy to see that
primitivity follows from primitivity in the connected case. In other
words it suffices to consider the connected groups of rank at most
two.

\smallskip

 We start by excluding the case where $\G$ is simple and without  central involution, i.e., we suppose $\G$ is one of
the groups $\SO(3), \SU(3), \SU(3)/\Z_3$, $\SO(5)$, or $\G_2$. The
Weyl group is generated by two involutions $w_-$ and $w_+$ in $\G$
and we claim that one can find elements $g\in \G$ arbitrarily close
to $e$ such that the group generated by $w_-$ and $gw_+g^{-1}$ is
infinite. This in turn implies that there are invariant metrics on
$M$ that are $C^\infty$ close to the given metric for which the
normal geodesic goes from $p_-$ to $gp_+$  and for which the Weyl
group is hence infinite. But this contradicts Lemma~\ref{weyl}. To
see the claim we assume, on the contrary, that it is false. Then we
could find a small neighborhood $U$ of $e\in \G$ and a map $k: U\to
\Z$ with $(w_- gw_+g^{-1})^{k(g)}=e$. Since for each integer $k$ the
set of all $g$ satisfying $(w_- gw_+g^{-1})^{k}=e$ is an algebraic
subvariety of $\G$, it follows that all $(w_- gw_+g^{-1})$ have a
common order independent of $g \in U$. However this is false for
each of the above groups. In all cases but $\SO(5)$, this follows
from the fact that all involutions are unique up to conjugacy, see Table
\ref{symm}.

\smallskip

The case $\G\cong \SO(3)\times \SO(3)$, where $\G$ is non-simple without central involutions is ruled out as well: As above, we can find a nearby metric with infinite Weyl group
unless  $w_-\in 1\times \SO(3)$ and $w_+\in \SO(3)\times 1$ (or vice
versa) and hence $\W\cong \Z_2\oplus \Z_2$. Since $\SO(3)\times \SO(3)$ contains no
subgroup isomorphic to $\S^3$ it follows that $\dim(\K_\pm)\le 1$, but this contradicts linear primitivity.

\smallskip

Now suppose $\G$ has central as well as non-central involutions,
i.e.,  $\G$ is one of the groups $\U(2), \S^1\times \SO(3)$, $
\SO(4)$, $ \S^3\times \SO(3)$, or $\Sp(2)$.  We can argue as before
unless one of the elements say $w_-$ is central in $\G$.  But then
$\W\cong \Z_2\oplus \Z_2$ or $\W\cong \Z_2$ and $\W$ normalizes the
group $\K_+$. From linear primitivity we see that the Lie algebras
of the groups $\K_-$, $\K_+$ and $w_+\K_-w_+$ generate the Lie
algebra of $\G$ as a vector space. Because of $\dim(\K_{\pm})\le 3$
this clearly rules out $\Symp(2)$. For the other groups it follows
that either $\K_-$ or $\K_+$ is three dimensional and thus
isomorphic to $\S^3$, so $\S^1\times \SO(3)$ is ruled out as well .
If $\G=\SO(4)$ or $\U(2)$, every $\S^3$ is normal and hence  $M$ is
fixed point homogeneous. Note that primitivity in these cases
immediately follows from linear primitivity since one of the groups
$\Kpm$ is a normal subgroup of $\G$.

\smallskip

If $\G=\S^3\times \SO(3)$ and one of  $\Kpm$ is  an
$\S^3$ factor,  $M$ is fixed point homogeneous as above, and
$\W=\Z_2$. If both $\Kpm$ are diagonal 3-spheres, we obtain a
contradiction to linear primitivity by observing that they must have
at least a one dimensional intersection. If $\Km$ is diagonal and $\Kp=\Z_2$,
the conjugates $\Km$ and $w_+\Km w_+$ also have a one dimensional
intersection.  In all other
cases, one of $\Kpm$, say $\Km$ is a diagonal $\S^3$ and  $\Kp$ is a
circle with slope $(p,q)$ in a maximal torus of $\G$. Notice that
linear primitivity also implies that $\W = \Z_2\times \Z_2$.
  We will
later determine what slopes $(p,q)$ are possible, and the
corresponding manifolds are Eschenburg spaces (cf. Section \ref{sec:
examples}). To prove primitivity in this case it is sufficient to show that $\Km$
nor a conjugate of $\Km$ can  be a subgroup of $\Kp$. But under such an assumption, we would have that $w_+=w_-$ and thus $W\cong \Z_2$,
 contradicting the Lower Weyl Group Bound.

\smallskip

It remains to consider the cases where all involutions of $\G$ are central, i.e., $\G$ is one of the groups $\S^1$, $\S^3$, $\S^1\times \S^1$,
$\S^1\times \S^3$, $\S^3\times \S^3$.  Clearly the order of the Weyl group is at most $2\rank
\G$. From linear primitivity it follows that the Lie algebras of
$\Kpm$ generate the Lie algebra of $\G$ as a vector space. This
implies that at least one of the groups $\Kpm$ is normal and $M$ is
fixed point homogenous and primitive.

\vspace{5pt}

\end{proof}

%

We can now use the above lemma
and the last paragraph of section~\ref{sec: general}
to prove the group primitivity stated
in \eqref{pim}:

\begin{cor}[Group Primitivity]\label{prim}
     Suppose that $M$ admits a positively curved cohomogeneity one
metric.
     Consider any other cohomogeneity one metric on $M$, then the
corresponding
     groups $\Km, \Kp$ generate $\G$ as a Lie group. Equivalently
     $\Km$ and
     $n \Kp n^{-1}$ generate $\G$ for any
 $n\in \N(\H)_{\subo}$.
\end{cor}

\begin{proof} Let $\Kpm$ denote the isotropy groups with respect
to
     a positively curved metric. By linear primitivity  $\Km$ and
     $\Kp$
     generate $\G$ as a group. We need to show that for any $a\in \N(\H)_{\subo}$,
 the groups
$\Km$ and
     $a \Kp a^{-1}$ generate $\G$ as well. But by primitivity of
     the core,
we know that  $\Km \cap \N(\H)_c$ and
     $a (\Kp  \cap \N(\H)_c) a^{-1}=a \Kp a^{-1} \cap \N(\H)_c $
generate the core group. In particular, the group generated by $\Km$
and
     $a \Kp a^{-1}$ contains $\N(\H)_{\subo}$, and hence is equal to
     the group generated by  $\Km$ and $\Kp$.
\end{proof}

\smallskip

We have the following  useful consequence of primitivity:

\begin{lem}\label{kernel}
Assume $\G$ acts effectively. Then the intersection $\H_- \cap \H_+$ of
the ineffective kernels $\H_{\pm}$ of $\Kpm/\H$ is trivial.
\end{lem}

\begin{proof}
We first observe the following: If for a connected homogeneous space
$\K/\H$, a normal subgroup $\L$ of $\H$ acts trivially on $\K/\H$,
then $\L$ is normal in $\K$ also. Indeed, first observe that
$\N(\L)$ acts transitively, since it in general acts with finite
quotient on the fixed point set of $\L$. Hence
$\K/\H=\N(\L)/(\N(\L)\cap \H)=\N(\L)/\H$  and thus $\K=\N(\L)$. In
our case, we can apply this to the normal subgroup $\H_- \cap \H_+$
of $\H$ which fixes both $\Sph^{l_{\pm}}$. Thus
 $\Kpm \subset\N(\H_- \cap
\H_+) $, and hence by primitivity $\N(\H_- \cap \H_+) = \G$. Since
the action is effective, $\H_- \cap \H_+$ is trivial.
\end{proof}

When $M$ is simply connected and $\G$ is connected, we recall from
\eqref{Hcomp} that $\Kpm$ and $\H$ are all connected as long as both
$l_{\pm} \ge 2$. If exactly one of $l_{\pm}$ is $1$, say $l_- = 1$
and $l_+\ge 2$,  $\Km$ is connected, $\H/\Ho = \Kp/\Kpo$ is cyclic, and
$\H = \H_-$. If in addition $\G$ is assumed to act effectively, it
follows from the above Lemma \ref{kernel} that $\Kp$ acts
effectively on $\Sph^{l_+}$. In the remaining situation where both
$l_{\pm} = 1$, Lemma \ref{kernel} and $|\H/\H_\pm|\le 2$  yield:

\begin{lem}\label{both-one}
Suppose $M$ is a 1-connected positively curved manifold on which the
connected group $\G$ acts effectively and isometrically with codimension
two singular orbits. Then one of the following holds:
\begin{enumerate}
\item[(a)]
$\H = \{1\}$ and both $\Kpm$ are isomorphic to $\SO(2)$.
\item[(b)]
$\H = \H_- = \Z_2$, $\Km = \SO(2)$ and $\Kp = \O(2)$.
\item[(c)]
$\H = \H_- \cdot \H_+ = \Z_2 \times \Z_2$, and both $\Kpm$ are
isomorphic to $\O(2)$.
\end{enumerate}
\end{lem}

\smallskip

Notice that part (a) of \eqref{both-one} is not possible  when
$\rank \G\ge 2$ since the action would then not be group primitive
due to the fact that both $\Kp$ and $\Km$ can be conjugated into a
common maximal torus.

As a consequence of the Core-Weyl Lemma one obtains an important
upper bound for the Weyl group:

\begin{prop}[Upper Weyl Group Bound]\label{upper}
Assume that $M$ is simply connected and $\G$ connected. Then
\begin{enumerate}
\item[(a)]

If  $\H/\Ho$ is  trivial or cyclic, we have $|\W|\le 8$ if the
corank of $\H$ in $\G$ is two, and \newline $|\W|\le 4$ if the
corank is one.

\item[(b)]
If $\H$ is connected and $l_\pm$ are both odd, $|\W|\le 4$ in the
corank two case and $|\W|\le 2$ in the corank one case.

\item[(c)] If none of $\N(\H)\cap\Kpm$ are  finite, $|\W|\le 4$ in the
corank two case and $|\W|\le 2$ in the corank one case.
\end{enumerate}
\end{prop}

\begin{proof}
We first consider the case where $\H/\Ho$  is non-trivial and
cyclic. Then \eqref{Hcomp} and \eqref{both-one}  imply that the
codimension of one of the orbits is  two and  one of the
corresponding  $\K$ groups is connected. Thus $\N(\H)/\H$ is not
finite since $\K\subset \N(\H)$. By passing to the reduction
$M^{\H}$, we deduce from the Core-Weyl Lemma \ref{weyl-core}
 that $|\W|\le 8$ ( $|\W|\le
4$ in the corank one case).

Now assume that $\H$ is connected. If
$\H=\{e\}$, the claim follows again from the Core-Weyl Lemma.
   If not, fix a maximal torus $\T \subset \H$. Clearly then
$M^{\T}$ has positive dimension. By Lemma~\ref{reduction}, the group
$\N(\T)_c$ acts on the reduction $M^{\T}_c$ with
 the same Weyl group. By \eqref{bundlecomp}, the Weyl
group of $\N(\T)_{\subo}/\T$ has index at most two in $\W(\G,M)$.

Next observe that for any torus $\T$ of a connected compact Lie
group $\G$, $\N(\T)_{\subo} \subset  \gC(\T)$, the centralizer of $\T$ in
$\G$. Because $\H$ is a connected Lie group $\T$ is maximal abelian
in $\H$ and thus $\gC(\T)\cap \H=\gT$.  Hence $\N(\T)_{\subo} \cap \H =
\T$ and thus $\N(\T)_{\subo}/\T$ acts with trivial principal isotropy
group on the reduction $M^{\T}_c$. It follows that $|\W|\le 8$ (
$|\W|\le 4$ in even dimensions) by the Core-Weyl Lemma.

Since the codimension of $\Sph_{\pm}^{\T} \subset \Sph^{l_{\pm}}$ is
always even, $\Sph_{\pm}^{\T}\ncong \Sph^0$ if both $\lpm$ are odd
and hence \eqref{bundlecomp} implies that   $\N(\T)_c/\T$ and
$\N(\T)_{\subo}/\T$ have the same Weyl group, which implies part (b).

 For part (c) just note that by assumption both normal spheres of the
 core $M^{\H}_c$ have positive dimension. As we have seen then $\N(\H)_c$
 and its identity component have the same orbits and Weyl group. Thus from  Core-Weyl Lemma  $|\W|\le 4$
 ( $|\W|\le 2$ in even dimensions).
\end{proof}

\begin{rem} The only cases where we have no bound on the Weyl group
are hence when $\H$ has corank zero, or when  $\H$ has corank one or
two and $\H=\Z_2\oplus\Z_2$.

In the equal rank case, $\N(\H)/\H$ is always finite, and hence the
Core-Weyl Lemma does not apply. However, in this case, information
about the Weyl group does not enter in the proof of Theorem A. It
will follow as a consequence of the proof that $\W$ is one of $
\D_1, \D_3,\D_4,\D_6$.

If $\H = \Z_2 \oplus \Z_2$, we note that $\N(\H)/\H$ is also finite
since each of $(\N(\H) \cap \Kpm)/\H$ is and $M^{\H}$ is primitive.
In fact this is the case where the Weyl group can become larger. One
easily sees that the Weyl groups  are $\D_3$ for  $P_{2k}$  and
$\D_6$ for $P_{2k+1}$, whereas for  $Q_k$ and  $R$ it  is always
$\D_4$. Hence, as a consequence of our classification, it follows
that the
  Weyl groups for simply
connected positively curved \co one manifolds are the same as
for linear actions on spheres. Notice also that there are many
actions among the linear actions on spheres, for example all tensor
product actions, where $\W=\D_4$, and some of those with $\lpm$ odd
and $\H$ not connected (see Table \ref{odd}).

\end{rem}


\section{Examples and Candidates.}\label{sec: examples}

To aid the induction step in our proof of Theorem A it is important
to know more details about the individual manifolds and actions that
occur. The linear actions are of course well known, and the
essential ones and their normal extensions are exhibited in Tables
\ref{odd} and \ref{even} in Appendix II. The corresponding details
for the remaining spaces and actions, i.e., for the known
non-spherical cohomogeneity one manifolds of positive curvature (the
second part of Theorem A), and for our new candidates (third part of
Theorem A), is provided in the following Table
\ref{examples}.
In the next seven sections we show that the list is complete.
Indeed all the cases in which
nonspherical examples occur are covered by
Lemma~\ref{s3s3pos}
and Proposition~\ref{rank3}.

In this section, we will explain which
of these actions
correspond to the known cohomogeneity one manifolds of positive curvature.
  The information in the Table is separated into
homogeneous examples, biquotients, and candidates (with some
overlap). Due to its special significance we have included as a
separate entry  the linear action of $\SO(4)$ on $\Sph^7$ and
separated the two \co one actions on the Aloff Wallach space $W^7$
by its lower index. All manifolds are assumed to be simply
connected.

\smallskip

For subgroups $\S^1\subset \S^3\times \S^3$ we have used the
notation $\gC^i_{(p,q)}=\{(e^{pi\gt},e^{qi\gt})\mid \theta \in \R\}$
and $\gC^j_{(p,q)}=\{(e^{pj\gt},e^{qj\gt})\mid \theta \in \R\}$ and
$\Q$ denotes the quaternion group $\{\pm 1,\pm i, \pm j , \pm k\}$.

\renewcommand{\thetable}{\Alph{table}}
\renewcommand{\arraystretch}{1.4}
\stepcounter{equation}
\begin{table}[!h]
      \begin{center}
          \begin{tabular}{|c|c|c|c|c|c|c|c|}
\hline $M^n$ & $\G$          &$\Km$  & $\Kp$ & $\H$ & $\Hb$
&$(l_-,l_+)$ & $\W$ \\
\hline \hline

$\Sph^7 $ &$\S^3\times\S^3$ &$ \gC^i_{(1,1)} \gH$ & $  \gC^j_{(1,3)} \gH$  & $
\Q$
  & $\Z_2\oplus\Z_2 $  &$(1,1)$ &$\D_6$ \\
\hline

$B^7 $ &$\S^3\times\S^3$ &$ \gC^i_{(3,1)} \gH$ & $ \gC^j_{(1,3)} \gH$ & $\Q
$
  & $\Z_2\oplus\Z_2$ &$(1,1)$ &$\D_3$ \\
\hline

$W^7_{(1)} $ &$\S^3\times\S^3$ & $\Delta\S^3\cdot\H$ & $\gC^i_{(1,2)}$
&
$\Z_2$ & $ 1 $ &$(3,1)$ &$\D_2$ \\
\hline

$W^7_{(2)} $ &$\S^3\times\S^3$ & $\gC^i_{(1,1)} \gH$ & $\gC^j_{(1,2)} \gH$
  & $\Z_4\oplus\Z_2$ & $  \Z_2 $ &$(1,1)$ &$\D_4$ \\
\hline

$B^{13}$ &$\SU(4)$ & $\Sp(2)\cdot\Z_2$ &
$\SU(2)\cdot \S^1_{1,2}$ & $\SU(2)\cdot\Z_2$ & $\SU(2)\cdot\Z_2  $
&$(7,1)$ &$\D_2$ \\
\hline\hline

$ \E_p^7, p\ge 1$ &$\S^3\times\S^3$ & $\Delta\S^3\cdot \H$ &
$\gC^i_{(p,p+1)}$ & $\Z_2$ & $ 1 $ &$(3,1)$ &$\D_2$ \\
\hline

$B_p^{13}, p\ge 1$ &$\SU(4)$ & $\Sp(2)\cdot\Z_2$ &
$\SU(2)\cdot \S^1_{p,p+1}$ & $\SU(2)\cdot\Z_2$ & $\SU(2)\cdot\Z_2  $
&$(7,1)$ &$\D_2$ \\
\hline\hline
$P_k, k\ge 1$ &$\S^3\times\S^3$ & $\gC^i_{(1,1)}\gH$ &
$\gC^j_{(2k-1,2k+1)} \gH$ & $\Q $ & $ \Z_2\oplus\Z_2  $ &$(1,1)$ &$\D_3$
or $\D_6$ \\
\hline

$Q_k, k\ge 1$ &$\S^3\times\S^3$ & $\gC^i_{(1,1)} \gH$ &
$\gC^j_{(k,k+1)}  \gH$ & $\Z_4\oplus\Z_2 $ & $ \Z_2  $ &$(1,1)$ &$\D_4$ \\
\hline

$R$ &$\S^3\times\S^3$ & $\gC^i_{(3,1)}\gH$ &
$\gC^j_{(1,2)} \gH$ & $\Z_4\oplus\Z_2  $ & $ \Z_2  $ &$(1,1)$  &$\D_4$ \\
\hline

\hline
          \end{tabular}
      \end{center}
      \vspace{0.1cm}
      \caption{Known examples  and
      candidates.}\label{examples}
\end{table}

Some explanations are in order. The embedding of $\H$ is not always
explicitly given, but can be determined in each case.
$\Z_4\oplus\Z_2 $ is always  embedded as $\{(\pm 1 ,\pm 1),(\pm i ,
\pm i)\}$. Otherwise, a $\Z_2$ inside $\H$ is always embedded in the
circle inside $\Kp$. The embedding of $\Q$ depends on the slopes,
although it is always embedded diagonally up to conjugacy. E.g. for
$B^7$ it must be of the form $\{\pm (1,1), \pm (i,-i) , \pm
(j,-j),\pm (k,k)\}$. The embedding of $\SU(2)$ is in a $2\times 2$
block in $\SU(4)$.

Most of these actions are only almost effective, i.e. $\G$ and $\H$
have a finite normal, hence central  subgroup in common. The
effective version can easily be determined in each case, and we
include in our Table the most important part, the effective group
$\Hb$. It is also important to notice that the full effective groups
for $P_k$ are $\Z_2\oplus\Z_2\subset \{\O(2),\O(2)\}\subset_k
\SO(4)$ and for $\Q_k$ (as well as for $R$) are $\Z_2\subset
\{\SO(2),\O(2)\}\subset_k \SO(3)\SO(3)$.  Here the groups $\Km$ and
$\Kp$ are embedded in different blocks in each component of
$\SO(3)\SO(3)$. The isomorphism types of these groups are consistent
with, and in fact determined, by Lemma \ref{both-one}.

There are obvious and important isomorphisms among some of these \co one
actions which are apparent from the tables: $P_1=\Sph^7 \, , \,
Q_1=W^7_{(2)} \, , \, E_1=W^7_{(1)}$ and $B^{13}_1=B^{13}$.

The Weyl groups can be  computed from the given isotropy groups. For
example in the case of $P_k$, one chooses $w_-=(e^{\pi i/4},e^{\pi
i/4})$ and $w_+=(e^{\pi j/4},(-1)^k e^{\pi j/4})$  as representatives. One then checks
that   $(w_-w_+)^3=1$  in $\N(\H)/\H$ for $k$ even, and $(w_-w_+)^6=1$
for  $k$ odd. Hence $\W=\D_3$ for $k$ even and $\W=\D_6$ for $k$
odd.

\bigskip

The cohomogeneity one actions on the known positively curved
manifolds were discovered by the first and last author in 1997, see
\cite{Z} and \cite{GSZ}. Although one can determine the group
diagrams for these actions directly, it will be much simpler for us
to use
 the classification.
More precisely we will use
 Lemma~\ref{s3s3pos} and Proposition~\ref{rank3} from below, whose proofs
are independent
 of this section.
\bigskip

\begin{center}
$\Sph^7$ with $\G=\SO(4)$
\end{center}

\smallskip

The 7-sphere has a cohomogeneity one action by $\SO(4)$ given by the
isotropy representation of the symmetric space $\G_2/\SO(4)$. A
normal subgroup $\SU(2)$ of $\SO(4)$ acts freely on $\Sph^7$ and
hence is given by  the Hopf action. If we divide by this action, we
obtain an induced action of $\SO(3)$ on $\Sph^4$, which must be
given by the usual action on trace free symmetric $3\times 3$
matrices. The isotropy groups of this action on $\Sph^4$ are given
by $\Km=\O(2)$, $\Kp=\O(2)$, and $\H=\Z_2\oplus\Z_2$
  and hence are the same for the
$\SO(4)$ action on $\Sph^7$. Since $\SU(2)$ acts freely, the slopes
for the circles $\Kpmo$, viewed  as subgroups of $\S^3\times\S^3$,
must have $\pm 1$ in the second component. Using Lemma
\ref{s3s3pos}, the slopes must be $(1,1)$ and $(3,1)$  and this
completely determines the group picture.

\bigskip

\begin{center}
$B^7 =\SO(5)/\SO(3)$ with $\G=\SO(4) $
\end{center}

\smallskip

In the positively curved homogeneous Berger space $\SO(5)/\SO(3)$
the subgroup $\SO(3)$ is embedded via the irreducible representation
of $\SO(3)$ on trace free symmetric $3\times 3$ matrices (see
\cite{Be}). Notice that $\SO(4)\setminus\SO(5)/\SO(3)=\Sph^4/\SO(3)$
is one dimensional and thus $\SO(4)$ acts on $\SO(5)/\SO(3)$
 by cohomogeneity one.
Next we observe that the extended $\O(4)$ action is not orbit
equivalent to the $\SO(4)$ action since for the $\SO(3)$ action on
$\Sph^4$ the antipodal map takes one singular orbit to the other.
This implies that the two singular isotropy groups $\gK^-$ and
$\gK^+$ are isomorphic up to
 an outer automorphism of $\SO(4)$.
 Combining this property with
Lemma~\ref{s3s3pos} we see that
the action is determined: both singular groups are
1 dimensional and that the slopes of the circles of the corresponding
ineffective
$\gS^3\times \gS^3$-action are given by
$\{(3,1), (1,3)\}$.
\bigskip

\begin{center}
$E^7_p$ with $\G=\SO(3)\times\SU(2) $
\end{center}

\smallskip

The Eschenburg space $E^7_p = \diag(z,z,z^p) \backslash
\SU(3)/\diag(1,1,\bar{z}^{p+2}) , p\ge 1 $ has positive curvature
according to \cite{Es2}. The group $\SU(2)\times\SU(2)$ acting from
left and right in the first two coordinates induces an action on
$E^7_p$ whose orbit through the identity is $
\SU(2)\times\SU(2)/(\triangle \S^3\cdot \H)=\RP^3$ with
$\H=\Z_2=\langle(1,-1)\rangle$ or $\langle(-1,1)\rangle$. One easily sees that the action of
$\Km$ on the slice is nontrivial and hence $E^7_p$ is \co one. The
group $\Kp$  is in this case not determined by this information. A
computation  shows it is a circle with slope $(p+1,p)$ and hence
$\H=((-1)^{p+1},(-1)^p)$, see \cite{GSZ}. For $p$ even, the left
hand side $\SU(2)$ acts effectively as $\SO(3)$ and for $p$ odd, the
right hand side one does. For $p=1$ we obtain the \co one picture
for $W^7_{(1)}$ and the right hand side $\SO(3)$ acts freely. For
$p=2$ the left hand side $\SO(3)$ acts freely, as one sees
immediately from the group picture.

\bigskip

\begin{center}
$W^7_{(2)}$ with $\G=\SO(3)\times\SO(3) $
\end{center}

\smallskip
 For the positively curved Aloff-Wallach space
$W^7=\SU(3)/\diag(z,z,\bar{z}^2)$ \cite{AW}, we have
$\N(\H)/\H=\U(2)/\H=\SO(3)$ which acts freely on the right and hence
we can write $B^7=\SU(3)\SO(3) / $ $\U(2)$ (see \cite{Wi}).
Furthermore the second factor acts freely on $W^7$, and the  action
descends to the natural cohomogeneity one action of $\SO(3)$ on
$\CP^2=W^7/\SO(3)$. Thus $\G$ acts by cohomogeneity $1$. From
Lemma~\ref{s3s3pos} it follows that there is only one cohomogeneity
one action of $\SO(3)^2$ on a positively curved simply connected
$7$-manifold for
which one the factors acts freely. Thus the action is determined,
both singular isotropy groups are one dimensional and that the
slopes are given by $\{(1,2),(1,1)\}$.

\bigskip

\begin{center}
$B^{13}_p$ with $\G=\SU(4) $
\end{center}

\smallskip

The Bazaikin space $B^{13}_p =
\diag(z,z,z,z,z^{2p-1})\backslash\SU(5)/\Symp(2)
\diag(1,1,1,1,\bar{z}^{2p+3}), p \ge 1$ has positive curvature by
\cite{Ba} (see also \cite{Z} and  \cite{DE}). The action of
$\SU(4)\subset\SU(5)$ on the left induces an action on $B^{13}_p$
whose orbit through the identity is $ \SU(4)/(\Sp(2)\cup
i\Sp(2))=\RP^5$. The action on the slice is easily seen to be
nontrivial and hence $B^{13}_p$ is \co one. From the proof of
Proposition \ref{rank3} in the case of $\G=\SU(4)$ it follows that
$\H=\SU(2)\cdot\Z_2$ and $\Kp=\SU(2)\cdot\S^1$ where $\S^1$ is
allowed to have slopes $(q,q+1)$ inside of a maximal (two
dimensional) torus of the centralizer of  $\H$. We can now consider
the fixed point set of the involution $\diag(-1,-1,1,1,1)\in\SU(5)$
as in \cite{Ta} and one shows that it's fixed point set is
$\diag(z,z,z^{2p-1})\backslash \SU(3)/\diag( {z}, {z}
,\bar{z}^{2p+3} ) = \diag(z,z,z^{p})\backslash \SU(3)/\diag( 1,1
,\bar{z}^{p+2} ) = E^7_{p}$ (see \cite{DE}). Hence the slopes of the
$\SU(4)$ action are determined (i.e. $q=p$). Because of
$B^{13}_1=B^{13}$, this group picture is determined as well.

\bigskip

We add the following information about these actions, needed in our
proof:

\begin{lem}[Extensions]\label{extensions}
The nonlinear actions in Table \ref{examples} have the following
extensions:
\begin{enumerate}
\item[(a)]

The manifolds $B^7, P_k,Q_k$, and $R$, with  their natural \co one
action, do not admit any connected normal extensions.

\item[(b)]

For the manifolds $E_p$ and $B^{13}_p$, the natural action has a
unique connected normal extension by $\S^1$.

\end{enumerate}
\end{lem}

\begin{proof}
For the spaces $B^7, P_k,Q_k$, and $R$, which have singular orbits
of codimension two, the identity component of the principal isotropy
group of the extended action would normalize both singular isotropy
groups contradicting primitivity.

For the spaces $E_p$ and $B^{13}_p$, the natural action has a
$\U(1)$ extension, since e.g. $\SU(4)\subset\SU(5)$ lies in $\U(4)$.
Since the group diagram of this extension can be derived from that
of $\G$, any two extensions are equivariantly diffeomorphic.
\end{proof}

One also easily derives the following information from the group
diagrams in Table \ref{examples} and Table \ref{odd}.

\begin{lem}[Free Actions]\label{free}
If $\G$ acts by \co one on  an odd dimensional simply connected
positively curved manifold  $M$  and there exists a  subgroup
$\L\subset\G$ with $\L=\SU(2)$ or $\L=\SO(3)$ which acts freely,
then

\begin{enumerate}
\item[(a)]

$M=E_1=W^7_{(1)}$ or $M=E_2$ with  $\L=\SO(3)\subset
\SO(3)\SU(2)=\G$.

\item[(b)]

$M=W^7_{(2)}$  with  $\L=\SO(3)\subset \SO(3)\SO(3)=\G$.

\item[(c)]

 $M$ is a sphere and the subaction of $\L\cong \S^3$ is given
by the Hopf action.

\end{enumerate}
\end{lem}

\begin{rem}\label{CP2Bundle}
 The existence of the free $\SO(3)$ actions on $E_1$ and $E_2$
 was first
observed by Shankar in \cite{shankar:chern}, in connection with his
discovery of counter examples to the so-called Chern conjecture. In
the case of $E_1=W^7_{(1)}$ and $W^7_{(2)}$ it is the natural free
action of $\N(\H)/\H$ on $W^7$.

Also notice that in all three  cases the quotient by $\SO(3)$ is
equal to $\CP^2$, which one can recognize from the induced \co one
diagram on the base. In the case of $E_1$ and $E_2$ it is the action
of $\SU(2)$ on $\CP^2$ which has a fixed point. In the case of
$W^7_{(2)}$ it is the \co one action by $\SO(3)$ with singular
orbits of codimension two.

\end{rem}

The proof of Theorem A will occupy the next 7 sections. As stated
earlier, this is achieved by classifying essential  cohomogeneity
one actions by compact connected groups on simply connected odd
dimensional manifolds with positive (sectional) curvature.

All partial
classification results will be formulated in Propositions, and

\begin{itemize}
\item
for simplicity we will abuse language and
assume from now on without
stating it explicitly, that the manifolds $M$  under consideration are simply
connected and positively curved.
\end{itemize}

When a manifold is recognized via its isotropy groups, we will often
say that we have ``recovered'' a particular action and manifold and
leave it up to the reader to find the corresponding entry in Tables
\ref{odd} or \ref{even} and to verify that the groups are indeed
recovered up to equivalence of their diagrams.


\section{Equal Rank Groups.}\label{sec: equalrank}

We are now ready to begin our classification of essential isometric
cohomogeneity one $\G$-actions on simply connected positively curved
manifolds $M$. This section is concerned with the simplest situation
of the rank lemma, where

\begin{itemize}
\item
$\rank(\H) = \rank(\Km) = \rank(\Kp) =\rank(\G)$
\end{itemize}

\no In particular, the normal spheres

\begin{itemize}
\item
$\Sph^{l_{\pm}} = \Kpm/\H$ are even dimensional
\end{itemize}

\no  and hence one of $\SO(2n+1)/\SO(2n)$ or $\G_2/\SU(3)$. Thus

\begin{itemize}
\item
$\H \subset \{\Km,\Kp\} \subset \G$ are all connected.
\end{itemize}

\no Since an equal rank subgroup of $\G_1\cdot \G_2$ is of the form
$\H_1 \cdot \H_2$ with $\H_i\subset \G_i$, $\G$ is clearly semisimple,
and hence by the product Lemma

\begin{itemize}
\item
$\G$ is simple.
\end{itemize}

\no Since the weights of the isotropy representation of an equal
rank subgroup are roots, we have

\begin{itemize}
\item
The irreducible subrepresentations $\fm_i$ of $\H$ are
      pairwise non-equivalent.
\end{itemize}

\bigskip

We will divide our analysis into the following three cases: (1) $\H$
is not semisimple, (2) $\H$ is semisimple, but not simple, and (3)
$\H$ is simple.

\begin{prop}\label{equalrank1}
If $\G$ acts essentially, with non-semisimple $\H$ of corank zero, then
$\G$ is one of $\SU(3), \Sp(2)$, or $\G_2$ and the action is the
adjoint representation restricted to the sphere.
\end{prop}

\begin{proof}
We  first show that in fact $\H$ is a maximal torus $\T$. If not,
let $\H' \triangleleft \H$ be a simple connected normal subgroup,
and $\S^1 \subset Z(\H)$. Since $\Kpm/\H$ are even dimensional
spheres, either $\H'$ or $\S^1$ must act trivially on the
irreducible subrepresentations of $\H$ in $\Kpm$. By the isotropy
lemma the same then holds for each irreducible subrepresentation of
$\H$ in $\G$ and we obtain a contradiction as in the proof of Lemma
\ref{Hrank2}.

Therefore $\H = \T$ and we conclude that $\Sph^{l_{\pm}} \cong \Sph^2$, and
$\H/\H_{\pm}$ both circles. By primitivity  we see that $\dim \T =
\rank \G \le 2$. If $\rank \G = 1$ the action is obviously a
suspension action which is non essential. It follows that $\G$ is
one of $\SU(3), \Sp(2)$, or $\G_2$.

      To unify the discussion of these three cases we will use
the well known fact (see e.g. \cite{wolf}) that the Weyl group,
$\N(\T)/\T$ of
$\G$ acts transitively on the set of roots of $\G$ of the same length.

The Weyl group of $\SU(3)$ is $\D_3$ acting transitively on its set
of three equal length roots. Each root corresponds to a $\U(2)
\subset \SU(3)$, and by primitivity the pair $(\Km,\Kp)$ must be a
pair of $\U(2)$ subgroups of $\SU(3)$ corresponding to different
roots. We have recovered the diagram for the adjoint action of
$\SU(3)$ on $\Sph^7$.

Both $\Sp(2)$ and $\G_2$ have roots of two lengths. From the
Isotropy Lemma it follows that the singular isotropy groups must
correspond to roots of different lengths.

The Weyl group of $\Sp(2)$ is $\D_4$ with two long roots
      $\Sp(1)\times \S^1 \subset \Sp(1)\times \Sp(1)
\subset \Sp(2)$  and  two short roots  $\U(2)\subset \Sp(2)$. All
pairs  $(\Km,\Kp)$ corresponding to a long and a short root define
the same manifold, namely  $\Sph^9$ with the adjoint action of
$\Sp(2)$.

The Weyl group of $\G_2$ is $\D_6$, and  has three long roots and
three short roots. A short root corresponds to $\U(2)\subset\SU(3)$.
There are two $\U(2)\subset\SO(4)$, one a long root and one a short
root. Since $\Kpm$ cannot both be in $\SO(4)$ by primitivity, this
leaves, modulo the action of the Weyl group, only one configuration
for the pairs $(\Km,\Kp)$ and we have recovered the adjoint action
of $\G_2$ on $\Sph^{13}$.
\end{proof}

\begin{prop}
If $\G$ acts essentially, with semisimple, nonsimple $\H$ of corank
zero, then $\G =\Sp(3)$ and the action is the unique linear action on
$\Sph^{13}$ with $\H = \Sp(1)^3$.
\end{prop}

\begin{proof}
Suppose $\H'$ is a simple normal subgroup of $\H$ with $\rank \H'
\ge 2$. Similarly to Lemma \ref{Hrank2},  we can find a
subrepresentation on which $\H'$ and $\H/\H'$ act non-trivially,
which can not degenerate since $\Kpm/\H$ are even dimensional
spheres. Hence, by assumption $\H$ is a semisimple group with rank
one factors only. In particular both $\Sph^{l_{\pm}}$ are
4-dimensional.

As above, we see that for any two different simple subgroups $\H_1$
and $\H_2$, the isotropy representation of
      $\H$ has an irreducible subrepresentation $\fm$ on which both $\H_i$
act non trivially. By the isotropy lemma, this representation has to
degenerate along
      the normal geodesic $c$ at some singular orbit, say
$\K/\H=\Sp(2)/\Sp(1)\Sp(1)$. Note that
there is an element $w \in \W$ represented by an element $w\in \K \cap
\N(\H)$, which  acts on $\H$
      by permuting the two factors $\H_1$ and $\H_2$, and leaving
      all other factors of $\H$ invariant.
      Thus the action of Weyl group on the factors of
      $\H$ contains all possible transpositions, and
      it is hence the full symmetric group. The only symmetric groups which are
      dihedral are $\S_2$ and $\S_3$. Hence $\H$ has at most three factors
      or equivalently $\rank(\G)\le 3$. If $\rank(\G) = 2$, $\G$ must contain
an
      $\Sp(2)$ or $\SO(5)$, which rules out  $\G = \SU(3)$ and $\G_2$,
       and for $\G = \Sp(2)$ the action must be a suspension
action, which is not essential.

If $\rank(\G) = 3$, $\G$ contains a semisimple
      $9-$dimensional subgroup $\H$ as well as an $\Sp(2)\Sp(1)$, which
      rules out $\SU(4)$ and $\SO(7)$, and in the case of
      $\G=\Sp(3)$ with $\H=\Sp(1)^3$ leaves, by primitivity, only
      one configuration for $\Kpm$ and we have recovered the action of
      $\Sp(3)$ on $\Sph^{13}$.

\end{proof}

\begin{prop}
If $\G$ acts essentially, with simple $\H$ of corank zero, then $\G =
\F_4$, and the action is the unique linear action on $\Sph^{25}$ with
$\H = \Spin(8)$.
\end{prop}

\begin{proof} Using that $\gH$ is a simple equal rank subgroup of $\G$
with a spherical isotropy representation, we can deduce from
Table~\ref{sph} that $(\G,\H)$ is either
$(\F_4,\Spin(8))$ or $(\F_4,\Spin(9))$.
The latter case can actually not occur since the 16-dimensional
representation of $\F_4/\Spin(9)$ can not possibly degenerate.
Recall that the isotropy representation of $\F_4/\Spin(8)$ decomposes
      into three pairwise nonequivalent
      $8$ dimensional representations of $\Spin(8)$, each contained in
a $\Spin(9)$. Clearly the action is determined by primitivity, and
we have recovered the unique cohomogeneity one action of $\G = \F_4$
on $\Sph^{25}$.
\end{proof}

We point out that for all actions classified in this section
the
 cohomogeneity one Weyl groups coincide with the core groups
 $\N(\H)/\H$ which are either $\D_3$, $\D_4$ or $\D_6$.


\section{Non Semisimple Groups.}\label{sec: not semi}

In this and the following five sections we assume that:
\begin{itemize}
\item
$M$ is simply
connected cohomogeneity one $\G$-manifold, with an invariant metric\\ of
positive curvature,
\item
$\G$ is connected acting essentially with principal isotropy group
$\H$ of corank two.
\end{itemize}

\no Based on the even dimensional classification
\cite{verdiani:1,verdiani:2}, the following is quite simple:

\begin{prop}
Suppose $\G$ is  not semisimple and acts essentially with corank 2.
Then either $\G=\S^1\cdot \L$, where $\L$ is one of $\SO(n),
\Spin(7)$, or $\G_2$, and the action is a tensor product action on
$\Sph^{2n-1}, \Sph^{15}$, or $\Sph^{13}$ respectively. Otherwise
$\G=\U(2)\SU(2)$ with its tensor product action on $\Sph^7$.
\end{prop}

\begin{proof}
After passing to a finite covering of $\G$ we may assume
      $\G=\S^1\times \L$. Since $\H\cap \S^1$ is in the ineffective kernel of
the action we can assume it is trivial. Moreover, $\H$ does not project
surjectively onto $\S^1$, since otherwise  the subaction of $\L$ would be
orbit equivalent to the $\G$-action, which would then not be essential.
Assume first that the subaction of the $\S^1$-factor is free. Then
$B=M/\S^1$ is
      an even dimensional simply connected
       manifold of positive sectional curvature with
      a cohomogeneity one action of $\L$, and $B$ is not $2$-connected.
      So Verdiani's classification implies that $B$ is a complex projective
space.
      Since $M$ is simply connected, the Euler class of the bundle
      $\S^1\rightarrow M\rightarrow B$ is a generator of $H^2(B,\Z)$.
      Using the Gysin sequence we deduce that $M $ is a homology sphere.

If the subaction of the $\S^1$-factor is not
free, we can assume
      without loss of generality that $\Km \cap \S^1\neq 1$.
      Since $\S^1 \cap \H = 1$, $\Km \cap \S^1$ acts freely
      on $\Km/\H$ and hence  $\G/\Km$ is a component of the
fixed point set $M^{(\Km \cap \S^1)}$.
      By assumption (cf. \ref{lem: normal}) $\Km$ is not normal in
$\G$, and
$\dim(\G/\Km)>1$. Moreover,
      $\Km$ must project surjectively to $\S^1$, since $\G/\Km$ has positive
curvature and hence finite fundamental group.
On the other hand, since $\H$ does not project surjectively to $\S^1$,
it follows that $\G/\Km$
      has codimension $2$, and thus $M$ is a homotopy sphere by the
connectedness lemma (cf. \ref{codim2}).

The actual determination of the action follows from Straume's
classification (see Table \ref{odd}).
\end{proof}


\section{Semisimple Rank $2$ Groups.} \label{sec: rank2}

In the next four sections we assume in addition to $M$ being a simply
connected cohomogeneity one $\G$-manifold, with an invariant metric of
positive curvature, that:
\begin{itemize}

\item
$\G$ is connected, simply connected
 and semisimple  acting essentially with principal isotropy group
$\H$ of corank two.
\end{itemize}

In this section we consider the case where $\rank \G = 2$, and hence
$\H$ is finite. Clearly then
      $\Kpmo = \S^1$ or  $\S^3$.

\smallskip

   We will first deal
       with the most interesting case, where $\G$ is not simple, i.e.,
$\G=\S^3\times \S^3$.

       \begin{prop}\label{s3s3}
If $\G=\S^3\times \S^3$ acts essentially with corank 2, $M$ is
equivariantly diffeomorphic to one of the following spaces: An
Eschenburg space $E_p, p\ge1$, a $ P_k, k\ge 1$, the Berger space
$B^7$, a $ Q_k, k\ge 1$, or $R$ with the actions described in
\text{Table \ref{examples}}.
  \end{prop}

Since our actions are not assumed to be effective, we will use the
notation $\bar{\G} , \bar{\K}$ and $ \bar{\H}$  if the action is
made effective.
In view of our description provided in Table \ref{examples} in
Section \ref{sec: examples}, the Proposition is easily seen to
follow from the following:

\begin{lem}\label{s3s3pos}
Under the condition of the above Proposition, there are three
possibilities:
\begin {enumerate}
\item $\Hb=1$,  $\Kbm\cong \S^3$ and $\Kbp \cong \S^1$. In $\S^3\times
\S^3$, $\Km = \triangle \S^3 \cdot \H$, $\Kp=
\gC^i_{(p,p+1)}$ with
$p\ge 1$, and $\H \cong \Z_2$.

\item $\Hb\cong\Z_2$, $\Kbm\cong\SO(2)\, and \,\Kbp\cong\O(2)$. In $\S^3\times
\S^3$, the groups are   $\Km=\gC^i_{(1,1)}\cdot\H$ , $\Kp=
\gC^j_{(p,p+1)}\cdot\H$ with $p\ge 1$ and $\H\cong\Z_4\oplus\Z_2$, or
the same kind of  space with slopes $\{(3,1),(1,2)\}$.

\item $\Hb\cong\Z_2\oplus\Z_2$, and  $\Kbm\cong\O(2)\cong\Kbp$. In $\S^3\times
\S^3$, the groups are  $\Km=\gC^i_{(1,1)}\cdot\H$ , $\Kp=
\gC^j_{(p,p+2)}\cdot\H$ with $p$ odd $\ge 1$ and $\H\cong\Q$, or the same
kind of space with slopes $\{(3,1),(1,3)\}$.

\end{enumerate}

\end{lem}

\begin{proof}

If $l_- = l_+=3$, the assumption that the action is essential means
that $\Ko$ cannot be one of the $\S^3$ factors. Hence  both $\Kpmo
\simeq \S^3$ are embedded diagonally in $\S^3\times \S^3$,
contradicting group primitivity since any two diagonal embeddings
are conjugate, and in the effective picture all groups are
connected, and in particular $ \bar{\H} = \{1\}$.

We now know that at least one of the singular orbits has codimension 2, which
for the moment we denote as $\G/\K$ and where we can assume that, up
to conjugacy, $\Ko= \gC^i_{(p,q)}$ for two relatively prime
nonnegative
 integers $p,q$.
 Moreover, note that the Product Lemma \ref{product} implies that
neither $p$ nor $q$ can be 0 since the normalizer of $\Ko$ in one of
the $\S^3$ factors is finite.

 In the following we will
      make use of a consequence of the equivariance of the
second fundamental form of $\G/\K$ regarded as a $\K$ equivariant linear map $B\colon S^2T\to
T^\perp$. The non-trivial irreducible representations
       of $\S^1=\{e^{i\gt } \mid \gt \in \R \}$ consist of two dimensional
      representations given by multiplication by $e^{in\gt}$ on $\C$,
      called a weight $n$ representation. The action of $\Ko$ on
$T^\perp=\R^2$
      will have weight $k$ if $\H\cap \Ko=\Z_k$ since $\Z_k$ is  the
      ineffective kernel. As we will show below, only the cases $k=2,4$ arise
      and we claim that
 $|p- q|=2$  or $(p,q)=(1,1)$ in the case of $k=4$,   and
$|p- q|=1$ in the case $k=2$.

\no To see this, we first observe that the action of $\Ko$ on $T$
has weights $0$ on $W_{\subo}$ spanned by $(-qi,pi)$, weight $2p$ on the
two plane $W_1$ spanned by $(j,0)$ and $(k,0)$ and weight $2q$ on
the two plane $W_2$ spanned by $(0,j)$ and $(0,k)$. The action on
$S^2(W_1\oplus W_2)$  has therefore weights $0$ and $4p$ on
$S^2W_1$, $0$ and $4q$ on $S^2W_2$  and $2p+2q$ and $2p-2q$ on
$W_1\otimes W_2$.

\no Next, we claim that for any homogeneous metric on $\G/\Ko$ there
exists a vector $w_1\in W_1$ and $w_2\in W_2$ such that the 2-plane
spanned by $w_1$ and $w_2$ tangent to $\G/\K$ has curvature 0
intrinsically. Indeed, if $(p,q)\neq (1,1)$ or equivalently $p\neq
q$, $\Ad(\Ko)$ invariance of the metric on $\G/\Ko$ implies that the
two planes span$\{(j,0) , (0,j)\}$ and span$\{(k,0) , (0,k)\}$ and the line
$W_{\subo}$  are orthogonal to each other.  Hence $\Ad((j,j))$ induces an
isometry on $\G/\Ko$,
 which implies that
 the two plane spanned by the commuting vectors
  $w_1=(j,0)\in W_1$ and $w_2=(0,j)\in W_2$
  is the tangent space of
 the fixed point set of $\Ad((j,j))$ and thus
has curvature 0. If $(p,q)=(1,1)$, $\Ad(\Ko)$ invariance implies
that the inner products between $W_1$ and $W_2$ are given by
$\langle (X,0),(0,Y)\rangle = \langle \phi(X),Y\rangle$ where
$\phi\colon W_1\to W_2$ is an $\Ad(\Ko)$ equivariant map.
Hence, if we choose $j'=\phi(j)$ and $k'=\phi(k)$, the two planes
span$\{(j,0) , (0,j')\}$ and span$\{(k,0) , (0,k')\}$ are orthogonal to each
other, so that by the same argument $w_1=(j,0)\in W_1$ and
$w_2=(0,j')\in W_2$ span a 2-plane with curvature 0.

\no  If we now assume that $(p,q)\neq (1,1)$
 at least one of the numbers
$4p$ or $ 4q$ is not equal to the normal weight $k>0$. The
 equivariance of the second fundamental
form then implies that  $B_{S^2W_i}$ vanishes for at least one $i$
and hence by the Gauss equations  $B(w_1,w_2)\neq 0$ for the above
vectors $w_1$ and $w_2$. If $(p,q)= (1,1)$ the same holds if $k=2$.
Using the equivariance of the second fundamental form once more we
see that $W_1\otimes W_2$ contains a subrepresentation whose weight
is equal to the normal weight $k$. Hence, $|2p+2q|=k$ or
$|2p-2q|=k$, which proves our claim.

\smallskip

In addition we  observe that  $\H$  cannot contain an element $h$ of
the form
      $(a,\pm 1)$ or $(\pm 1,a)$ with $a$ being a noncentral element.
 Indeed, this would
      imply that
      $\N(h)_{\subo}=\S^1\times \S^3$ or $\S^3\times \S^1$ and hence $M^h$ would be a
      totally geodesic submanifold of
      codimension 2 in $M$.  By \eqref{codim2}
        $M$ would be $\Sph^7$ with a linear
action.
       But there is only one
       action on $\Sph^7$  with $\Kpmo=\S^1$, see Table \ref{odd}, and
      for that action $\H$ does indeed not contain such elements (cf. Table \ref{examples}).

\smallskip

Now let us consider the case where say $(l_-,l_+) = (3,1)$. Since
the action is assumed essential we have $\Kmo=\triangle  \S^3$ and
$\Kpo = \S^1$. From the fact that $\triangle  \S^3$ can be extended
only by the central element $(1,-1)$, we see that $\Kbm$ is
connected and $\bar{\H}=1$. Thus $\H=\Z_2$ since $\H=1$, and hence
$k=1$, contradicts the above equivariance argument.
 Thus $\H=\{(1,\pm 1)\}$ or
$\{(\pm 1,1)\}$, and $\Kp \supset \H$ is connected since $M$ is
simply connected (cf. \eqref{Hcomp}). We can assume that, up to conjugacy and switching the two factors
in $\S^3\times \S^3$,
 $\Kp=\Kpo=(e^{ip\gt},e^{iq\gt})$
for two relatively prime positive
 integers $p,q$ such that $q\ge p$.
 Using $k=2$,
      the above equivariance argument implies that
      $q-p=1$ and hence $(p,q)=(p,p+1)$ with $p>0$.

\smallskip

     It remains to consider the cases where $(l_-,l_+) = (1,1)$,
i.e., $\Kpmo=\S^1$. By Lemma~\ref{both-one} $\bgH$ contains only
elements of order two, which implies that $\H$ can only contain
elements of order two or four.
 This in turn implies that the normal
weights of the two singular orbits are $2$ or
$4$.

 We now have slopes $p_-,q_-$ on the left and $p_+,q_+$
on the right. We next proceed to derive the following
strong restrictions:  $1=\min\{|q_+|,|q_-|\}=\min\{|p_+|,|p_-|\}$.
The first step utilizes the Alexandrov geometry of the
quotients $M/\S^3\times 1$ and $M/1\times
\S^3$.

      In general, for an isometric $\G$ action on $M$, it is a
consequence of the slice theorem, that the \emph{strata}, i.e.,
components in $M/\G$ of orbits of the same type are (locally)
totally geodesic (cf. \cite{grove:survey}). In the case of
$M/\S^3\times 1$, the isotropy groups are effectively trivial on the
regular part since $(a,1)$ cannot lie in $\H$ unless it lies in the
center. Along $B_\pm$ the isotropy groups are $\Z_{q_-}$ and
$\Z_{q_+}$. This implies that the image of both $B_\pm$  in
$M/\S^3\times 1$ are totally geodesic if $\min\{|q_+|,|q_-|\}>2$.
Since these strata are two dimensional and $M/\S^3$ is four
dimensional, both strata cannot be totally geodesic according to
Petrunin's analogue \cite{Pe} of Frankel's theorem for Alexandrov
spaces.
 Hence
we have, $\min\{|q_+|,|q_-|\}\le 2$  and $\min\{|p_+|,|p_-|\}\le 2$.
Furthermore, if equality holds in one of these
inequalities, then
$\G$ acts effectively as $\SO(3)\times \gS^3$.

\smallskip

According to Lemma \ref{both-one}, two  cases remain
corresponding to $\bar{\H}= \Z_2$ or
      $\Z_2\oplus\Z_2$ since $\bar{\H}= 1$  and $l_\pm =1$ contradicts
 group primitivity. In
either case $\H$ contains an element $h$ of order four. Combined
with the above restrictions on $h$, we have $h^2=(-1,-1)$. Thus
$\bar{\G}\neq  \SO(3)\times \gS^3$ and
$1=\min\{|q_+|,|q_-|\}=\min\{|p_+|,|p_-|\}$ as claimed above.

\smallskip

      If $\bar{\H}=\Z_2$, we can assume that
      $\Kbm=\SO(2)\, ,\, \Kbp=\O(2)$ and the non-trivial element
      $\bar{h} \in \bar{\H}$ is in the second component of
$\Kbp$. Clearly, $\H$ contains an
      element $h$, whose image in
      $\bar{\H}$ is $\bar{h}$, and by the above each component in
      $h$
      is an unit imaginary quaternion. Since $\bar{h}$ acts trivially on
      $\Sph^{l_-}$
and by reflection on $\Sph^{l_+}$, so does $h$. In particular, $h$
commutes with $\Kmo$ and  we can
  arrange w.l.o.g. that
$\Kmo=\gC^i_{(p_-,q_-)}$ for two relatively prime positive
 integers $p_-,q_-$ with $q_-\ge p_-$. Then $h$ is one of $(i,\pm i)$, and
hence $p_-,q_-$ are both odd. Also, since conjugation by $h$ must
preserve $\Kpo$ and induce a reflection on it, we can assume, after
possibly conjugating with an element in $\N(h)$, that
$\Kpo=\gC^j_{(p_+,q_+)}$ with positive integers $p_+$ and $q_+$
which are relatively prime.

For the precise group picture in $\S^3\times \S^3$, there are two
possible subcases. Either $\H=\Z_4= \langle h \rangle = \{\pm (1,1) ,
\pm h \}$ or
$\H=
\Z_4\oplus\Z_2= \langle h, (1,-1)\rangle = \{ (\pm 1,\pm 1) , (\pm
i,\pm i)\}$. To rule out
$\H=\Z_4$, assume first that $p_+$ and $q_+$ are both odd. In this
case $\H\cap \Kpo=\Z_2$. Thus the normal weight is $2$ and
 equivariance implies that $|p_+ \pm q_+|=1$,
 a contradiction. If one is even and the other odd
      $\H\cap \Kpo=1$,   which contradicts again the above
equivariance argument.
 Now assume that $\H= \Z_4\oplus\Z_2=\{ (\pm 1,\pm 1)
, (\pm i,\pm i)$, which implies that $\H\cap \Kpo=\Z_2$ and hence
 $ q_+ - p_+  =  \pm 1$. On the left, we have that
 $\Kmo\cap \H =\langle h \rangle =\Z_4$
and hence the normal weight is $4$, which implies that $ q_- - p_- =
2$, or $(p_-,q_-)=(1,1)$.
Together with the above Frankel argument, this implies that we have
the possibility $(p_-,q_-)=(1,1)$ and $q_+-p_+ = \pm 1$ or
$(p_-,q_-)=(1,3)$ and $(p_+,q_+) =(2,1)$. In the first case we can
also assume that $q_+
>p_+$ by interchanging the two factors if necessary, and hence
$(p_+,q_+)=(p,p+1)$ with $p$  $\ge 1$.

\smallskip

Finally, we assume that $\bar{\H}=\Z_2\oplus\Z_2$. In this case
there are up to sign
 two noncentral order $4$ elements $h_-$ and $h_+$ in $\H$,
whose images $\bar{h}_-$ and $\bar{h}_+$ in $\bar{\H}$ are in the
second components of $\Kbp$ and of $\Kbm$ respectively, as well as
in the identity components $\Kmo$ and $\Kpo$ respectively. Notice
that $h_-$ and $h_+$ must anticommute in $\G$ since both components
of $h_-$ and $h_+$ as well as $h_-h_+$ are unit imaginary
quaternions. Since $\bar{h}_{\pm}$ act on $\Sph^{l_{\pm}}$ as
expected from the previous case, we can arrange that $\Kpmo$ are of
the form $\Kmo=\gC^i_{(p_-,q_-)}$ and $\Kpo= \gC^j_{(p_+,q_+)}$
respectively, and correspondingly $h_- = (\pm i,\pm i)$ and $h_+ =
(\pm j,\pm j)$ and thus all $p_i,q_i$ are odd. We can also arrange,
as above, that $q_-\ge p_->0$ and $p_+ , q_+>0$.

There are now two possibilities for $\H$. Either $\H=\triangle \Q$
(up to signs of the components) or $\H=\triangle \Q \oplus
\langle(1,-1)\rangle$. In the latter case, since $(1,-1)$ generates
another component for $\Km$ and for $\Kp$, $M$ is not simply
connected by Lemma \ref{isocomp}. Thus $\H=\triangle \Q$, the
weights on both normal spaces are $4$ and hence $q_\pm - p_\pm=\pm
2$  or $(p_\pm,q_\pm)=(1,1)$.
 Combining all of the above now yields only two
possibilities. Either $\{ (p_-,q_-)\, , \, (p_+,q_+)\} = \{(1,3) ,
(3,1)\} $ or $\{(1,1) , (p_+,q_+)\}$ with  $q_+ -p_+=2$, where we
used the fact that $\{ (p_-,q_-)\, , \, (p_+,q_+)\} = \{(1,1) ,
(1,1)\} $ would not be group primitive.
\end{proof}

We now turn to the simple rank two groups:

      \begin{prop}
There are no  actions of corank two of any of the groups $\SU(3)$,
$\Sp(2)$ or $\G_2$.
       \end{prop}

\begin{proof}

  From the Core-Weyl Lemma, we see that for the effective versions
$\Hb \ne 1$. In particular,  \eqref{Hcomp} implies that $l_{\pm}$
cannot both be $3$.

Now suppose one of $l_{\pm}$ is $3$, and w.l.o.g. then $\Kbm=\S^1$,
and $\Kbp=\S^3\cdot \Hb$ and hence $\Hb$ is cyclic by\eqref{Hcomp}.
It follows that  $\N(\H) \cap \Kpm$ are both at least 1-dimensional
and by part c) the  Upper  Weyl Group Bound
 $|\gW| \le 4$. But this
yields a contradiction to the Lower Weyl Group Bound if $\G=\Sp(2)$,
or $\G_2$. If $\G=\SU(3)$, then
$\N(\H)_{\subo}=\U(2)$ or $\T^2$. In either case it
follows that $w_+$ may be represented by  a central element in
$\N(\H)_{\subo}$. Using $\S^1=\Kbm\subset
\N(\H)_{\subo}$ it follows that the Weyl group
normalizes $\gK_-$. But then linear primitivity implies that
equality can not hold in the lower Weyl group bound -- a
contradiction.

It remains to consider the situation where both $l_{\pm} = 1$, and
thus,  by Lemma~\eqref{both-one}, either $\Hb=\Z_2$ or $\Z_2\oplus\Z_2$.
 In
the latter case we know that $\N(\H)/\H$ must be finite since each
of $(\N(\H) \cap \Kpm)/\H$ are and $M^{\H}$ is primitive. However,
for $\bar{\G} = \SO(5)$ we can diagonalize both involutions
simultaneously. In one case, then $\Z_2\oplus\Z_2$ is contained in
an $\SO(3)$ block and the normalizer contains a circle. In the other
case, $\Z_2\oplus\Z_2$ is contained in an $\SO(4)$ block, and the
normalizer contains a torus. Similar arguments can be applied to all
the other groups individually as well. These, however, are also all
covered by the a general result due to Borel \cite{borel}, which
asserts in particular that any $\Z_2\oplus\Z_2 \subset \bar{\G}$ is
contained in a torus unless $\pi_1(\bar{\G})$ has $2$-torsion.

If  $ \Hb=\Z_2$ and hence $\Kbm=\S^1 , \Kbp=\O(2) $ the Lower Weyl
Group Bound implies that $|W|\ge \dim \G/\H=\dim \G$ and  $|W|\le 8$
by the Upper Weyl Group Bound. Hence $\G=\SU(3)$,
and it follows that
      $\N(\H)=\U(2)$ since this is the only equal rank symmetric subgroup of
$\SU(3)$. In particular $N(\H)$ is connected
 and the Core-Weyl Lemma gives the contradiction
$|W|\le 4$.
\end{proof}


\section{Semisimple Rank $3$ Groups.} \label{sec: rank3}

If $\G$ has rank 3 and $\H$ has corank 2, one has the two subcases
$\Ho=\S^1$, or $\Ho$ is one of $\S^3$ or $\SO(3)$. Also recall that
$\H/\Ho$ is cyclic. By the Isotropy Lemma $\max\{l_-,l_+\}\ge 2$, and
by the rank Lemma $l_\pm$ cannot both be even.

In the case of $\Ho=\S^1$, one has the possibilities
$(\l_-,\l_+)=(1,2),(1,3),(2,3),(3,3)$ (up to order) and in the
latter two cases all groups are connected.  Furthermore,
     $\Ko=\T^2$ if $l_\pm =1$ ,  $\Ko=\SO(3)$, or $\S^3$ if $l_\pm=2$
      and  $\Ko=\U(2)$, or $\S^1\times\S^3$
     if $l_\pm=3$.

If $\Ho$ is 3-dimensional, one has the possibilities $\l_\pm=
1,3,5,7$ and $\Ko=\U(2)$, $\S^3\times\S^1$, or $\SO(3) \times\S^1$
if $l_\pm=1$, $\Ko=\SO(4)$, or $\S^3\times\S^3$ if $l_\pm=3$,
$\Ko=\SU(3)$ if $l_\pm=5$ and $\Sp(2)$ if $l_\pm=7$.  If
$\Ho=\SU(2)$ (in every effective version), the lowest dimension of a
representation is 4, which must degenerate somewhere and hence one
of $\Kpmo=\SU(3)$ or $\Sp(2)$.

\bigskip

We will first deal with the case where $\G$ has a normal subgroup of rank
one, i.e., almost effectively $\G = \S^3\times \L$, where $\rank \L = 2$.

\begin{prop}
     If $\rank \G =3$ and
$\G$ has a normal subgroup of rank one, an essential action of $\G$
with corank 2 is the tensor product action of $\SU(2)\SU(3)$.

       \end{prop}

\begin{proof} Before we start with the four possible subcases, let
us notice that a three dimensional subgroup $\Ho$ of $\S^3\times \L$
must be contained in $\L$ since the action is almost effective and
essential.

\smallskip

\begin{center}
Case 1. $\G=\S^3\times \S^3 \times \S^3$
\end{center}

     If $\Ho$ were three dimensional,
the projection onto one of the factors would be onto and hence the
action would be inessential. Thus $\Ho=\S^1$, and one of $l_\pm$, is
$2$, or $3$. First suppose, e.g., $l_- = 3$. Then the semisimple
part of $\Km$ is $\S^3$ whose involution is a Weyl group element.
Being central in $\G$,  it has $\G/\Km$ as a fixed point component,
contradicting the fact that it cannot have positive curvature. Hence
we are left with $l_-=2$ and $l_+=1$. In particular $\Kmo=\S^3$ and
$\Kp=\T^2$. By the Product Lemma it follows that we can assume that
$\Kmo=\{(q,q,q)\vert q\in \S^3\}$ and hence $\Ho=\{(z,z,z)\vert z\in
\gS^1\}$. Clearly then the cyclic group $\Km/\Kmo=\H/\Ho$  has at most
two elements. Since $\Kp \cong \T^2 \subset \N(\Ho)_{\subo} \cong \T^3$ we
can represent the Weyl group element $w_+$ by an element of the form
$\iota = (\iota_1,\iota_2,\iota_3)$ of order 2 if $\H=\Ho$, and
order 4 otherwise. Since we can also replace $\iota$ by
$\iota(i,i,i)$ we can arrange that $\iota_p^2 = 1$ holds for at least two
indices $p$. But then a component of $M^{w_+}$ is a totally geodesic
submanifold of $\G/\Kp$ of the form $\S^3\times \S^3\times
\S^3/\T^2$ or $\S^3\times \S^3\times \S^1/\T^2$ neither one of which
can have positive curvature.

\smallskip

\begin{center}
Case 2. $\G =\S^3\times \SU(3)$
\end{center}\smallskip

     We first settle the case that $\H$ is 3-dimensional. The only three
dimensional spherical  subgroup of $\SU(3)$ is $\SU(2)$ (cf. Table
\ref{sph} in Appendix II). Since its normalizer is $\S^3\times
\U(2)$, the action by $\Ho$ is fixed point homogeneous,  $M$
is a sphere, and the action is inessential.

Now suppose $\Ho=\S^1$. We can then assume
     that $\Ho$ is not contained in the $\S^3$ factor since otherwise $M$
     would again be fixed point homogeneous.
     We distinguish between
     two subcases:

a)  The involution $\iota\in \Ho$
     is not in the center of $\G$, i.e. $\iota=(\pm 1,b)$, and we can
assume $b =\diag(-1,-1,1)$.

b)  The involution of $\Ho$ is central in $\G$.

\smallskip

Subcase a). Then
     $\N(\iota)_{\subo}=\S^3\times \U(2)$ acts on
$M^\iota_c$ by \co one with one
     dimensional principal isotropy group. Thus $M^\iota_c$ has dimension 7 and
     $M$ dimension 11 and hence $M^\iota_c$ is simply connected by the
     Connectedness Lemma.

Let us first assume that  $M^\iota_c$ is a sphere. The Connectedness
Lemma implies  that $M$ is $4$-connected.
     We may assume that the action of $\S^3 \times \U(2)$ on $M^\iota_c$ has
      finite kernel, since otherwise we can deduce from part (b) of
     the Connectedness Lemma that $M$ is $5-$connected and hence a sphere.
     By assumption $\N(\iota)_{\subo}$ acts linearly on
     $M^\iota_c$. There are two types of
     linear actions by $\S^3 \times \U(2)$ on the
     7-sphere: one is a sum action and the other
the tensor product
     action. If it were a sum action,
     the $\S^3$ factor would have a fixed point and hence would be
     contained in some $\Kpm$, contradicting the assumption that the
     action on $M$ is essential.

     Hence it is the tensor product action and thus $\S^3$
     acts freely on $M^\iota_c$. This implies that the action of
$\S^3$ on
     $M$ is also free since all $\G$ orbits meet $M^\iota_c$ and
$\S^3$ is
     normal in $\G$. Since $M$ is 4-connected, the quotient $M/\SU(2)$ is
     two connected but not 4-connected and by Verdiani's classification
     in even dimensions $M/\SU(2) = \HP^2$. From the Gysin sequence it
     follows first
     that the Euler class of the
      bundle
     $\S^3\to M\to \HP^2$ is a generator of $H^4(\HP^2,\Z)$
     (again since $M$ is 4-connected),
     and then that  $M$ is a homology sphere. From Table \ref{odd} we then that it must be  the tensor product action of $\SU(2)\SU(3)$.

Next we exclude the case  that $M^\iota_c$ is not a sphere. Since
any two involutions in $\SU(3)$ are conjugate, we can choose an
element
     $g\in \SU(3)$ such that $\iota$ and
      $g\iota g^{-1}$ span a dihedral group
      $\D_2=\Z_2^2$. By Frankel,
$M^\iota_c\cap g M^\iota_c$ is non-empty and by transversality at
least 3-dimensional. Since $\D_2$ is contained in a torus the
codimension is even.   From the assumption that $M^\iota_c$ is not a
sphere, we conclude that it cannot have dimension 5 by
\eqref{codim2}, and hence it is 3-dimensional.  Since $M^\iota_c\cap
g M^\iota_c\to M^\iota_c$ is 3-connected by part (c) of the
Connectedness Lemma, $M^\iota_c\cap g M^\iota_c$ is simply connected
and hence must be $\Sph^3$. In particular $M^\iota_c$ is
2-connected. The only 2-connected positively curved 7-manifolds in
our classification theorem are $B^7$ and $P_k$. However, as we have
seen in Lemma \ref{extensions}, for these manifolds the group does
not have a connected normal extension. It follows that the $\S^3
\times \U(2)$ action has a one dimensional kernel, which must be the
center of $\U(2)$, and hence this is actually an action by
$\SU(2)\SO(3)$. But this group does not act on $B^7$ or $P_k$ or any
of its subcovers, see Table \ref{examples}.

\smallskip

Subcase b). In this case $\Ho$ has only one involution, namely
$(-1,\diag(1,1,1))$.

Consider the cyclic subgroup $\gC_4$ of order four in $\gH_{\subo}$. We
may assume $\gC_4\not\subset \gS^3$
 and thus $\N(\gC_4)= \Pin(2)\times \U(2)\supset \N(\gH)$.
 Let $M'$ be a component of $\Fix(\gC_4)$ on which
 $\N(\gC_4)_{\subo}$ acts with cohomogeneity one.
 By induction assumption $M'$ is up to covering
 a 5-sphere endowed with a linear action.
 This shows that $\gK_-$ ( or $\gK_+$) is a
 $4$-dimensional subgroup of $\N(\gC_4)_{\subo}$.

Clearly the semisimple part $\SU(2)$ of $\gK_-$ is normal in
 $\N(\gC_4)\supset \N(\gH)$ and  $\SU(2)\cdot \gH=\gK_-$.
 Hence $\N(\gH)$ and thereby the Weyl group
 normalizes $\gK_-$. Because of
 $\rank(\gK_-)=2$ it is clear that $\N(\gH)_{\subo}\not \subset \gK_-$.
 Combining this with linear primitivity we see that
 $\gK_+/\gH$ contains a trivial subrepresentation.
 Therefore $\gK_+/\gH\cong \Sph^3$, or $ \Sph^1$.
 The latter case would imply $\gK_\pm\subset \N(\gC_4)$ which
 contradicts primitivity.
 In the former case the Weyl group has order
 at most $4$, by the upper Weyl group bound, Proposition~\ref{upper}.
 Since $\gK_-$ is normalized by the Weyl group, linear primitivity says
 that the Lie algebras of $\gK_-$, $\gK_+$ and $w_-\gK_+ w_-$
 span the Lie algebra of $\G$. But this is clearly impossible as these groups
 have $\gH$ in common.

\smallskip

\begin{center}
Case 3. $\G=\S^3\times \Sp(2)$
\end{center}

     Again we first settle the case that $\H$ is 3-dimensional.
There are two spherical 3 dimensional subgroups of $\Sp(2)$ :
$\Sp(1)\times 1$ \, and \, $\triangle \Sp(1)$ (cf. Table \ref{sph}).
In the first case $\Ho$ acts transitively in the unit sphere
orthogonal to $M^{\Ho}$ since $\N(\Ho)=\S^3\times \Sp(1)\times
\Sp(1)$ and is hence fixed point homogeneous. In the second case
$\G/\H$ effectively becomes $\S^3\times\SO(5)/\SO(3)$ and the Chain
Theorem applies.

We can now assume $\Ho=\S^1$ and one, say $\Km$ has rank 2. If
$\Kmo$ contains one of the involutions $\iota = (\pm 1 ,\pm
\diag(1,-1))$, up to conjugation, we obtain a contradiction as
follows. If $\iota$ lies in $\H$, $M^\iota_c$ is \co one under
$\N(\iota)_{\subo}=(\S^3)^3$ with one dimensional principal isotropy
group. As we saw in Case 1, such an action does not exist. If
$\iota$ does not lie in $\H$, it has $(\S^3)^3/(\Km\cap \N(\iota))$
as a fixed point component, which cannot have positive curvature.

We may assume that $\Km$ contains the center of $\S^3\times\Sp(2)$.
Since $\G/\Km$ cannot be totally geodesic, it follows that the
center of $\S^3\times\Sp(2)$ is contained in $\H$. Therefore $\H$ is
not connected and we may assume that $\Km \cong \T^2$ (see
Lemma~\ref{Hcomp}).
 By the product
lemma $\Km$ projects to a maximal torus of $\Symp(2)$. Since $\gH$
contains no involution as above, it follows that the Weyl group
element $w_-$
 can be represented by an element
$\iota:=(*,\diag(\pm 1,  \pm 1))\in \Km$. Clearly the fixed point
set of $\iota$ would be a homogeneous space which does not have
positive sectional curvature.

\smallskip

\begin{center}
Case 4. $\G=\S^3\times \G_2$
\end{center}

     We first rule out the case that $\H$ is 3-dimensional.
     The only 3-dimensional spherical  subgroup of $\G_2$
      is $\SU(2)\subset\SU(3)\subset\G_2$. Although this does not
      immediately follow from Table \ref{sph}, it is easily
      verified by considering the four three dimensional subgroups
      of $\G_2$. Since a four dimensional representation of $\Ho=\SU(2)$ must
degenerate, one of $\Kpmo=\SU(3)\subset \G_2$ (no $\Sp(2)$ exists in
$\G_2$), which contradicts the Product Lemma.

Hence $\Ho=\S^1$ and we can assume that $\rank \Km=2$. Among the
involutions in $\Kmo$ there is one of the form $\iota=(\pm 1, b)$
with $b$ an non-trivial involution, which has normalizer $\SO(4)$
(see Table \ref{symm}). Thus $\N(\iota)_{\subo}=\S^3\times \SO(4)$. If
$\iota$ lies in a principal isotropy group, the reduction
$M^{\iota}_c$ has $\S^3\times \SO(4)$ acting by \co  one with a one
dimensional principal isotropy group, but such an action does not exist as we
saw in the first case. Otherwise $\iota$ has a homogeneous fixed
point component $\S^3\times\SO(4)/(\Km\cap \N(\iota)_{\subo})$ which
cannot have positive curvature.
\end{proof}

\smallskip

It remains to deal with the cases where $\G$ is simple.

\begin{prop}\label{rank3}
If $\G$ is  simple with  $\rank \G =3$ acting essentially and with
corank 2, then it either the linear reducible representation of  $
\SU(4)$ on $\Sph^{13}$ or the \co one action of $\SU(4)$ on one of the
Bazaikin spaces $B^{13}_p, p \ge 1$ \emph{(see Table
\ref{examples})}.
       \end{prop}

\begin{proof}
There are three cases to consider, corresponding to $\G=\SU(4),\Sp(3)$ or
$\Spin(7)$. We first consider the most interesting case where
$\G=\SU(4)$.

\begin{center}
Case 1. $\G=\SU(4)$
\end{center}

We will first  rule out the case  that $\Ho=\S^1.$ We can assume
that $\Ho=\diag(z^{p_1}, z^{p_2},z^{p_3},z^{p_4})\subset \SU(4)$ and
hence the isotropy representation of $\G/\H$ has weights $p_i-p_j$.
By the Isotropy Lemma there can be at most two distinct non-zero
weights and one easily sees that this leaves only four possibilities
$(p_1,p_2,p_3,p_4)=(1,-1,0,0), (1,1,-1,-1),(1,1,1,-3),$ and
$(3,3,-1,-5)$. In the last two cases $\N(a)=\U(3)$ for some element
$a\in \Ho$ corresponding to $z$ with $z^8 = 1$. But then  the
reduction $M^a_c$ is a \com manifold under $\U(3)$ with  one
dimensional principal isotropy group, which does not exist by
induction.

   If $(p_1,p_2,p_3,p_4)=(1,-1,0,0)$ we choose the
involution $\iota=\diag(-1,-1,1,1)\in \Ho $. Then
$\N(\iota)_{\subo}/\iota = \S(\U(2)\U(2))/\iota=\SO(3)\U(2)$ acts by \co
one on the seven dimensional reduction $M^\iota_c$ with one
dimensional principal isotropy group.  By induction, up to covers,
such a 7-dimensional \com  could be only  a sphere with  a sum
action
     or the Eschenburg space $E_p$. But in both cases, the isotropy group is
     not contained in the $\SO(3)$ factor as it is for $M^{\iota}_c$.

If $(p_1,p_2,p_3,p_4)=(1,1,-1,-1)$, we   observe that
$\N(\Ho)_{\subo}/\Ho=\S(\U(2)\U(2))/\diag(z,z,\bar{z},\bar{z})$ is equal
to $\SO(4)$ since $\SU(2)\SU(2))$ acts transitively with isotropy
$\diag(-1,-1,-1,-1)$. In the full normalizer $\N(\Ho)/\Ho$ we have a
second component corresponding to the element that interchanges the
two normal $\SU(2)$ subgroups of $\S(\U(2)\U(2))$. Hence
$\N(\Ho)/\Ho=\O(4)$. Furthermore, $M^{\Ho}$ has only one seven
dimensional component  since the
inclusion $M^{\Ho}\subset M$ is 1-connected by part (b) of the
Connectedness Lemma. Hence $\O(4)$ acts by \co one on $M^{\Ho}$ with
cyclic principal isotropy group. Such a manifold is either $Q_k$ or
a spaceform. But in  $Q_k$ the slopes of $\Kp=\S^1$ are $(k,k+1)$
and hence its (ineffective) $\SO(4)$ action does not extend to
$\O(4)$. It also cannot be a space form, since the action on its
cover would be a sum or modified sum action and hence $|\W|\le 4$,
which gives a contradiction to the lower Weyl group bound $l_-+l_+\ge 7$.

\vspace{10pt}

We can now assume that $\Ho$ is three dimensional. But the only
spherical 3-dimensional subgroups of $\SU(4)$ are $\SU(2)\subset
\SU(3)\subset \SU(4)$ or $\Delta\SU(2)\subset \SU(2)\SU(2)\subset
\SU(4)$ , (cf. Table \ref{sph}). In the latter case
$\G/\Ho=\SO(6)/\SO(3)$ and the Chain Theorem applies.

Hence we can assume $\Ho=\SU(2)$ embedded as the lower $2\times 2$
block. By the isotropy lemma
one of $\Kpmo$ is equal to $\SU(3)$ or $\Sp(2)$. It is important to
observe that $\N(\Ho)/\Ho=\U(2)$ acts transitively on all possible
embeddings of $\SU(3)$ or $\Sp(2)$ in $\SU(4)$ containing the same
$\Ho$ (In the case of $\Sp(2)$ this is best seen in $\SO(6)$).

Assume first that $\Kmo=\SU(3)$. If $l_+=1$, $\Kp=\SU(2)\cdot\S^1$
is connected and, modulo $\N(\Ho)/\Ho$, both $\Kpm$ are contained in
$\U(3)$ which contradicts primitivity. If $l_+=3$ and hence
$\Kp=\SU(2)\SU(2)$, the element $-\Id\in \SU(4)$ is in $\Kp$ and
represents a Weyl group element. Since it is central, $B_+$ is
totally geodesic, but it cannot have positive curvature. If $l_+=5$
and hence $\Kp=\SU(3)$, the action is not primitive. If $l_+=7$ we
have $\Kp=\Sp(2)$. All embeddings are determined, modulo
$\N(\Ho)/\Ho$, and we have the linear action of $\SU(4)$ on
$\Sph^{13}$.

This leaves $\Kmo=\Sp(2)$. If also $\Kp=\Sp(2)$,  the action is not
primitive. The case of  $\Kp=\SU(2)\SU(2)$ is dealt with as above,
and $\Kp=\SU(3)$ was already considered. It only remains to consider
the case where $l_+=1$. Since $\Kp = \Kpo \subset \S(\U(2)\U(2))$ we
can assume up to conjugacy that $\Kp =
\Ho\diag(z^k,z^l,\bar{z}^{(k+l)/2},\bar{z}^{(k+l)/2})$. Notice that
 $-\Id\in \SU(4)$ cannot be in $\H$ since it is in $\Kmo=\Sp(2)$
and $\Km/\H=\Sph^7$. But $-\Id$ can also not be in $\Kp$, since then
it would represent $w_+$ in contradiction to the fact that $B_+$ has
zero curvatures and hence cannot be totally geodesic. This implies
that $(k,l)=(2p,2q)$ with $(p,q)=1$, $p$ and $q$ not both odd.
Choosing $z = i$ and multiplying by $\diag(1,1,\pm(i,-i))$ we see
that $\iota=\diag(1,-1,-1,1)$ or $\iota=\diag(-1,1,-1,1)$ is in
$\Kp$. If it does not lie in $\H$, it has
  $\U(2)\SU(2)/\T^2\subset
\Bp$ as a  fixed point component,  which does not have positive
curvature. Hence $\H$ is not connected. Since $\Sp(2)\subset\SU(4)$
can only be extended by $\Z_2$, $\H/\Ho=\Z_2$, with $\iota$
representing a second component. Thus $M^\iota_c$ is \co one under
the action of $\S(\U(2)\U(2))/\langle \iota \rangle$ with $\N^{\H}
(\iota)/\langle \iota \rangle =\diag(1,1,z,\bar{z})$ as its
principal isotropy group. Moreover, the subaction by
$\SU(2)\SU(2)/\langle \iota \rangle = \SO(3)\SU(2)$ is again \co one
with trivial principal isotropy group. In this reduction
$\Kp=(z^{2p},z^{2q})$  which effectively becomes $(z^p,z^q)$. This
reduction must be an Eschenburg space, and hence $(p,q)=(q+1,q)$
with $q\ge 1$. Hence our original manifold must be a Bazaikin space
$B^{13}_p$ (cf. Table \ref{examples}).

\begin{center}
Case 2. $\G=\Sp(3)$
\end{center}

   The symmetric subgroups of $\Sp(3)$ are
$\Sp(2)\Sp(1)$ and $\U(3)$ where the latter is only a normalizer
under an order 4 element (e.g. $i \Id$). If $\Ho=\S^1$ we  have, for
appropriate $a$, $\N(a)=\Sp(2)\Sp(1)$ or $\U(3)$ with one
dimensional principal isotropy group which does not exist by
induction.

 Now assume that $\Ho$ is three
dimensional. The 3 dimensional spherical subgroups of $\Sp(3)$ are,
according to Table \ref{sph}, $\diag(q,q,q),\diag(q,q,1)$ or $
\diag(q,1,1)$ with $q\in\Sp(1)$. In the first case, we can choose
$\iota=i\text{Id}\in\Ho$ and hence $\N(\iota)=\U(3)$  acts by \co
one on $M^{\iota}_c$ with one dimensional principal isotropy group,
which does not exist by induction. In the second and third case we
can choose an involution $\iota\in\Ho$ with $\N(\iota)=\Sp(2)\Sp(1)$
which acts by \co one on the reduction $M^{\iota}_c$ with three
dimensional  principle isotropy group.
 By induction it
must be a linear sum or modified sum action which contains a
standard $\Sp(1)\subset \Sp(2)$ in its principal isotropy group.
Thus $\Ho=\diag(1,1,q)$ and hence $\Sp(2)$ acts with finite
principal isotropy group on the reduction $M_c^{\H_{\subo}}$, which, as we
saw in Section 6, is not possible.

\begin{center}
Case 3. $\G=\Spin(7)$
\end{center}
   The  symmetric subgroups of $\SO(7)$
are $\SO(6), \SO(5)\SO(2)$ and $\SO(4)\SO(3)$, and correspondingly
for $\Spin(7)$. If $\Ho=\S^1$, we can choose $\iota \in \Ho$ with
$\iota^2$ but not $\iota$ in the center of $\Spin(7)$, and
$\N(\iota)_{\subo}$ is one of the groups $\Spin(6), \Spin(5)\Spin(2)$ or
$\Spin(4)\Spin(3)$. Hence they act by \co one on the reduction
$M^\iota_c$ with one dimensional principle isotropy group. But such
a manifold does not exist by induction.

Now suppose $\H$ is $3$-dimensional. If $\Ho$ is  a $3\times 3$
block in $\Spin(7)$ we are done by the Chain Theorem. Thus by Table
\ref{sph} we can assume that $\Ho=\SU(2)$ is embedded as a normal
subgroup of a $4 \times 4 $ bock. By the Isotropy Lemma
  a four  dimensional representation of $\Ho$
must degenerate, which means that one of $\Kpmo$ must be $\SU(3)$ or
$\Sp(2)$. There is only one embedding of $\Sp(2)$ and, since it
corresponds to $\SO(5)\subset \SO(7)$, its central element is
central in $\Spin(7)$. It  then has   $\G/\K=\Spin(7)/\Sp(2)$ as its
fixed point set which does not admit positive curvature.

We can therefore assume that $\Kmo=\SU(3)$. Observe now that
 $\N(\Ho)_{\subo}/\Ho=(\Spin(4) \times
\Spin(3)/\triangle \Z_2) / \SU(2)=\S^3\times \S^3/(-1,-1)=\SO(4)$
acts by \co one on the reduction $M^{\Ho}_c$ with cyclic principal
isotropy group $\H/\Ho$. All non-spherical examples and candidates
in dimension $7$, as well as their subcovers, do either not admit a
cohomogeneity one action of $\SO(4)$, or only allow for actions with
a noncyclic principal isotropy group. Thus $M^{\gH}_c$ is a space
form. Using once more that the principal isotropy group is cyclic we
see that the action is inessential and thus both singular orbit have
codimension $4$, a contradiction as the left singular orbit has
codimension $2$.
\end{proof}


\section{Semisimple Groups with a Rank $1$ Normal
Subgroup.}\label{sec: rank one factor}

In this section we will complete the analysis of simply
connected, positively curved cohomogeneity one $\G$-manifolds, where
$\G$ has a normal
subgroup of rank one:

\begin{prop}
Suppose a semi-simple $\G$ of rank at least four has a normal
subgroup of rank one, and acts essentially with corank 2. Then $\G =
\SU(2)\SU(n) $, $ M = \Sph^{4n-1}$ and the action is the tensor
product action.
\end{prop}

\begin{proof}

      Let $\G=\S^3 \times \L$, where
$\L$ is a simply connected semisimple group with $\rank \L \ge 3$
and hence $\rank\H\ge 2$ and  $\rank \H\cap\L \ge 1$.

      First observe that if
$\H \cap \S^3 \triangleleft \H$ is not contained
      in the center of $\S^3$, then the reduction $M^{\H \cap \S^3}_c$ has
codimension 2 in $M$, and hence $M$ is a sphere, and we are done by the
classification of essential actions on spheres. Thus, if we set
$\S =
\S^3$ if $\H \cap \S^3$ is trivial, and $\S
=
\SO(3)$ if $\H \cap \S^3$ is non-trivial, we can assume
that $\G = \S \times \L$ and $\H \cap
\S$ is trivial.

In the proof we will use the following useful notation for the
groups $\Kpm$ and $\H$: $\K_{\S}=\K\cap \S $ and $\K_{\L}=\K\cap
\L$.
 Furthermore, there exists a connected normal subgroup $\K_\Delta
$ of $\Ko$  embedded diagonally in $\S\times\L$ such that $\Ko =
(\K_{\S}\cdot\K_\Delta\cdot\K_{\L})_{\subo}$. It follows that $\K_\Delta$
is a rank one group and, by  the Product Lemma, $\K_{\S}$ is finite,
if non-empty.

We divide the proof into three subcases: (1) $\S = \S^3$ acts
freely, (2) $\S = \SO(3)$ acts freely, and (3) $\S$ does not act
freely. As it turns out, only the first case can occur.

\smallskip

\smallskip

\begin{center}
Case 1.  $\S^3$ acts freely
\end{center}

In this case $B:=M/\S^3$ is an even dimensional simply connected
      cohomogeneity one $\L$ - manifold of positive curvature.
      By Verdiani's classification, $B$ is a rank one symmetric space
      and the action of $\L$ on $B$ is linear.

Fix a maximal torus $\T = \T^h$ of $\H_{\L}=\H\cap\L\subset \L$,
which has positive dimension by assumption, and consider the
reduction $M' = M^{\T}_c$. Since $\N(\T)/\T = \S^3 \times
\N^{\L}(\T)/\T$, the reduction $M'$ supports a cohomogeneity one
action by a group $\S^3\times \L'$,
      where $\L'$ has rank 1 if $\Ho\subset\L$, or rank 2 if $\H_\Delta$
      is non-trivial.  The group $\L'$  also acts on the reduction $B^{\T}_c$
      as well as on $M'/\S^3\subset B^{\T}_c$ and in both cases with
      principal isotropy group $\N^{\H}(\T)$. Hence  $B':=M'/\S^3= B^{\T}_c$

 The totally geodesic fixed point set $B'$ is again a rank one symmetric
      space and must be simply connected since it is orientable.
       This in turn implies that $M'$
      is simply connected.

Since $\T$ is a maximal torus in $\H_\L$, the principal isotropy
group of the
 $\S^3\times \L'$ action on $M'$
       has
      at most finite intersection with the $\L'$ factor. As the subaction
      of the $\S^3$-factor is free, our results in the previous two sections
       combined with
      Lemma \ref{free} implies that
      $M'=\Sph^{4k+3}$  and $B'=\HP^k$.

The Euler class of the $\S^3$ bundle $M\to B$ pulls back to the
Euler class of $M'\to B'$ which is a generator in $H^4(B',\Z)=\Z$.
This is only possible if
       $B\cong \HP^l$. The Euler class of $M\to B=\HP^l$ is therefore
       also a generator of  $H^4(\HP^l,\Z)$, and the
Gysin sequence
      implies that $M$ is a homology sphere. Table \ref{odd} now
      shows that it is the tensor product action of $\SU(2)\SU(n)$.

\smallskip

\smallskip

\begin{center}
Case 2.  $\SO(3)$ acts freely
\end{center}

In this case $B = M/\SO(3)$ is an even dimensional positively curved
cohomogeneity one $\L$-manifold. Since $M \rightarrow B$ is a
principal $\SO(3)$ bundle and $M$ is simply connected we see that
$B$ is simply connected, but not $2$-connected. By Verdiani's
classification $B$ is a complex projective space. In the long
homotopy sequence $\pi_2(M)\to \pi_2(B)\to\pi_1(\SO(3))=\Z_2\to
\pi_1(M)$ the map in the middle can be regarded as representing the
second Stiefel Whitney class in $H^2(B,\Z_2)$. Hence it is
non-trivial for the bundle $M\to B$.

Consider as above a maximal torus $\T =\T^h$ of $\H_{\L}$
      and the corresponding reductions $M' \subset M$ and $B' \subset B$.
       Since the $\L$ action on $B$ is linear, it follows that $B'$
      is a complex projective space as well, and by naturality, the principle
$\SO(3)$ bundle
      $M'\rightarrow B'$ has a non vanishing
      second Stiefel Whitney class also. This in turn implies that
       $M'$ is simply connected.

Also as above, we note that $M'$ comes with a cohomogeneity one action of
      $\SO(3)\times \L'$ where $\rank(\L')\in \{1,2\}$.
      Since $\SO(3)$ acts freely, it follows from our previous sections and
Lemma \ref{free}  that $M'=E_1\, , E_2$ with $\L'=\SU(2)$ or
$M'=W^7_{(2)}$ with $\L'=\SO(3)$. In all three cases
      $B'\cong \CP^2$ (see Remark \ref{CP2Bundle}) and the action of $\L'$ on
      $\CP^2$ is the action of $\SU(2)$ with a fixed point in the first two cases
       and in the third case
       the action of $\SO(3)$ on $\CP^2$  induced by
        the tensor product action of $\SO(2)\times \SO(3)$ on $\Sph^5$.

      Consider first the case that  $B'$ is endowed with the standard
$\SU(2)$
      cohomogeneity one action which has a fixed point. Clearly  only
another "sum" action on a higher dimensional complex projective space
can have this as
a reduction. Because of
$\rank(\L)\ge 3$,
      it follows  that a normal simple subgroup
$\L'\subset
\L$ of rank
      at least $2$ has non-empty fixed point set in $B$, and in fact acts fixed
point homogeneously. Since the action of $\SO(3)$ on the fibers only
extends to an action of $\SO(4)$ and the action of $\L'$ fixes one
$\SO(3)$ orbit in $M$, it follows that $M$ is fixed point
homogeneous. Clearly this is not possible since spheres do not
support free actions of $\SO(3)$.

Assume now that  $B'\cong \CP^2$ is equipped  with the
      cohomogeneity one action of $\SO(3)$ with
      both singular orbits of codimension two.
 The only way this is a reduction of an
$\L$-action on a higher dimensional complex projective space,   is
that up to orbit equivalence the $\L$ action  is given by an
      $\SO(h+1)$-action on $\CP^h$ for some $h\ge 5$. Indeed, one
      sees that for all other actions in Table \ref{even}, one of the
      normal spheres has odd codimension, which is preserved under a
      reduction by a torus.

      The codimension of the singular orbits of the
      $\SO(h+1)$-action are $2$ and $h-1$.
      The singular isotropy group for the orbit of codimension $h-1$ has a
simple identity
component
      of $\rank\ge 2$ and $\Km=\SO(2)\cdot \H$ (see Table \ref{even}). For
the lifted picture
      upstairs in $M$, i.e., in the diagram $\H \subset \{\Km, \Kp\} \subset
\SO(3)\times \L$, we see that the projections of $\Kp$ and $\H$
 to the
$\SO(3)$ factor are  trivial and the projection of $\Km$ is one
dimensional. But this contradicts
       group primitivity.

\smallskip

\begin{center}
Case 3.  $\S^3$ or $\SO(3)$ does not act freely.
\end{center}

In this subsection $\S$ is one of $\S^3$ or $\SO(3)$, and we assume
that $ \H_{\S} =\H\cap \S= 1$, but $\S$ does not act freely on $M$.
In particular one of $\K^\pm_{\S}$, say $\K^-_{\S}$ is non-trivial.

Choose  an element $\iota \in \K^-_{\S}$. Since $\iota$ is not in
$\H$, the component $V$ of $M^{\iota}$ containing $c(-1)$ is an odd
dimensional positively curved homogeneous space $\N(\iota)_{\subo}/\Km
\cap \N(\iota)_{\subo}$. From the classification of positively curved
homogeneous spaces we deduce that

\begin{itemize}
\item
$V = \L/ \K^-_{\L} $.
\end{itemize}

 Since
 $\Km \cap \N(\iota)_{\subo}$ has corank one in $ \N(\iota)_{\subo}$ and
$\rank \N(\iota)_{\subo} = \rank \S \times \L$, it follows that $\Km$ has
corank one in $\G = \S\times \L$. The Product Lemma hence implies
that $(\K^-_\Delta)_{\subo}$ is non-empty. Indeed, since $\S\times\L$ and
$\Km$ do not have a normal subgroup in common, we have either
$(\K^-_{\S})_{\subo}=\S^1$, which has finite normalizer in $\S$, or
$\Kmo=(\K^-_{\L})_{\subo}$ is of equal rank in $\L$ which has finite
normalizer in $\L$. Thus it also follows that the projection of
$\Km$ into $\L\subset \S\times\L$, which is isomorphic to
$\K^-_\Delta\cdot\K^-_{\L}$, has equal rank in $\L$ and hence
$\N^\L((\K^-_{\L})_{\subo})$  has equal rank also, i.e. $(\K^-_{\L})_{\subo}$ is
a regular subgroup of $\L$.

 The cover $\tilde V =
\L/(\K^-_{\L})_{\subo}$ of $V$
      is hence an odd dimensional homogeneous space of positive curvature
       with $\L$
semisimple
      of rank $\ge 3$ and $(\K^-_{\L})_{\subo}$  regular.
  From the classification of 1-connected, positively curved
homogeneous spaces (Table \ref{trans} and Table \ref{homog}), we see
that

\begin{itemize}
\item
The pair
      $(\L,(\K^-_{\L})_{\subo})$ is one of $(\Sp(d),\Sp(d-1))$ or $(\SU(d+1),\SU(d))$
      with $d\ge 3$.
\end{itemize}

 Note that since $(\K^-_{\L})_{\subo}$ is
simple, $\Km_{\S}$ is finite and $\Km_{\Delta}$ of rank one, it
follows that $\K^-_{\L}$ acts transitively on $\Sph^{l_-}$, unless
$\K^-_{\L} = \H_{\L}$. In the latter case we can apply the Chain
Theorem, and hence we can assume that $\K^-_{\L}$ indeed acts
transitively on $\Sph^{l_-}$.

\smallskip

Consider the case $(\L,(\K^-_{\L})_{\subo})=(\Sp(d),\Sp(d-1))$. Clearly,
the odd dimensional sphere $\Sph^{l_-} = \Sp(d-1)/(\H_{\L})_{\subo}$ is
equal to  $\Sp(d-1)/\Sp(d-2)$.
      If $d \ge 4$, we can again apply the Chain Theorem.
      In the remaining case consider the reduction $M^{\Sp(1)}_c$
corresponding to a standard $\Sp(1) \subset \H_{\L} \subset \Sp(3) =
\L$ which has a \co one action by  $\S \times \Sp(2)$.  From our
classification in the previous section it follows that it must be a
sum action or a modified sum action. But in that case both $\Kpm\cap
\S$ are  either trivial or all of $\S$. This is a contradiction
since $\K^-_{\S}$ is nontrivial and   finite.

\smallskip

In the case of $(\L,(\K^-_{\L})_{\subo}) =(\SU(d+1),\SU(d))$, $d\ge 3$, we
see as above  that $\Sph^{l_-} = \SU(d)/(\H_{\L})_{\subo}$ is one of
$\SU(d)/\SU(d-1)$, or $\SU(4)/\Sp(2)$. In particular, we can appeal
to the Chain Theorem when $d \ge 4$.

 If $ \Sph^{l_-} =\SU(4)/\Sp(2)$, we obtain a contradiction
to the Isotropy Lemma since the $8$-dimensional
      representation  of  $\SU(5)/\Sp(2)$
      on the orthogonal complement of $\U(4)$
      can only degenerate in $\Sp(3)/\Sp(2 )$, but $\Sp(3) \nsubseteq
\SU(5)$.

\smallskip

It remains to consider the case $\Sph^{l_-} =\SU(3)/\SU(2)$. Since
$\Kmo \supset \SU(3)$, the group $(\K^-_\Delta)_{\subo}$ must be $\S^1$
and hence $\Kmo=\Delta \S^1\cdot \SU(3)$ and $\Ho=\S^1\cdot \SU(2)$,
although the precise embedding of $\S^1\subset\Ho$ is still to be
determined. In any case, the projection of $\Ho$ onto the first
factor $\S$ is also given by a circle and hence $\Ho$ has a two
dimensional representation (inside $\S$) which necessarily
degenerates in $\Kp$. Hence $\Sph^{l_+}$ is either
$\Sph^2=\S^3/\S^1$ or $\Sph^3=\S^3\cdot\S^1/\S^1$ and all groups are
connected. In both cases primitivity implies that  $\Kp$ projects
onto $\S$ and hence in both cases $\Delta\SU(2)\cdot\SU(2)$ must be
contained in $\Kp$.

\smallskip

If $\Kp/\H = \Sph^2$, we have $\Kp=\Delta\SU(2)\cdot\SU(2)$ which
determines the embedding of $\H$
       and hence the whole group diagram is determined. The
action is the tensor product action of $\SU(2)\times\SU(4)$ on
$\Sph^{15}$, but this contradicts the fact that the action of
$\S=\SU(2)$  was assumed to be non free on the left singular orbit.

\smallskip

If $\Kp/\H = \Sph^3$, we have $\Kp=\Delta\SU(2)\cdot\SU(2)\cdot\S^1$
and hence $w_+$ can be represented by a central element in $\G$. But
then $\G/\Kp$  is totally geodesic, which is not possible.
\end{proof}


\section{Non Simple Groups without Rank $1$ Normal
Subgroups.}\label{sec: not simple}

It remains to consider semisimple groups $\G$  without normal
subgroups of rank one. In this section we deal with the non simple
case, and prove the following

\begin{prop}
Let $\G$ be a non simple semisimple group without normal subgroups
of rank one. If $\G$ acts essentially with corank 2, it is the
tensor product action of $\Sp(2)\Sp(k)$ on $\Sph^{8n-1}$.
\end{prop}

\begin{proof}

Allowing a finite kernel $\F \subset \H$ for the action, we can assume
that
$\G=\L_1\times \L_2$ with
$\rank(\L_i)\ge 2$, and none of the
$\L_i$ have normal subgroups of rank one. We let
      $\pr_i\colon \G \rightarrow \L_i$ denote the projections, and set
$\Kpm_i = \Kpm \cap \L_i$, and $\H_i = \H \cap \L_i$. There  are
connected normal subgroups $\K_\Delta $ of $\Ko$ and $\H_\Delta$ of
$\Ho$  embedded diagonally in $\L_1\times \L_2$ such that $\Kpmo =
(\Kpm_1\cdot\K_\Delta\cdot\Kpm_2)_{\subo}$ and
$\Ho=(\H_1\cdot\H_\Delta\cdot\H_2)_{\subo}$.

\smallskip

We first claim that at least one of  the four groups $\Kpm_i$ acts
transitively on $\Sph^{l_{\pm}}$:

\smallskip

\no If one of $\H_i$, say $\H_1$ is non trivial when the action is
made effective, then one of $\Kpm_1$ acts transitively, since
otherwise they both act freely or trivially which implies that
$\H_1$ would be a subset of $\H_-\cap \H_+ = \F$, contradicting
primitivity \eqref{kernel}.

\no If both $\H_i \subset \F$, we see that $\Ho=\H_\Delta$ embeds
diagonally  in $\L_1\times \L_2$, and as a consequence $\rank \H =
\rank \L_i = 2$. Now assume w.l.o.g. that $\Km$ has corank one in
$\G$. From the Product Lemma it follows as before that $\K^-_\Delta$
is not trivial of rank one and hence each of $\Km_i$ has rank one.
Thus all simple subgroups of $\Km$, and hence of $\H$ as well, have
rank one. In particular $\Sph^{l_-}$ is one of $\Sph^1=\T^2/\S^1$,
 $\Sph^3=\S^3\cdot\S^1/\Delta\S^1$ , or
$\Sph^3=\S^3\cdot\S^3/\Delta\S^3$. If one of $\Km_i$ is three
dimensional, it clearly must act transitively on $\Km/\H$ and the
same is true if $\Km$ and hence $\H$ are abelian. Hence we need to
rule out the case $\K_{\Delta}^-=\S^3$ and $(\Km_1)_{\subo} \cong
(\Km_2)_{\subo} \cong \S^1$, with $\Ho=\T^2$ embedded into the maximal
torus of $\Km$, such that it is onto $\Km_1\cdot\Km_2$.  Since
$\rank(\pr_i\Km) = \rank(\L_i) = 2$, we see that the isotropy
representation of $\L_1\times \L_2/\Kmo$ consists of a 3-dimensional
representation and all other  irreducible subrepresentations  are
even dimensional and  pairwise inequivalent. It follows that there
is an induced Riemannian submersion
$$\pi: \L_1\times\L_2/\Kmo \to
\L_1/\pr_1(\Kmo) \times \L_2/\pr_2(\Kmo)$$

\no where the latter is equipped with a product metric.   Let $\iota
= (\iota_1, \iota_2)$ denote the central element in $\Km_\Delta\cong
\S^3$. Since $\iota$ acts by the antipodal map on the slice, the
fixed point component $V$ of $M^{\iota}$ containing $c(-1)$ is the
positively curved homogeneous manifold
$(\N(\iota_1)\times\N(\iota_2))/\Km\subset\L_1\times\L_2/\Km$. Since
$\Km\cong \S^3\times\T^2$, the classification of positively curved
homogeneous spaces(cf. Table \ref{trans} and \ref{homog}) implies
that $V=\S^3\times\S^3/\Delta\S^3$ effectively. Hence neither
$\iota_i$ can be central in $\L_i$ and we  let $U_i \subset T(
\L_i/\Km_i) \subset T (\L_1\times\L_2/\Km)$ be the proper subspaces
on which $\iota$ acts by $-id$. Then $U_1 \oplus U_2$ is horizontal
with respect to the submersion $\pi$. But in the base, any plane
spanned by $u_i \in U_i, i = 1,2$ has curvature zero, so in the
total space it has nonpositive curvature intrinsically. This,
however, yields the desired contradiction since by equivariance of
the second fundamental form, $U_1 \oplus U_2$ is totally geodesic.

\bigskip

All in all it is no loss of generality to assume that say

\begin{itemize}
\item
\text{$\Km_1$ acts transitively $\Sph^{l_-}$}
\end{itemize}

Since in this case $\Km = \Km_1 \cdot \H$, the Weyl group element $w_-$ may
be represented by an element in $\L_1$. Thus $\pr_2(w_- \Kp
w_-)=\pr_2(\Kp)$ and since
      $\pr_2(\Km)=\pr_2(\H)\subset \pr_2(\Kp)$,
      we can employ Linear Primitivity to see that
      $\pr_2(\Kp)= \L_2$. In particular $\Kp_2 \triangleleft \G$, and hence
$\Kp_2 = \{1\}$ since the action is essential. It follows that
 $\Kp_\Delta\cong \L_2$ has rank two and thus:

\begin{itemize}
\item
\text{$\Kp$ has corank two in $\G$, and $\rank \L_2 = 2$}
\end{itemize}

Since $\Kp$ and $\H$ have the same rank, either $\Kp_1 = \H_1$, or
$\Kp_\Delta=\H_\Delta$. The latter would imply that the subaction by
$\L_1$ is \co one. Hence we can assume that $\Kp_1 = \H_1$, and
$\Kp_\Delta$ acts transitively on $\Sph^{l_+}$. Since $l_+$ is even,
$\L_2$ is either $\Sp(2)$ or $\G_2$, corresponding to $\Sph^{l_+}$
either $\Sp(2)/\Sp(1)\Sp(1)$ or $\G_2/\SU(3)$. The latter, however,
is impossible since then $\H$ would contain $\SU(3)$ embedded
diagonally in $\L_1\times \G_2$ in contradiction to the Isotropy
Lemma. In summary, using in addition the fact that $\Km$ must be of
corank one and $\H_2=\{1\}$,  we have:

\begin{itemize}
\item
\text{$\L_1 \times \L_2= \L_1 \times \Sp(2)$}
\item
\text{$\Kp= \H_1\Delta\Sp(2)$} and $\H=\H_1\Delta\Sp(1)^2$

\item
\text{$\Km_1=\Km_1\Km_\Delta\Km_2$ with $\Km_\Delta$ of rank one and
$\Km_2$ acting freely.}
\end{itemize}

Since   $\Sp(1)^2$ in $\H$ is  embedded diagonally, one $\Sp(1)$
must agree with $\Km_\Delta$ and the other must be embedded
diagonally in $\Km_1\Km_2$. From the classification of transitive
actions on spheres, it follows that $\Km_2=\Sp(1)$ and
$\Km_1=\Sp(k)$ with $k\ge 1$ and hence $\H_1=\Sp(k-1)$. It remains
to determine $\L_1$. From our group diagram we have so far, it
follows that $\pr_1(\Km)=\Sp(k)\Sp(1)$ and
$\pr_1(\Kp)=\Sp(k-1)\Sp(2)$ are equal rank subgroups of $\L_1$.
This implies
that $\L_1=\Sp(k+1)$. The group diagram is now determined and the
action is the tensor product action of $\Sp(k+1)\Sp(2)$ on
$\Sph^{4k+11}$.
\end{proof}


\section{Simple Groups.}\label{sec: simple}

In this section we will show that a simple group of rank at least
four either does not act isometrically on an odd dimensional
positively curved 1-connected manifold, or that it acts linearly on
a sphere.

\begin{prop}
There are no  actions of corank two for $\G=\Sp(k), k \ge 4$.
\end{prop}

\begin{proof} Recall that we already saw that $\G=\Sp(2)$ and
$\G=\Sp(3)$ do not act with corank two on a positively curved \com.

  If $\H$ contains $\Sp(1)$
embedded as a standard $1\times 1$ block, then the reduction
$M^{\Sp(1)}_c$ is odd dimensional, and $\Sp(k-1)$ acts by
cohomogeneity one on it. By induction, such an action does not
exist. Thus we may assume that
     $\H$ does not contain a $1\times 1$ block.

Since  $\rank(\H)\ge 2$, we can find an involution $\iota_1\in \T
\subset \Ho$
     that is not central in $\Sp(k)$. The reduction $M^{\iota_1}_c$ is
odd dimensional and supports a cohomogeneity one action of
     $\Sp(k-l)\cdot \Sp(l)$. From our induction hypothesis, this action
is a tensor product or a sum action and hence $\H$ contains a
$1\times 1$ block unless $(k,l)=(4,2)$. It remains to consider  the
tensor product action of $\Sp(2)\times\Sp(2)$, whose principal
isotropy group and hence also $\H$ contains $\Delta (\Sp(1)\times
\Sp(1))\subset \Delta \Sp(2)\subset \Sp(2)\times \Sp(2) \subset
\Sp(4)$.
     Now pick
     $\iota_2=\diag(-1,1,-1,1) \in \H \subset \Sp(4)$, and note that the
reduction
$M^{\iota_2}_c$ supports
     a cohomogeneity one action of $\Sp(2)\times \Sp(2)$ corresponding to
the $(1,3)$ and $(2,4)$ blocks, but with principal isotropy
containing the above  $\Delta (\Sp(1)\times \Sp(1))$ since $\iota_2$
is central in it. In particular, the principal isotropy group of
this action has a three dimensional intersection with either of the
two $\Sp(2)$ factors. But such a linear action does not exist.
\end{proof}

The case of $\G=\SU(k)$ with $k \ge 5$ is harder since  there is an
exceptional cohomogeneity one action of $\SU(5)$  on
    $\Sph^{19}$, and the fact that $\SU(4)$ acts essentially on both
$\Sph^{13}$ and on the Bazaikin spaces $B_p$, which can hence occur
in a reduction.

\begin{prop}
The linear action of $\SU(5)$ on $\Sph^{19}$  is the only essential
cohomogeneity one action by $\SU(k), k \ge 5$ of corank two.
\end{prop}

\begin{proof}

We first claim that $\H$ contains $\SU(2)$ embedded as
    a standard $2\times 2$--block. To see this, choose an element $\iota\in\Ho$
    of order $2$ that is not central in $\SU(k)$.
     Then $ \S(\U(k-2l)\U(2l))$
acts by
    cohomogeneity one on the reduction
$M^{\iota}_c$.  For $\max\{k-2l,2l\}\ge 4$
    we see that either the kernel of the action and in particular $\gH$ contains
    a $2\times 2$ block, or else the action must be a tensor product action,
    a sum action, or
the action of $\U(5)$ on $\Sph^{19}$ or $\U(4)$ on $\Sph^{13}$. In
either case we again obtain a $2\times 2$ block in $\gH$ . Thus we
may assume $(k,l)=(5,1)$ and the universal cover of $M^{\iota}_c$
    is $\Sph^{11}$ endowed with the tensor product action of
$\SU(3)\U(2)$ with principal isotropy group $\T^2$. Since in this case $\iota =
\diag(1,1,1,-1,-1)$, it follows that $M^{\iota}_c$ admits an action
of
    $\SU(3)\cdot \SO(3)\cdot \S^1$, and is therefore $\RP^{11}$.
    From the connectedness lemma we deduce that the codimension of $M^{\iota}_c$
    is strictly larger than $10$. Thus $\dim(M)=23$ and $\Ho\cong \T^2$.
    The singular orbits in $M^{\iota}_c$ have codimensions $3$ and $4$.
    Since $\Ho\cong \T^2$, these codimensions necessarily coincide with the
codimensions in $M$, all groups are connected, and we see that
$\Km,\Kp\subset N(\iota)$ -- a contradiction to primitivity.
    \\[1ex]
   From the fact that  $\H \subset \SU(k)$ contains $\SU(2)$ embedded
    as a standard $2\times 2$--block we proceed as follows:
    The reduction $M^{\SU(2)}_c$ supports a
cohomogeneity one action
    by $\SU(k-2)\cdot \S^1$. By induction, this corank two action
satisfies one of the following

      \begin{enumerate}
       \item[$\bullet$] The action is a sum action and
         $\SU(k-3)\subset \SU(k-2)$ is contained in the principal
isotropy group, or $k=6$ and the action is a sum action of
$\Spin(6)\cdot\S^1$ which contains $\Sp(2)$ in its principal
isotropy group.
       \item[$\bullet$] The action is orbit equivalent to the subaction
            of the $\SU(k-2)$-factor. This can only
             occur for $k= 6$ for the exceptional actions
            on $\Sph^{13}$ or $B_p$, and for  $k=7$ for the exceptional action
            on $\Sph^{19}$. In all cases, the
isotropy group
            contains an $\SU(2)$ embedded as  a $2\times 2$ block, and
            in the last case $\SU(2)^2$ embedded as two
            $2\times 2 $--blocks.
       \item[$\bullet$] $k = 6$ and the action is given as the tensor
product action of
            $\S^1\cdot \Spin(6)$ on $\Sph^{11}$ and the principal
            isotropy group contains $\SU(2)^2\subset \SU(4)$ embedded as two
            $2\times 2 $--blocks.
      \end{enumerate}

\no Clearly then for $k \ge 8$, we see that $\H$ contains
    $\SU(k-3)$ embedded as a standard $(k-3)\times (k-3)$
    block and we are done by the Chain Theorem. It remains to deal
with the cases $k=5,6,7$.

\smallskip

\begin{center}
  $\G=\SU(5)$
\end{center}

  By the above reduction argument we see that
$\H$ contains another $\SU(2)$ block.
    If $\dim(\H)>6$, then $\Ho$ is an equal rank extension of $\SU(2)^2
\subset \SU(5)$
    and hence $\Ho=\Sp(2) \subset \SU(4) \subset
\SU(5)$. But the irreducible $8$-dimensional
    representation of $\SU(4) \subset \SU(5)$ restricted to
$\Ho = \Sp(2)$ can not  degenerate since $\Sp(3)$ is not contained in $\SU(5)$.
Thus $\Ho=\SU(2)^2$.

Note that the $8$-dimensional
    representation of $\SU(4) \subset \SU(5)$ restricted to
$\Ho = \SU(2)\times \SU(2)$ splits as a sum of two  four dimensional
    representations each of which is acted on non trivially by exactly one of the
$\SU(2)$ factors.
    We may assume that such a  representation degenerates in $\Km$, and hence
$\Kmo=\SU(3)\cdot \SU(2)\subset \SU(5)$.
    There is also a $4$-dimensional irreducible subrepresentation of $\Ho =
\SU(2)\times \SU(2) \subset \Sp(2)$
    and  the Isotropy Lemma  implies that $\Kpo= \Sp(2)$. All groups
    are connected
    and we have recovered the picture of $\Sph^{19}$.

\smallskip

\begin{center}
  $\G=\SU(6)$
\end{center}

First suppose that the rank three group $\H$ contains
    $\Sp(2)\subset\SU(4)$. We can assume that
$\Sp(2)$ is a normal subgroup of $\H$, since otherwise $\H$ is
$\SU(4)$ and the chain theorem applies, or $\H$ is $\Sp(3)$, which
is maximal and thus $\G$ has a fixed point.
    Since the isotropy representation of
    $\SU(6)/\Sp(2)$ has an irreducible $8$-dimensional subrepresentation
coming from $\Sp(2) \subset \SU(4) \subset \SU(5)$,
    we can employ the Isotropy Lemma to see that one of the
    isotropy groups, say $\Km$, contains $\Symp(3)$ as a normal subgroup.
   But this is impossible since  we also have $\rank(\Km)=4$ and
     $\Sp(3) \subset \SU(6)$ is a maximal connected subgroup.

 Now we can assume   that $\Ho$ contains another
$\SU(2)$ block.
    Let $\iota $ be the product of the central elements of the $2$ blocks, i.e.,
    up to conjugacy $\iota=\diag(1,1,-1,-1,-1,-1) \in \S(\U(2)\U(4))$
lies in $\Ho$. The reduction $M^{\iota}_c$ is an odd dimensional
manifolds which supports a cohomogeneity one action by $
\S(\U(2)\U(4))/\iota = \SU(2)\cdot \S^1\cdot \SO(6)  $
    whose principal
    isotropy group contains the lower $4\times 4$-block $\SO(4) =
\SU(2)\SU(2)/\iota$ of $\SO(6)$. If the action is a sum action $\H$
contains $\Sp(2)$, which we already dealt with.

If the action is  the tensor product action, it is $\SU(2)$
ineffective and $\H$ contains the third $2\times 2$-block. Then
$\Ho=\SU(2)^3$, since otherwise $\Ho=\Sp(1)\Sp(2)$, which we already
dealt with. At one singular orbit say
     $\Km/\H$ the trivial representation of $\gH_{\subo}$ has to
     degenerate,
    which can only happen in a codimension $2$
    orbit. Thus $\gH_{\subo}$ is normal in  $\Km$. Also, at least one the
three $\SU(2)$
factors of $\H$
    is also normal in $\Kp$, contradicting primitivity.

\smallskip

\begin{center}
  $\G=\SU(7)$
\end{center}

  From the reduction argument above,
it follows that $\H$ contains $\SU(2)^3$ embedded as three
    $2\times 2$-blocks. Hence the element
    $\iota=\diag(1,-1,-1,-1,-1,-1,-1)$ lies in $\H$ up to conjugacy.
The reduction $M^{\iota}_c$ admits a cohomogeneity one action of
$\SU(6)\times \S^1$
    which must be a sum action. Hence $\H$ contains $\SU(5)$ and
     the chain theorem applies.
\end{proof}

For $\G = \Spin(k), k \ge 8$ we have:

\begin{prop}
There are no essential  cohomogeneity one actions of corank two by
$\Spin(k), k \ge 8$, other than the exceptional linear actions of
$\Spin(8)$ on $\Sph^{15}$  and $\Spin(10)$ on $\Sph^{31}$.
\end{prop}

\begin{proof}
We will separately treat the cases $k=8,9,10$, and  $k\ge 11$.

\smallskip

\begin{center}
  $\G=\Spin(8)$
\end{center}

  In the case of $\Spin(8)$  we can assume,  by the Chain Theorem, that
     $\gH$ even up to an outer automorphism of
    $\Spin(8)$ does not contain a $3\times 3$ block. This is
    particularly useful since there
     exists an outer automorphism which takes the standard
    $\SU(4)\subset \Spin(8)$ into the standard
    $\Spin(6)\subset\Spin(8)$ and $\Sp(2)$ into $\Spin(5)$.

    Since $\rank(\Ho)=2$, $\Ho$ is one of $\G_2 \, ,\,  \Sp(2) \, ,\,
\SU(3)$,
    $ \S^3\cdot \S^3\, ,\,  \S^1\cdot \SU(2)\,$ or \, $\T^2$. We
deal with each case  separately, and
   we apply Table \ref{sph} to determine the embeddings.

If $\Ho=\G_2$, the groups $\Kpm$ must be $\Spin(7)$. There are 3
such $\Spin(7)$ in $\Spin(8)$ which are taken into each other by the
outer automorphisms of $\Spin(8)$. Primitivity then determines the
group diagram and $M$ is $\Sph^{15}$.

If $\Ho=\Sp(2)\subset \SU(4)\subset\Spin(8)$, an outer automorphism
takes $\Sp(2)$ into a $5\times5$ block, and the Chain Theorem
applies.

If $\Ho=\SU(3)\subset \SU(4) = \Spin(6)\subset\Spin(8)$ the subgroup
$\L=\SU(2)$ in $\Ho$ is normal in $\SU(2)\SU(2)\subset\SU(4)$ which also,
via an outer automorphism, is a $4\times4$ block in $\Spin(8)$. The
normalizer of this $\SU(2)$ is therefore $(\S^3)^4$, and hence $(\S^3)^3$ acts
by \co one on the reduction $M^{\L}_c$   with a one dimensional principal
isotropy group. As we know such an action does not exist.

If $\Ho=\S^3\cdot\S^3$ we see from Table \ref{sph} that the $\S^3$
factors either sit as a $3\times 3$ block, as a Hopf action on
$\R^8$, or as a normal subgroup of a $4\times 4$ block.  In the
second case, up to an outer automorphism, the embedding is also
given by a  $3\times 3$ block. By the Chain Theorem it suffices to
consider the case  that both $\S^3$ factors are given as normal
subgroups of a $4\times 4$ block. But then up to an automorphism
$\Ho$ is a $4\times 4$ block.

If $\Ho=\S^3\cdot\S^1$, we can assume as before that $\S^3$ is given
by a normal subgroup of a $4\times 4$ block. Then $M^{\S^3}_c$
admits a cohomogeneity one action of $\Spin(4)\times \S^3$ with one
dimensional principal isotropy group. But such an action does not
exist.

If $\Ho = \T^2$ is abelian, choose an element
    $\iota\in\gH_{\subo}$ for which $\iota^2$ but not $\iota$ is in the center of
    $\Spin(8)$.
    Then the reduction $M^{\iota}_c$ admits a cohomogeneity one action
    of $\Spin(4)\cdot \Spin(4)$ or $\Spin(6)\cdot \Spin(2)=\SU(4)\cdot \S^1$
    with a $2$-dimensional principal isotropy group.
    By our induction assumption such an action does not exist.

\smallskip

\begin{center}
  $\G=\Spin(9)$
\end{center}

We can think of the maximal torus $\T^2$ in  $\Ho$ as a subtorus in
$\S^1\cdot \SU(4)\subset \Spin(8)$. Choose an involution $\iota \in
\T^2\cap \SU(4)$.
 The normalizer
    $\N(\iota)_{\subo}$  is then either $\Spin(8)$ or
    $\Spin(5)\cdot \Spin(4)$, and
    the reduction $M^{\iota}_c$ supports a cohomogeneity one action by
    $\N(\iota)_{\subo}/\ml\iota\mr$  with  principal isotropy group of corank 2.

It is easy to rule out the possibility $\N(\iota)_{\subo}=\Spin(8)$.
Indeed,  the reduction $M^{\iota}_c$ clearly has codimension $\le 8$
and $\dim M \ge 22$ since $\dim\H \le 14$. Thus $M^{\iota}_c$ is
simply connected by the Connectedness Lemma. Hence the action of
$\Spin(8)$ would have to be the exceptional action on $\Sph^{15}$,
which contradicts the fact that the action is  by
$\Spin(8)/\ml\iota\mr \cong \SO(8)$.

Thus we may assume that $\N(\iota)_{\subo}=\Spin(4)\cdot \Spin(5)$. If the
action on $M^{\iota}_c$ were almost effective or $\Spin(4)$ or
$\Spin(5)$ its ineffective kernel, $\H$ would contain a $3\times
3$-block. Hence we can assume that a normal subgroup of $\Spin(4)$
is contained in $\H$ and that the action is a sum action of
$\Spin(3)\cdot \Spin(5)$. If the second factor acts as $\SO(5)$,
$\H$ again contains a $3\times 3$-block. If on the other hand the
second factor acts as $\Sp(2)$, $\H$ contains $\Sp(1)\subset
\Sp(1)\times \Sp(1) \subset
  \Sp(2)$ which is a normal subgroup in $\Spin(4)\subset \Spin(5)$.
  In this case, the involution $(-1,-1)\in \Sp(1)\times \Sp(1)\in\H$
has $\Spin(8)$
  as its normalizer. As seen above, this is impossible.

\smallskip

\begin{center}
  $\G=\Spin(10)$
\end{center}

  We choose an involution $\iota\in \gH$ that
is not central in $\Spin(10)$.
    Then $\N(\iota)_{\subo}$ is given by $\Spin(2)\cdot \Spin(8)$ or by
    $\Spin(4)\cdot \Spin(6)$, and it acts on the reduction $M^{\iota}_c$ with
cohomogeneity one and with principal isotropy group of corank $2$.\\[1ex]
If $\N(\iota)_{\subo}=\Spin(4)\cdot \Spin(6)$, then we argue as in the case
of $\Spin(4)\cdot \Spin(5)\subset \Spin(9)$ that $\H$ contains an
$\SU(2)$ normal in $\Spin(4)$ and an
$\SU(2)\subset\SU(2)\SU(2)\subset\SU(4)$ from the sum or tensor
product action of $\SU(2)\SU(4)$. This $\SU(2)$ is a normal subgroup
of $\Spin(4)\subset\Spin(6)$ and we can find a different $\iota$
with $N(\iota)_{\subo}= \Spin(2)\cdot \Spin(8)$.

Assume now that $N(\iota)_{\subo}=\Spin(8)\cdot \Spin(2)$.
    If $\H$ contains the $\Spin(2)$-factor, then by induction it
    must also contain
$\G_2\subset \Spin(8)$. It follows that the isotropy
    representation of $\G/\gH_{\subo}$ contains a nontrivial tensor product
    of $\Spin(2)$ and $\G_2$ coming from the tensor product representation of
   $\Spin(8)\cdot \Spin(2)$ in $\Spin(10)$.  But then $\G/\Ho$
    is not spherical.

The only other possibilities for the action of $\Spin(8)\cdot
\Spin(2)$ on the reduction $M^{\iota}_c$
    is that up to an outer automorphism and possibly a covering
     it is a
    tensor product or sum action. By the Chain Theorem we
     can also assume that $\gH$
contains no  $6\times 6$-block. Hence, if it is a tensor product
action,  we can assume that $\gH$  contains $\SU(4)$, and since
$\SU(4)$ is not of equal rank in any group, it follows that
$\Ho=\SU(4)$. Similarly, if the reduction comes from a sum action,
$\Ho=\Spin(7)\subset\Spin(8)$ via the spin representation.

If $\Ho=\SU(4)$, then $\Ho$ has a six dimensional representation
from $\Ho=\Spin(6)\subset\Spin(8)$ and an eight dimensional
representation orthogonal to $\Spin(8)$. They necessarily have to
    degenerate in different orbits  and hence $\H$ is connected,
$\Km=\Spin(7)$,
    $\Kp=\SU(5)$ and we have recovered the action of $\Spin(10)$ on
    $\Sph^{31}$.

If $\Ho=\Spin(7)$, then $\Ho$ has a $7$-dimensional,
    two $8$-dimensional and a trivial representation. The
$8-$dimensional
    representation can only degenerate in $\Kmo=\Spin(9)$ and the
    trivial representation in
    $\Kp=\Spin(2)\cdot \Spin(7)$. The order two element in the center
    of $\Spin(10)$ is contained in $\Spin(9)$ and hence not in $\H$.
    Since $\Ho=\Spin(7)$ has a one dimensional centralizer in $\Spin(10)$,
      $\Kp = \Spin(2)\Spin(7) \subset
      \Spin(2)\Spin(8)\subset\Spin(10)$. It follows that
     the central element of $\Spin(10)$ must also be contained in
    $\Spin(2)\subset\Kp$ and hence $\G/\Kp$ is totally geodesic -- a contradiction.

\smallskip

\begin{center}
  $\G=\Spin(k)$ with $k \ge 11$
\end{center}

We let $\gC$ denote the center of $\Spin(k)$. We first consider
 the special case that the subaction of $\gC$ on $M$ has more than
 one orbit type. Then we may assume $\gK_-\cap \gC\neq \gH\cap \gC$.
 Clearly $\gK_-\cap \gC$ acts freely on the normal sphere and hence
 $\G/\gK_-$ is totally geodesic. This implies $\gK_-$ contains $\Spin(k-1)$
 and $\gH$ contains $\Spin(k-2)$ -- a contradiction.

Thus $\gC$ acts with one orbit type and $M/\gC$ is a manifold.
 We now drop the assumption that $M$ is simply connected and replace
 $M$ by $M/\gC$. We also replace $\Spin(k)$ by $\SO(k)$ and
 $\gC$ by the center of $\SO(k)$.

Choose an involution $\iota\in \gH\subset \SO(k)$ which is not contained
 in $\gC$. Then $N(\iota)=\SO(2h)\cdot \SO(k-2h)$.
 Given that $\rank(\gH)\ge 3$ for $k\ge 11$ and $\rank(\gH)\ge 4$
 for $k\ge 12$ we can arrange for $h\ge 2$ and $k-2h\ge 3$.

Notice that $\Fix(\iota)$ has a component $M'$ with a cohomogeneity one
 action of $\SO(2h)\cdot \SO(k-2h)$. The kernel of the action
 contains $\iota$ as well as $\gC$. Thus the center of
 $\SO(2h)\cdot \SO(k-2h)$ is contained in kernel of the action.
 We can assume that up to a covering this action is induced by
 a representation of $\Spin(2h)\times \Spin(k-2h)$ on a sphere
 (with principal isotropy group of corank 2).
 Furthermore the center of $\Spin(2h)\times \Spin(k-2h)$ acts on the sphere
 with one orbit type. It is easy to see that such a representation does not
 exist.

\end{proof}

\begin{prop}
There are no cohomogeneity one actions with corank two of any of $\G
= \F_4, \E_6, \E_7,$ or $\E_8$.
\end{prop}

\begin{proof}
If  $\G=\F_4$, choose an involution $\iota_1\in \Ho$. Then
$\N(\iota)_{\subo}=\Spin(9)$ or
    $\Sp(1)\cdot \Sp(3)$ acts by \co one on the reduction $M^{\iota}_c$
    with corank two.
     As we have seen, this rules out
$\N(\iota)_{\subo}=\Spin(9)$.
    If $\N(\iota)_{\subo}=\Sp(1)\cdot \Sp(3)$ then $\H$ contains $\Sp(2) \subset
\Sp(1)\cdot \Sp(3) \subset \F_4$,
    and there is a different involution $\iota_2 = \diag(-1,-1) \in \Sp(2)$.
    Its normalizer cannot be another $\Sp(1)\cdot \Sp(3)$
     since $\iota_2$ is central in $\Sp(2)$  and hence  cannot be
    central in the new $\Sp(3)$. Therefore we again have
     $\N(\iota)_{\subo}=\Spin(9)$ and  we
    obtain a contradiction.

If $\G=\E_6$, choose an involution $\iota\in \Ho$. Then
    $\N(\iota)_{\subo}= \SU(6)\cdot \SU(2)$ or $\Spin(10)\cdot \S^1$ and
    by induction we see that $\gH$ for any of the possible  actions of
these groups on the reduction $M^{\iota_1}_c$ must contain $\SU(4)$.

    Choose next $\iota_2=\diag(1,1,-1,-1)\in \SU(4)$.
Since $\N(\iota_2)_{\subo} \cap \SU(4) = \S(\U(2)\U(2))$ it follows that
$\gH$ contains another  $\SU(4)$ whose intersection with the first
$\SU(4)$ is at most seven dimensional.
    Thus $\dim(\H)\ge 23$.
 Using Table \ref{sph},
     it follows that
    $\Ho=\Spin(8) \subset \Spin(9)\subset \F_4 \subset \E_6$, where we
have used the fact that $\Ho= \Spin(9)$ is not allowed since the 16
dimensional spin representation cannot degenerate.
      The
centralizer of $\Ho$ in $\E_6$ is at least two dimensional since the
dimension of $\E_6/\Spin(8)$ equals $50\equiv 2$ mod $8$
and $\E_6/\Spin(8)$ has a
spherical isotropy representation.

    At one of the singular orbits the trivial representation
    has to degenerate. This can only occur in a codimension
    $2$ orbit.
    At the other singular orbit one of the
    $8$-dimensional representation
    has to degenerate.  But $l_-=1$ and $l_+=8$ is a contradiction to
    the Lower and Upper Weyl Group Bound.

If $\G=\E_7$ or $\G=\E_8$, choose a noncentral involution $\iota\in
\Ho$.
   Then
$\N(\iota)_{\subo}=\SU(8) \, ,\,\Spin(12)/\Z_2\cdot \S^3$ or $\E_6\cdot
\S^1$ in the case of $\E_7$ and $\N(\iota)_{\subo}=\Spin(16)/\Z_2$ or
$\E_7\cdot \S^3$ in the case of $\E_8$.
  But by induction we know that none of these groups can act isometrically by
cohomogeneity one on a positively curved manifold with corank two.
\end{proof}


\section{$3$-Sasakian Structure of the Exceptional Families.}\label{sec: Hitchin}

In this section we
establish
the relationship (Theorem \ref{Hitchin}) between the manifolds $P_k$ and $Q_k$ and the
interesting orbifold examples due to Hitchin \cite{Hi1}:

\begin{thm}[Hitchin] There exists a unique self dual Einstein
orbifold metric $O_k$ on $\Sph^4$ with the following properties:
\begin{itemize}
\item[a)] It is invariant under the \co one action by
 $\G=\SO(3)$ with singular
orbits of codimension two.
\item[b)] It is smooth on $M \setminus B_+$.
\item[c)] Along the right hand side singular orbit $B_+=\RP^2$ it is
smooth in the orbit direction and  has angle equal to $2\pi /k$
perpendicular to it.
\end{itemize}
\end{thm}

For the
      \co one action of $\SO(3)$ on $\Sph^4$ the
      isotropy groups are given by $\Kpm=\O(2)$ embedded in two
      different blocks and $\H=\Z_2\oplus\Z_2$. There exists a similar
      action by $\SO(3)$ on $\CP^2$ given by multiplication with real
      matrices on homogenous coordinates in $\CP^2$. One easily shows that in
      this case $\Km=\SO(2) \, , \, \Kp=\O(2)$, again in two different
      blocks, and $\H=\Z_2$ generating the second component in $\Kp$.
      Conjugation in $\CP^2$ then gives rise to an $\SO(3)$ equivariant
       two fold branched
      cover $\CP^2\to \Sph^4$ with branching locus the real points
$\G/\Kp=\RP^2$
       and a two fold cover along the left hand side singular orbits.
       When $k = 2\ell$ is even, one can thus pull back  the metric
       $O_{2\ell}$
       to become an orbifold metric on $\CP^2$ with normal angle
       $2\pi/\ell$.

We now describe the relationship with 3-sasakian geometry, see
\cite{BG} for a general reference. Among the equivalent definitions, we will use the following: A metric is called 3-sasakian if
$\SU(2)$ acts isometrically and almost freely with totally geodesic
orbits of curvature 1. Moreover, for $U$ tangent to the $\SU(2)$
orbits and $X$ perpendicular to them,  $X\wedge U$ is required to be
 an eigenvector
of the curvature operator $\hat{R}$  with eigenvalue 1, in
particular  the sectional curvature $\sec(X,U)$ is equal to 1. The
dimension of the base is a multiple of 4, and its induced metric is
quaternionic K\"{a}hler with positive scalar curvature, although it
is in general only an orbifold metric. Conversely, given a
quaternionic K\"{a}hler orbifold metric on $M$ with positive scalar
curvature, one constructs the so-called Konishi bundle whose total
space has a 3-sasakian  orbifold metric, such that the quotient
gives back the original metric on $M$.  In this fashion one obtains
a one-to-one correspondence between 3-sasakian orbifold metrics and
quaternionic K\"{a}hler orbifold metrics with positive scalar
curvature. If the base has dimension 4, quaternionic K\"{a}hler is
equivalent to being self-dual Einstein and the Konishi bundle is the
$\SO(3)$ principle orbifold bundle of self dual 2-forms on the base
with the metric  given by the naturally defined connection metric.
Hence the Hitchin metrics give rise to 3-sasakian orbifold metrics
on a seven dimensional orbifold $H^7_k$. The cohomogeneity one
action by $\SO(3)$ on the base admits a lift to the total space
$H^7_k$ which commutes with the almost free principal orbifold
$\SO(3)$ action. The joint action by $\SO(3)\SO(3)$ on $H^7_k$ is
hence an isometric \co one action. In general, one would expect the
metric on $H^7_k$ to have orbifold singularities since the base
does. However, we first observe that this is not the case. Although
the claim also follows from the proof of Theorem \ref{PQ-H}, we give
a simple and  more geometric proof.

\begin{thm}\label{Hitchinsmooth}
For each $k$, the total space $H_k^7$ of the  Konishi bundle
corresponding to the selfdual Hitchin orbifold $O_k$ is a smooth
3-Sasakian manifold.
\end{thm}

\begin{proof}
Notice  that the singular orbit $B_+$ in $O_k^4 , k>2$ must be
totally geodesic. Indeed, being an orbifold singularity, one can
locally lift the metric on a normal slice $\Disc^2$ to $\RP^2$ to
its k-fold branched cover $\hat{\Disc}\to \Disc$ with an isometric
action by $\Z_k$ such that $\hat{\Disc}/Z_k= \Disc$. Hence the
singular orbit is a fixed point set of a locally defined group
action and thus totally geodesic.

      The $\SO(3)$
principle bundle $H^7_k$ is smooth over all smooth orbits in
$H^4_k$. If it has orbifold singularities, they must  consist of an
$\SO(3)\SO(3)$ orbit which projects to $B_+$, and is  again totally
geodesic by the same argument as above. This five dimensional orbit
is now 3-sasakian with respect to the natural semi-free $\SO(3)$
action on $H^7_k$, since it is totally geodesic and contains all
$\SO(3)$ orbits. But the quotient is 2-dimensional which contradicts
the fact that the base of such a manifold has dimension divisible by
4.
\end{proof}

As  mentioned in the Introduction, except for $\Sph^7 = P_1$, the
manifolds $P_k$ are the first  $2$-connected seven dimensional
$3$-Sasakian manifolds. The cohomology rings of the manifolds $Q_k$ happen to coincide with the cohomology rings of all the previously known 3-Sasakian 7-manifolds with first Betti number one. These are exactly the
Eschenburg spaces $\diag(z^a,z^b,z^c)\backslash $
$\SU(3)/\diag(1,1,\bar{z}^{a+b+c})$ with $a,b,c$ positive pairwise
relatively prime integers \cite{BGM}. They contain the 3-Sasakian
manifolds $E_k$ as a special case. All of these, as well as those with second Betti number at least two \cite{BGMR},  are constructed from the constant
curvature $3$-Sasakian metric on $\Sph^{4n+3}$, equipped with the
Hopf action, as 3-Sasakian reductions with  respect to an isometric
abelian group action commuting with the Hopf action. As a
consequence all of them are toric,  i.e.  admit an isometric action
by a $2$-torus commuting with the $\SU(2)$ action. In contrast, the
examples $P_k$ and $Q_k$, for $k \ge 2$ are not toric, since  the
orbifolds $O_k$, for  $k \ge 3$,  have $\SO(3)$ as the identity
component of their isometry group.

\smallskip

 Before verifying that the above  manifolds $H_k^7$
coincide with the ones described  in the introduction, we first
discuss a general framework for \co one orbifolds.

Observe that a group diagram as in
(\ref{diagram}), where we  assume that  $h_\pm$  are embeddings, but
$j_\pm$ are only homomorphisms with finite kernel  and $j_-\circ h_-
= j_+\circ h_+=j_{\subo}$ with $\K_\pm /
\H=\Sph^{l_{\pm}}$, defines a \co one orbifold $O$:
 The regular
orbits, being hypersurfaces, have no orbifold singularities, and we
can therefore assume that $j_{\subo}$ is an
embedding, although we still allow the action of $\G$ to be
ineffective otherwise. A neighborhood of a singular orbit is given
by $D(B_\pm)=\G\times_{\Kpm} \Disc^{\l_\pm +1}$ where $\Kpm$ acts on
$\G$ via right multiplication: $g\cdot k=gj_{\pm}(k)$ and on
$\Disc^{\l_\pm +1}$ via the natural linear extension of the action
of $\Kpm$ on $\Sph^{l_{\pm}}$. This then can be written as
$D(B_{\pm})=\G\times_{(\Kpm /\ker j_\pm)} (\Disc^{\l_\pm +1}/\ker
j_\pm)$ and the singularity normal to the smooth singular orbit
$\G/j_\pm(\Kpm)$ is  $\Sph^{l_\pm}/\ker j_\pm$. It is easy to see
that any \co one orbifold can be described in this fashion. In fact
this  follows since the frame bundle of a \co one orbifold is
a \co one manifold, and thus orbifolds inherit \co one diagrams as
described.
 In all the cases of interest here, we note that both $l_{\pm} = 1$,
 and the orbifolds  are therefore (topologically)  manifolds.

 We are now
ready to prove:

\begin{thm}\label{PQ-H}
Our manifolds  $P_k$ and $Q_k$ are
     equivariantly diffeomorphic to  the universal
covers of the 3-sasakian manifolds $H_{2k-1}$ and $H_{2k}$
respectively.
\end{thm}

\begin{proof}
Since the metrics in the Hitchin examples are smooth near $B_-$, it
follows that $\Km\cong \O(2)$ and hence $\H\cong \Z_2\oplus\Z_2$.
 Hence we can
 assume that $j_-$ is an embedding of $\Km\cong \O(2)$ into the lower
block in
     $\SO(3)$, $h_-$ the diagonal embedding $\Z_2\oplus\Z_2\subset\O(2)$,
     and via $j_-\circ h_-$ the group $\H $ is embedded as the set of diagonal
matrices in $\SO(3)$. As in the case of smooth $\SO(3)$ invariant
 metrics on $\Sph^4$, the Hitchin orbifold metrics collapse in
 different directions corresponding to $\Kpmo$, and the normal
angle along $B_+$ is $2\pi/k$ .  If we define the homomorphism  $\phi_k\colon
\SO(2)\to \SO(2)$ by $A\to A^k$, we see that
 $j_+(A) \in \SO(3)$ is $\phi_k(A)$ for $A\in \Kpo$ followed by an embedding into
 $\SO(3)$,
 which we can assume is
  in the upper block in order  to be consistent with the
   $\H$-irreducible $1$-dimensional subspaces of
$\fso(3)$.

 On the right hand side a neighborhood of the singular
     orbit is given by
     $D(B_+)=\SO(3)\times_{\Kp} \Disc_+^2$
where $\Kpo$ acts on $\SO(3)$ via $\phi_k$ and on  $\Disc_+^2$ via
$\phi_2$ since $\Kpo\cap\H=\Z_2$. The description of the disc bundle
$D(B_+)$ gives rise to a description of the corresponding (smooth)
$\SO(4)$ principle orbifold bundle $\SO(3)\times_{\Kp} \SO(4)$ where
the action of $\Kpo$ on $\SO(3)$ is given by $\phi_k$ as above, and
the action on $\SO(4)$ is given via $\SO(2)\subset\SO(4)\colon
A\in\SO(2)\to (\phi_k(A),\phi_2(A))$ acting on the splitting
$T\oplus T^\perp$ into tangent space and normal space of the
singular orbit. Similarly for the left hand side where $k=1$. In
order to take orientations into account and their consistent match
for the gluing in the middle, we start with an oriented basis
$\dot{c}(t),i,j,k$ for the regular orbits, where we have used for
simplicity the isomorphism $\fso(3)\cong \fsu(2)$. On the left hand
side  the $i$ direction collapses, $T$ is oriented by $j,k$  and
$T^\perp$ by $ \dot{c}(-1),i$. Here $i$ corresponds to the
derivative of the Jacobi field along $c$ induced by $i$. On the
right hand side the $j$ direction collapses, $T$ is oriented by
$k,i$ and $T^\perp$ by $ \dot{c}(1),j$. Here $j$ corresponds to the
negative of the derivative of the Jacobi field along $c$ induced by
$j$. Furthermore, one easily checks that $\SO(2)\subset\O(2)$ has a
positive weight on $T$ where we have endowed the isotropy groups on
the left and on the right with orientations induced by $i$ and $j$
respectively. Hence $\Kpmo\subset\SO(3)\SO(4)$ sits inside the
natural maximal torus in $\SO(3)\SO(4)$ with slopes $(1,1,2)$ on the
left, and $(k,k,-2)$ on the right.

We can now  determine the group picture for the $\SO(3)\SO(3)$
action on the principle bundle of the vector bundle of self dual two
forms. This vector bundle can also be viewed as follows: If P is the
$\SO(4)$ principle bundle of the orbifold tangent bundle of
$\Sph^4$,
 then the quotient  $P/\SU(2)$  under a
normal $\SU(2)$ in $\SO(4)$ is an $\SO(3)$ principle bundle and by
dividing by the two normal subgroups, one obtains the principle
bundles for the vector bundle of self dual and the vector bundle of
anti self dual 2 forms. This is due to the fact that the splitting
$\Lambda^2V\cong \Lambda_+^2V\oplus\Lambda_-^2V$ for an oriented
four dimensional vector space corresponds to the splitting of Lie
algebra ideals $\fso(4)\cong \fso(3)\oplus\fso(3)$ under the
isomorphism $\Lambda^2V\cong\fso(4) $. Alternatively we can first
project under the two fold cover $\SO(4)\to\SO(3)\SO(3)$ and then
divide by one of the $\SO(3)$ factors. Under the homomorphism
$\SO(4)\to\SO(3)\SO(3)$ and the natural maximal tori in $\SO(4)$ and
in $\SO(3)\SO(3)$, a slope $(p,q)$ circle goes into one with slope
$(p+q,p-q)$. Hence the slopes of $\Kpmo$ in $\SO(3)\SO(3)\SO(3)$ are
$(1,3,-1)$ on the left, and $(k,k-2,k+2)$ on the right. This also
implies that both $\SO(3)$ factors act freely on $P$. If we divide
by one of the $\SO(3)$ factors to obtain the two $\SO(3)$ (orbifold) principal
bundles, the slopes of the circles $\Kpmo$ viewed inside
$\SO(3)\SO(3)$ become $(1,3)$ on the left and $(k,k-2)$ on the right
for one principal bundle, and $(1,-1)$ on the left and $(k,k+2)$ on
the right for the other.

To see which principal bundle is the correct one for the Hitchin
metric, recall that in \cite{Hi1} one chooses an orientation on the
regular orbits in order to derive the correct differential equation
for  Eintein metrics which are self dual with respect to the given
orientation. For this fixed orientation Hitchin constructs the
solution for the self dual Einstein metrics and checks smoothness at
the singular orbits. Hence either the family of
principal bundles with slopes $\{(1,3),(k,k-2)\}$  or the one with
slopes $\{(1,-1),(k,k+2)\}$ are the desired $\SO(3)$ principle
bundle for all $k$. But we know that for $k=1$ the principle bundle
$H_1=\RP^7$ has slopes $\{(1,-1),(1,3)\}$ and for $k=2$ the bundle
$H_2=W^7/\Z_2$ has slopes $\{(1,-1),(1,2)=(2,4)\}$ (see section
\ref{sec: examples}). Hence the slopes for the principle bundles
defined by the Hitchin metric are $\{(1,-1),(k,k+2)\}$, which is, up
to covers, the $P$ family for $k$ odd and the $Q$ family for $k$
even.
\end{proof}

The second family of $\SO(3)$ principle bundles in the above proof
are the principal bundles of the vector bundle of anti-self dual two
forms, which the proof shows are also smooth. We note that in the
case of $k=3$ one obtains the slopes for the exceptional manifold
$B^7$ and in the case of $k=4$ the ones for $R$.

\smallskip

We note that in order for the \co group diagram on the frame bundle
or the principle bundle $H_k$ to be consistent, $\Kp\cong \O(2)$ for
$k$ odd, and $\Kp\cong \O(2)\times \Z_2$ for $k$ even. This also
determines the embedding of $\Kp$ into $\SO(3)\SO(3)$ and hence the
orbifold group diagram for the Hitchin metrics. The manifolds $H_k$
are two-fold subcovers of $P_k$ and $Q_k$. In the case of $P_k$ we
divide by the full center and in the case of $Q_k$ we add another
component to all three isotropy groups (see Lemma \ref{isocomp}). We
also point out that for $k=2\ell$ the total space of the Konishi
bundle associated to the lifted orbifold metric on $\CP^2$ is equal
to $Q_\ell$.

\smallskip

There is another interesting connection between self-dual Einstein
metrics and positive curvature. O.Dearricott \cite{De1} proved that
if one allows to scale  the 3-sasakian metric on a 7-manifold with
arbitrarily small scale in the direction of the $\SU(2)$ orbits,
then the metric on the total space has positive sectional curvature
if and only if the base self dual Einstein metric does. One can
apply this to the Boyer-Galicki-Mann  3-sasakian metrics \cite{BGM}
on the Eschenburg spaces $E_{a,b,c}=\diag(z^a,z^b,z^c)\backslash $
$\SU(3)/\diag(1,1,\bar{z}^{a+b+c})$ whose self dual Einstein
orbifold quotient  is a weighted projective space
$\CP^2[a+b,a+c,b+c]$. O.Dearricott showed in \cite{De2}, that many
(but not all) of the weighted projective spaces have positive
sectional curvature. The total space also admits an Eschenburg
metric with positive curvature, but the Dearricott metrics are
different in that the projection is a Riemannian submersion with
totally geodesic fibers, whereas in the Eschenburg metric the fibers
are not totally geodesic. It is hence natural to ask if the Hitchin
metrics have positive sectional curvature for some $k$ besides the
values $k=1,2$ where this is true by construction. Hitchin  gave an
explicit formula for the functions  describing his metrics for
$k=4,6$ in \cite{Hi1} and for $k=3$ in \cite{Hi2}, which  are simply
rational functions  of a parameter t along the geodesic in the first
two cases and algebraic functions in the third case. One can now
compute the sectional curvatures of the self dual Einstein metrics
in these special cases and one shows, surprisingly, that the
curvatures near the non-smooth singular orbit are all positive, but
some become negative near the smooth singular orbit. On the other
hand, it is not hard to construct 4-dimensional positively curved
orbifold metrics with these prescribed orbifold singularities.
Nevertheless, it is natural to suggest that there could be some
significance in the existence of the 3-sasakian metrics on $P_k$ and
$Q_k$ and the question whether these spaces have a metric where all
sectional curvatures are positive.

\smallskip


\section{ Topology of the Exceptional Examples.}\label{sec: topology}

      In order to prove Theorem \ref{topology}, we study the corresponding larger classes of \co
      one manifolds with arbitrary slopes.

      The class containing $P_k$
       consists of  the \co one manifolds  $M =
       M_{(p_-,q_-),(p_+,q_+)}$ where
       $
\H\subset\{\Km,\Kp\} \subset \G$ is given by $\G=\S^3 \times \S^3 ,
\{\Km , \Kp\} = \{ \gC^i_{(p_-,q_-)} \cdot \H ,  \gC^j_{(p_+,q_+)}
\cdot \H \}$ and $\H=\Q=\{\pm 1 , \pm  i, \pm j ,\pm k\}$, where
$(p_-,q_-)$ as well as $(p_+,q_+)$ are relatively prime odd
integers, and $\gC^k_{(p_,q_)}\subset \S^3\times \S^3$ is the
subgroup of elements $\{(e^{kp\gt},e^{kq\gt})\}$ as in section
\ref{sec: rank2}. It follows  that $\Kpm/\Kpmo=\Z_2$, where the
second component is generated by $(j,j)$ on the left and $(i,i)$ on
the right, up to signs (of both coordinates).
 The
embedding of $\Q$ is  determined by the slopes and is $\Delta \Q$,
up to sign changes in both
coordinates.      All
      cohomology groups, unless otherwise stated,
       are understood to be with $\Z$ coefficients.

      \begin{thm}
The manifolds $M=M_{(p_-,q_-),(p_+,q_+)}$ are 2-connected. If $\frac{p_-}{q_-} \ne \pm
\frac{p_+}{q_+} $ their cohomology ring is determined by
$\pi_3(M)=\Z_k$ with $ k=(p_-^2q_+^2-p_+^2q_-^2)/8$. Otherwise
$H^3(M)=H^4(M)=\Z$.
      \end{thm}

      \begin{proof}
We will use the same method as in \cite[Poposition 3.3]{GZ} although
the details will be significantly more difficult. In order to show
that $M$ is simply connected, one uses Van Kampen on the cover
$U_\pm = D(B_{\pm}) = \G\times_{\Kpm}\Disc^{\ell_{\pm}+1}$, which
deformation retract to $B_\pm = \G/\Kpm$, and $U_-\cap U_+ = \G/\H$.
We denote the projections of the sphere bundles by
      $\pi_\pm\colon
\G/\H= \G\times_{\Kpm}\Sph^{\ell_{\pm}}= \partial D(B_\pm) \to
B_\pm=\G/\Kpm$. For a homogeneous space $\G/\L$ with $\G$ simply
connected, the fundamental group is given by the group of components
$\L/\L_{\subo}$. This determines the homomorphisms $\pi_\pm\colon
\pi_1(\G/\H)\to \pi_1(\G/\Kpm)$ and it follows that $\pi_1(M)=0$.
For the  cohomology groups, we use the Mayer-Vietoris sequence on
the same decomposition, which gives a long exact sequence
\begin{equation} \label{MV}
\to H^{i-1}(B_-)\oplus H^{i-1}(B_+) \overset{\pi_-^* -
\pi_+^*}{\Lllongrightarrow} H^{i-1}(\G/\H) \to H^i(M)\to
H^i(B_-)\oplus H^i(B_+)\to
\end{equation}

We  first determine the cohomology groups of the singular and
regular orbits. Denote by $\mu_\pm\colon B_\pm^0=\G/\Kpmo \to B_\pm$
the natural projections, which are two fold covers. One knows that
$B^0_\pm$ are diffeomorphic to $\Sph^3\times \Sph^2 $, independent
of the slopes, see e.g. \cite[Proposition 2.3]{WZ}. For $B_\pm$ we
will show that it has the same cohomology as that of
$\Sph^3\times\RP^2$, although we do not known if they are
diffeomorphic.

\begin{lem} The cohomology of the $\G$ orbits are given by
\begin{itemize}
\item[(a)] $B_\pm$ is non-orientable with $\pi_1(B
_\pm)=\Z_2$,
 $H^0(B_\pm)=H^3(B_\pm)=\Z\,
,\, H^1(B_\pm)=H^4(B_\pm)=0$ and $H^2(B_\pm)=H^5(B_\pm)=\Z_2$.
Furthermore, $\mu_\pm^*\colon H^3(B_\pm) \to H^3(B_\pm^0)$ are
isomorphisms.

\item[(b)] The principal orbit is orientable  with $\pi_1(\G/\H)=\Q$,
$H^0(\G/\H)=H^6(\G/\H)=\Z\, ,\, H^1(\G/\H)=H^4(\G/\H)=0\, ,\,
H^2(\G/\H)=H^5(\G/\H)=\Z_2\oplus\Z_2$ and $H^3(\G/\H)=\Z\oplus\Z$.
\end{itemize}
\end{lem}

\begin{proof}
For the principal orbits, we observe that a normal subgroup
$\S^3\subset \G$ acts freely and hence  give rise to a principal
$\S^3$ bundle $\G/\H\to \S^3/Q$. This bundle must be trivial since
the classifying space $\HP^\infty$ is 4-connected. Hence $\G/\H\cong
\Sph^3\times (\Sph^3/\Q)$ and the cohomology groups of $\G/\H$
easily follow.

We first note that a singular orbit $B=\G/\K$ with $\K=\gC^i_{(p,q)}\H$ is non-orientable. This follows, since the action of $\K/\Ko$ on the tangent space of $\G/\K$ does not preserve orientation.
Considering the
projection onto the second factor in $\S^3\times \S^3$, we obtain
fibrations $L_q\to B  \overset{\sigma}{\longrightarrow} \RP^2$ and
$L_q\to B^0\overset{\sigma_{\subo}}{\longrightarrow} \Sph^2$ where the
fibers of these homogeneous fibrations  are lens spaces, since they
are of the form $((\S^3\times \S^1 )\H)/\K=\S^3/\{z^p\}$ with
$z^q=1$. It is well known that $H^*(B,\Z_{p'})=H^*(B^0,\Z_{p'})^\rho$
 for ${p'}$ a prime  different from 2, where $\rho$ is the deck
transformation of the two fold cover $\mu \colon B^0\to B$.
 The spectral sequence
for $\sigma_{\subo}$ implies that $\sigma_{\subo}^*\colon H^2(\Sph^2,\Z)\to
H^2(B^0,\Z)$ is an isomorphism. Since the deck groups of $B_{\subo}\to B$
and $\Sph^2\to \RP^2$ are compatible with the fibrations $\sigma$
and $\sigma_{\subo}$, it follows that   $\rho^*= -Id$ on $H^2(B^0,\Z)$.
Since $\rho$ reverses orientation, $\rho^* = +Id $ on
$H^3(B^0,\Z)=\Z$. Thus $H^i(B,\Z_p)=\Z_p$ for $i=0,3$ and $0$
otherwise.
 Since $q$ is odd, $H^*(L_q,\Z_2)=H^*(\Sph^3,\Z_2)$
and hence in the spectral sequence for $\sigma$ with $\Z_2$
coefficients all differentials necessarily  vanish. Thus
$H^i(B,\Z_2)=\Z_2$ for every $i$. This, together with the universal
coefficient theorem, easily determines the cohomology of $B_\pm$.

We finally show that $\mu^*\colon H^3(B)=\Z\to H^3(B^0)=\Z$ is an
isomorphism. By the universal coefficient theorem,  it suffices to
show that $\mu^*\colon H^3(B,\Z_{p'})\to H^3(B^0,\Z_{p'})$ is an
isomorphism for every prime ${p'}$. If ${p'}$ is odd, this is clearly the
case by what we proved above. For $p'=2$ we use the observation that
all differential in the spectral sequence for $\sigma$ with $\Z_2$
coefficients vanish. This implies that the edge homomorphism
$H^3(B,\Z_2)=\Z_2\to H^3(L_q,\Z_2)=\Z_2$ is onto and hence an
isomorphism. The same argument applies to the fibration $\sigma_{\subo}$
and hence $\mu^*\colon H^3(B,\Z_2)\to H^3(B^0,\Z_2)$ is an
isomorphism as well.
\end{proof}

The homomorphisms $\pi_\pm^*\colon \pi_1(\G/\H)\to \pi_1(B_\pm)$
determine $\pi_\pm^*\colon H^2(B_\pm)=\Z_2 \to
H^2(\G/\H)=\Z_2\oplus\Z_2$ via the universal coefficient theorem,
and show that $H^2(B_-)\oplus H^2(B_+)\to
H^2(\G/\H)$ is an isomorphism. Hence $H^2(M)=0$ and $M$ is
2-connected. Since we also have $H^4(B_\pm)=0$, the Mayer Vietoris
sequence implies that $H^3(M)$ is the kernel and $H^4(M)$ the
cokernel of $\pi_-^* - \pi_+^*\colon H^3(B_-)\oplus
H^3(B_+)=\Z\oplus\Z \to H^3(\G/\H)=\Z\oplus\Z $

To determine $\pi_\pm^*$, we consider the commutative diagram,
dropping the signs for the moment:
\begin{equation} \label{coh}
\begin{split}
\xymatrix{ \S^3\times \S^3 \ar[r]^{\tau} \ar[d]^{\eta} & \S^3\times
\S^3/\Ko \ar[d]^{\mu} \\
\S^3\times \S^3/\H \ar[r]^{\pi} & \S^3\times \S^3/\K }
\end{split}
\end{equation}
\no where all arrows are given by their natural projections. In
\cite[(3.6)]{GZ} it was shown that the image of a generator in
$H^3(\G/\Ko)=\Z$ is equal to $(-q^2,p^2)$, using the natural basis
in $H^3(\G)=\Z\oplus\Z$. Since $\mu^*$ is an isomorphism in degree
3, $\pi^*$ is determined as soon as we know the integral lattice
$\Im (\eta)^*\subset H^3(\S^3\times \S^3)$. Since $\S^3\times
\S^3/\H\cong S^3\times(\S^3/\Q)$ and $\S^3\to \S^3/\Q$ is an 8-fold
cover, this sublattice must have index 8. Using (\ref{coh}) for the
slopes $(1,1)$ and $(1,3)$, we see that $(-1,1)$ and $(-9,1)$ lie in
the lattice and must be a basis, since the element $(1,0)$ has order
8 in the quotient group. Using the basis $(-1,1)$ and $(4,4)$ the
matrix of $\pi_-^* - \pi_+^*$ becomes:
$$\left(\begin{array}{rr} \frac{1}{2}(p_-^2+q_-^2) & -\frac{1}{2}
(p_+^2+q_+^2) \\
\frac{1}{8}(p_-^2-q_-^2) & -\frac{1}{8}(p_+^2-q_+^2)
\end{array}\right) $$

Since $(p_-,q_-)$ are relatively prime, one easily sees that
$(\frac{1}{2}(p_\pm^2+q_\pm^2), \frac{1}{8}(p_\pm^2-q_\pm^2))$ are
relatively prime as well, which implies that the cokernel of
$\pi_-^* - \pi_+^*$ is a cyclic group. If we assume that
$\frac{p_-}{q_-} \ne \pm \frac{p_+}{q_+} $, the kernel is 0 and the
cokernel is  cyclic with order $\det (\pi_-^* - \pi_+^*)=
((p_-^2q_+^2-p_+^2q_-^2)/8$. Otherwise kernel and cokernel are equal
to $\Z$.
      \end{proof}

      \smallskip

      Next we consider the extension $N=N_{(p_-,q_-),(p_+,q_+)}$ of
      the $Q$ family, given by $\H=\{(\pm
1,\pm 1),$ $ (\pm i , \pm i )\} \subset \{\Km , \Kp\} =
\{\gC^i_{(p_-,q_-)}\H ,\gC^j_{(p_+,q_+)} \H\} \subset
   \G= \S^3\times \S^3$
 with $(p_-,q_-)$ as well as
$(p_+,q_+)$ relatively prime, $p_+$ even and $p_-,q_-,q_+$ odd.
Notice that the component groups $\Kpm/\Kpmo$ are determined by the
fact that $(i,i)\in \Kmo$ and $(1,-1)\in \Kpo$.

      \begin{thm}\label{mfdN}
The manifolds $N=N_{(p_-,q_-),(p_+,q_+)}$ are simply connected with $H^2(N)=\Z \, ,\,
H^3(N)=0$ and $H^4(N)=\Z_k$ with $k=p_-^2q_+^2-p_+^2q_-^2$.
      \end{thm}

      \begin{proof}
We indicate the changes in the proof which are necessary, and start
with the cohomology of the orbits. They contain torsion groups $S,T$
and  integers $c,d$ which are to be determined later.

\begin{lem}\label{orbits} The cohomology of the $\G$ orbits are given by
\begin{itemize}
\item[(a)] $B_-$ is orientable with $\pi_1(B_-)=\Z_2$,
      $H^0(B_-)=H^3(B_-)=H^5(B_-)=\Z\,
,\, H^1(B_-)=0 \, , \, H^2(B_-)=\Z\oplus\Z_2$ and $H^4(B_-)=\Z_2$.
Furthermore, $\mu_-^*\colon H^3(B_-) \to H^3(B_-^0)$ is
multiplication by $c$, a power of $2$.

\item[(b)] $B_+$ is non-orientable with $\pi_1(B_+)=\Z_4 \, , \,
H^0(B_+)=\Z\, ,\, H^1(B_+)=H^4(B_+)=0 \, , \, H^2(B_+)=\Z_4 \, , \,
H^3(B_+)=\Z\oplus S $ and $ H^5(B_+)=\Z_2$. Furthermore,
$\mu_+^*\colon H^3(B_+) \to H^3(B_+^0)$ is multiplication by $d$, a power of $2$, on
the free part.
\item[(c)] The principal orbit is orientable  with
$\pi_1(\G/\H)=\Z_2\oplus\Z_4$, $H^0(\G/\H)=H^6(\G/\H)=\Z\, ,\, H^1(\G/\H)=0$
, $ H^2(\G/\H)=H^5(\G/\H)=\Z_2\oplus\Z_4  \, ,\,
H^3(\G/\H)=\Z\oplus\Z\oplus T$ and $H^4(\G/\H)= T$.
\end{itemize}
     where $S$ and $T$ are  torsion groups of the form $(\Z_2)^m$.
\end{lem}

\begin{proof}
For $B=B_- = \S^3\times \S^3 /\gC^i_{(p,q)}\H$ one has
$(i,i)\in \Kmo$ and $(1,-1)$ generates the second component. Hence
$B_-$ is orientable with $\pi_1=\Z_2$.
 Projection onto the second coordinate in
$\S^3\times \S^3$ gives rise to fibrations $L_{2q}\to
B\overset{\sigma}{\longrightarrow} \Sph^2$ and $L_q\to
B^0\overset{\sigma_{\subo}}{\longrightarrow} \Sph^2$.
 Notice that the
fiber for the first fibration is $((\S^3\times\S^1)\H)/\K=\S^3/\{z^p,-1\}$
with $z^q=1$, which is $\S^3/\Z_{2q}$
since $p$ and $q$ are odd. As before, one now shows that for any
prime $p'$ different from 2, $H^i(B,\Z_{p'})=\Z_{p'}$ for $i=0,2,3,5$ and
$0$ otherwise. Since $H^*(L_{2q},\Z_2)=H^*(\RP^3,\Z_2)$ and
$H^1(B,\Z_2)=\Z_2$ it follows that all differentials vanish in the
spectral sequence for $\sigma$ with $\Z_2$ coefficients. This
determines $H^*(B,\Z_2)$ and the cohomology groups of $B$ easily
follow. Since $\mu^*\colon H^3(B,\Z_{p'})=\Z_{p'}\to H^3(B^0,\Z_{p'})=\Z_{p'}$
is an isomorphism for every prime ${p'}\ne 2$, it must be
multiplication by  a power of two over the integers.

\vspace{5pt}

For $B=B_+ = \S^3\times \S^3 /\gC^i_{(p,q)}\H$ with $p$ even
$q$ odd, the element $(i,i)$ generates the 4 components of $\K$.
 Hence $B$ is non-orientable with
$\pi_1(B)=\Z_4$. We also have  fibrations $L_{2q}\to B
\overset{\sigma}{\longrightarrow} \RP^2$ and $L_{q}\to
B^0\overset{\sigma_{\subo}}{\longrightarrow} \Sph^2$. Using the 4-fold
cover $\mu\colon B^0\to B$, it follows as before that for any prime
$p' \ne 2$  we have  $H^i(B,\Z_{p'})=\Z_{p'}$  for $i=0,3$ and $0$ otherwise.

 We now consider the spectral sequence of the fibration
$\sigma$ with $\Z_2$ coefficients. Since $H^1(B,\Z_2)=\Z_2$, it
follows that $d_2\colon E_2^{0,1}=\Z_2\to E_2^{2,0}=\Z_2$ is an
isomorphism and hence $d_2\colon E_2^{0,2}=\Z_2\to E_2^{2,1}=\Z_2$
vanishes and $d_2\colon E_2^{0,3}=\Z_2\to E_2^{2,2}=\Z_2$ is an
isomorphism as well. This determines $H^*(B,\Z_2)$ and one easily
derives the cohomology of $B$, up to a non-zero torsion group
$S=(\Z_2)^k$ in dimension three. It also follows that $\mu^*\colon
H^3(B,\Z_{p'})=\Z_{p'}\to H^3(B^0,\Z_{p'})=\Z_{p'}$ is an isomorphism for $p'\ne
2$, and hence $\mu^*\colon H^3(B,\Z)=\Z\oplus S\to H^3(B^0,\Z)=\Z$
is multiplication by a power of two on the free part.

\vspace{5pt}

     $\G/\H$ is clearly orientable with $\pi_1(\G/\H)=\Z_2\oplus\Z_4$.
     Using the 8-fold cover $\eta\colon \S^3\times\S^3\to \G/\H$
      and the fact that the deck transformations are
     homotopic to the identity, it follows  that
     the non-zero groups in  $H^i(B,\Z_p)$ for  $p\ne 2$ are $\Z_p$
for $i=0,6$ and  $\Z_p\oplus \Z_p$ for $i=3$. In the spectral
sequence for the
      fibration $\S^1\times \S^1\to \G/\H\to \Sph^2\times \Sph^2$ with $\Z_2$
      coefficients, all differential must be 0 since $H^1(\G/\H,\Z_2)=
\Z_2\oplus\Z_2$. This determines $H^*(\G/\H,\Z_2)$ and hence
$H^*(\G/\H,\Z)$, up to a non-zero torsion group $T=(\Z_2)^k$ in
dimension three and four.
\end{proof}

The homomorphisms on the group of components again show that $N$ is
simply connected and that the homomorphism $H^2(B_-)\oplus
H^2(B_+)=\Z_4\oplus\Z\oplus\Z_2\to H^2(\G/\H)=\Z_4\oplus\Z_2$ is an
isomorphism on  the torsion part, since it is determined by the
corresponding homomorphisms on the fundamental groups. Therefore
Mayer Vietoris implies that $H^2(N)=\Z$. By Poincare duality
$H^5(N)=\Z\to H^5(B_-)=\Z$ is  an isomorphism, which means that the
homomorphism $H^4(B_-)\oplus H^4(B_+)=\Z_2\to H^4(\G/\H)=T$ is onto
and hence an isomorphism with $T=\Z_2$. Since the universal
coefficient theorem  for $N$ implies that $H^3(N)$ cannot have any
torsion, $\pi_-^* - \pi_+^*\colon H^3(B_-)\oplus
H^3(B_+)=\Z\oplus\Z\oplus S \to H^3(\G/\H)=\Z\oplus\Z\oplus T $ is
injective on its torsion part, and hence an isomorphism with
$S=T=\Z_2$. Thus $H^3(N)$ is the kernel on the free part and
$H^4(N)$ its cokernel.

We next  determine the lattice generated by the image of $\eta^*$ in
$ H^3(S^3\times S^3)$. In the spectral sequence for the fibration
$\S^3\times \S^3\to \S^3\times \S^3/\Z_2\oplus\Z_4\to B_{\Z_2}\times
B_{\Z_4}$ the fundamental group $ \Z_2\oplus\Z_4 $ of the base acts
trivial in homology and hence the local coefficients become ordinary
$\Z$ coefficients. The only non-zero differential is $d_2\colon
E_2^{0,3}=\Z\oplus\Z \to E_2^{4,0}=H^4(B_{\Z_2}\times
B_{\Z_4},\Z)=\Z_2\oplus\Z_4\oplus\Z_2$ whose kernel is equal to the
image of the edge homomorphism, which can be viewed as $\eta^*$.
Clearly, $(-2,2)$ and $(2,2)$ lie in this kernel and must be a basis
of the lattice since the spectral sequence also implies that its
quotient group has order 8.
 In this
basis,  the matrix of $\pi_-^* - \pi_+^*$ on the free part is given
by:
$$\left(\begin{array}{rr} \frac{c}{4}(p_-^2+q_-^2) & -\frac{d}{4}(p_+^2+q_+^2) \\
\frac{c}{4}(p_-^2-q_-^2) & -\frac{d}{4}(p_+^2-q_+^2)
\end{array}\right) $$
\no where $c,d$ are the integers from   Lemma \ref{orbits}, which we
showed are powers of two. We now claim that $c=2$ and $d=4$, which
then implies (\ref{mfdN}) as before.

If this were not the  case, the order of $H^4(N,Z)$ would be even
since $(p_-^2+q_-^2,p_-^2-q_-^2)=2$ and
$(p_+^2+q_+^2,p_+^2-q_+^2)=1$, and we now show that it must in fact
be odd. To see this, we repeat the above Mayer Vietoris argument for
cohomology with $\Z_2$ coefficients. First observe that in the Gysin
sequence of $\S^1\to \G/\H\to B_\pm$ one has
$$\cdots \to H^3(B_\pm)=
\Z_2^2\overset{\pi_\pm^*}{\longrightarrow} H^3(\G/\H)=\Z_2^4\to
H^2(B_\pm)=\Z_2^2\to \cdots$$

\no  and hence $\pi_\pm^*$ are injective. Thus from the Mayer
Vietoris sequence
 $$0\to H^3(N)\to H^3(B_-)\oplus
H^3(B_+)=\Z_2^4 \overset{\pi_-^* - \pi_+^*}{\Lllongrightarrow}
H^3(\G/\H)=\Z_2^4 \to H^4(N)\to 0
$$
\no it follows that $H^3(N,\Z_2)=H^4(N,\Z_2)=0$ which completes the
proof.
\end{proof}

\bigskip


\section{Appendix I: Classification in Even Dimensions.}

In this appendix we give a relatively short proof of  Verdiani's theorem based on the obstructions, ideas and strategy presented here to handle the odd dimensional case.

\begin{thm}[Verdiani]
  Suppose $\G$ acts on an even dimensional positively curved simply
  connected
  compact manifold $M$ with cohomogeneity one. Then $M$ is equivariantly
  diffeomorphic to a rank 1 symmetric space.
\end{thm}

One of the reasons that make the even dimensional case less involved
is that one of the groups $\Km$ or $\Kp$ has the same rank as $\G$,
and thus $\rank(\G)-\rank(\H)=1$ as we saw in the Rank Lemma. Moreover, the
Upper Weyl Group Bound now says that $|W|=2$ or $4$, and $|W|=2$ if
$\H$ is connected and $l_-$ and $l_+$ are both odd. This becomes
especially powerful if combined with the Lower Weyl Group Bound.
 Another noteworthy difference is that the main body of work is
 confined to the case of simple groups, and that induction is
 only used occasionally.

In  \cite{iwata:1,iwata:2,uchida}  it was shown that any \com whose
rational cohomology ring is like that of $\CP^n \, , \, \HP^n$ or
$\CaP$, is equivariantly diffeomorphic to one of the known linear
\co one actions. Hence it is again sufficient to recover $M$ up to
homotopy type.

\bigskip

We treat the cases $\G$ simple, or not separately. For $\G$ simple
we distinguish among the subcases:  $\H$ contains a normal simple
subgroup of $\rank\ge 2$, or not.

\bigskip

\begin{center}
$\G$ is not simple.
\end{center}

\begin{prop}\label{evennotsimple}
If $\G$ is not simple and acts essentially with corank one, it is the
tensor product action of $\SU(2)\SU(k)$ on $\CP^{2k-1}$.
\end{prop}

\begin{proof}

  We can assume that $\G=\L_1\times \L_2$, and say $\rank(\Km)=\rank(\G)$, and
hence $\Kmo=\K_1\cdot \K_2$. Since by assumption $\Km$ does not
contain a normal subgroup of
$\G$, we see that $\G$ is semisimple and
  $\G/\Kmo$ is necessarily isometric to
a product. If each of the factors has dimension $>2$, then we can
derive a contradiction as in the proof of the Product Lemma in odd
dimensions. Using the conventions therein, we may assume that
$\L_2/\K_2=\Sph^2$, and that $\K_1$ acts transitively on the normal
sphere. It follows again that $\L_1/\K_1$ is isometric to a rank one
projective space, and as before, we derive a contradiction  if the
isotropy representation of $\L_1/\K_1$ is of real type. Hence the
isotropy representation is of complex type and $\L_1/\K_1=
\SU(k)/\U(k-1)$, $k\ge 2$ or $\L_1/\K_1= \G_2/\SU(3)$.

Because of primitivity $\Kp$ necessarily projects surjectively onto
$\L_2=\S^3$. The projection of $\Ho$ onto $\L_2$ cannot be 3
dimensional since then the subaction of $\L_1$ would be orbit equivalent to
the $\G$ action. If the projection is trivial,
  $\Km/\H=\Sph^1$ and
$\Kp/\H=\Sph^3$. Furthermore, $\H$ is connected since $\K_1$ and
$\K_2$ are both connected. But then the Upper Weyl Group Bound
implies that $|W|=2$, which combined with the Lower Weyl Group
Bound, gives $\dim \G/\H \le 4$, a contradiction. Hence the
projection is one dimensional and $\Kp/\H=\Sph^2$.

This completely determines the group picture. Indeed, if $\L_1/\K_1=
\SU(k)/\U(k-1)$, we have $\Kp=\L\times \S^3$ with the second factor
embedded diagonally in $\L_1\cdot\L_2$. Hence  $\L\subset \U(k-2)$
and $\H=\L\times\S^1$. In order for $\Km/\H$ to be a sphere, we need
$\L=\U(k-2)$ and we recover the tensor product action of
$\SU(k)\times\SU(2)$ on $\CP^{2k-1}$.

  In the case of $\L_1\cdot\L_2=\G_2\times\S^3$ we have
  $\Kp=\S^3\times\S^3$. $\Kp$ projects onto   $\SO(4)\subset\G_2$,
   and the second factor of $\Kp$ projects onto $\L_2$.
The tangentspace $T^+$ of the orbit $B_+=\G_2 x \S^3/\gK_+$
decomposes as an $8$ dimensional and a $3$ dimensional invariant
 subspace.   The natural representation in $S^2T^+$ splits into
two trivial two $5$-dimensional
and into subrepresentations on which $\G_2\cap \gK_+\cong \S^3\subset \H$
 acts nontrivially.

This in turn implies that $B_+$ is totally geodesic, a contradiction.
\end{proof}

\smallskip

\begin{center}
$\G$ is simple.
\end{center}

\begin{prop}
Assume  $\G$ is simple with corank 1  and all simple factors in $\H$
have rank one.  If the action is essential, then $(M,\G)$ is one of
the following pairs: $(\CP^6, \G_2), (\CP^9,\SU(5))$,
$(\HP^2,\SU(3)), (\HP^3,\SU(4))$, or $ (\CaP,\Sp(3))$, and the
actions are given by Table \ref{even}
\end{prop}

\begin{proof}

Since all normal factors in $\Ho$ are one and three dimensional, we
have $l_\pm = 1,2,3,4,5$, or $7$ and at least one of them is odd.
The Weyl group has order at most $4$, and the order is $2$ if $l_+$
and $l_-$ are both one of $3,5$ or $7$.

We first treat the case where $\rank\G\le 2$. For $\G=\S^1$ clearly
$M=\Sph^2$. If  $\G=\S^3$,  $\Kpm\cong\S^3$ corresponds to
$M=\Sph^4$, $\Km\cong\S^3, \Kp\cong\S^1$ corresponds to $M=\CP^2$,
and $\Kpm\cong \S^1$ is not primitive. In all three cases, the
action is not essential.

If $\G=\SU(3)$, then $\H$ cannot be three dimensional since
otherwise $\G$ is not primitive or has a fixed point. Therefore
$\Ho$ is a circle. The Lower Weyl Group Bound  implies that one of
$l_\pm$, say $l_-$  is 3 and hence $\Kmo=\U(2)$. Since $\U(2)$ is
maximal, $\Kp$ and hence $\H$ are connected. From the Upper and
Lower Weyl Group Bound it now follows that $l_+=1$ or $3$ is not
possible.
 Thus
$\Kp=\SO(3)$ or $\SU(2)$. If $\Kp=\SO(3)$ and hence $\Sph_+ =
\Sph^2$, we note there is only one embedding into $\SU(3)$ and its
isotropy representation is irreducible. One then easily shows
(Clebsch Gordon formula) that the representation of $\SO(3)$ in $S^2
T_+$ has no three dimensional subrepresentation, which implies that
$\SU(3)/\SO(3)$ is totally geodesic, a contradiction. If
$\Kp=\SU(2)$, we have recovered the group picture for $\HP^2$.

If  $\G=\Symp(2)$ and $\H$ is three dimensional, then Table
\ref{sph} implies that $\G/\Ho\cong \Sph^7$ or $\SO(5)/\SO(3)$. In
the former case $\G$ either has a fixed point or is not primitive.
In the latter case the Chain Theorem applies. So we may assume
$\dim(\H)=1$. The Lower Weyl Group Bound implies that $l_-=2\, ,\,
l_+=3$ and hence all groups are connected. If $\H = \{\diag(z^k,z^l)\mid z\in \gS^1\}$,
the isotropy representation has weights $2k, 2l$, and $2k \pm 2l$.
But since $\Symp(2)/\H$ by the Isotropy Lemma can have at most two
nonequivalent nontrivial subrepresentations,
  it follows that  $\H = \{\diag(z,z)\mid z\in\gS^1\}$ or
 $\{\diag(1,z)\mid z\in\gS^1\}$.
 In the former case $M^{\H}_c$ admits a cohomogeneity
one action of $\N(\H)/\H\cong \SO(3)$ with trivial principal isotropy
group, which contradicts the Core-Weyl Lemma. If
$\H=\{\diag(1,z)\}$,  we may assume that $\Km=\{\diag(1,g)\mid g\in
\S^3\}$.  Hence $\Km$ is normalized by the normalizer of $\H$ and in
particular by the Weyl group. By Linear Primitivity, the Lie
algebras of $\Km,\Kp, w_+ \Kp w_+$ span the Lie algebra of $\G$,
which is not possible since $\dim \Kp=4$.

Finally, if  $\G=\G_2$ and $\H$ is one dimensional, we obtain a
contradiction to the Lower Weyl Group Bound since  $l_\pm\le 3$. The
only three dimensional spherical subgroup is $\H=\SU(2)\subset
\SU(3)\subset \G_2$, as one easily verifies. The subgroups of $\G_2$
only allow $\l_\pm=1$ or $5$, and using the Lower Weyl Group Bound,
we see that $\Kmo=\SU(3)$, $\Kpo=\SU(2)\S^1$ and $\H$ is not
connected by the Upper Weyl Group Bound. Because of
$\G_2/\SU(3)\cong \Sph^6$,
  $\Km$  and thereby $\H$ has at most two, and hence two
  components. Now all groups $\Kpm$, $\H$
  are determined, and we  have recovered the picture of $\CP^6$ endowed
with the linear action of $\G_2$.

If the rank of $\G$ is 3 or larger, we first observe that  $\dim \H
\le 3\rank \H = 3(\rank \G -1)$ since all simple factors of $\H$
have rank one. Hence:

$$\dim(\G)-3\rank \G \le \dim(\G/\H)-3$$

\no By the Lower Weyl Group Bound
\[
\dim(\G/\H)\le 2(7+4)= 22
\]
and hence $\dim \G -3\rank(\G) \le  19 $.

First we consider the case that there is  an orbit, say $\G/\Km$,
  of codimension $8$. Then $\Symp(2)$ is a normal subgroup of $\Km$
and $\rank(\Km)=\rank(\G)$ since $\rank\Km - \rank \H =1$. The only
simple Lie groups satisfying the above dimension estimate and
containing $\Symp(2)$ as a regular subgroup, besides $\Symp(2)$
itself, are $\Spin(7)$ and  $\Symp(3)$. In case of $\G=\Spin(7)$ the
central element in $\Symp(2)$ is necessarily central in $\Spin(7)$,
but does not lie in $\H$. But then $\Spin(7)/\Km$ is totally
geodesic, which is not possible. In the case of $\G=\Symp(3)$, $\H$
contains an $\Symp(1)$-block and $M^{\Sp(1)}_c$ admits a
cohomogeneity one action by $\Symp(2)$ whose
  principal isotropy group has rank one.
As we saw earlier, this isotropy group must be three dimensional and
hence $\H$ is $6$-dimensional, which implies $\Kmo=\Symp(1)\cdot
\Symp(2)$. Since this group is maximal in $\Sp(3)$, it follows that
$\Km$ and hence $\H$ are connected. The Lower Weyl Group Bound
implies that $|W|=4$ and hence, by the Upper Weyl Group Bound, one
of the codimensions must be odd. Hence $\Kp=\Symp(2)$,
$\H=\Symp(1)^2\subset \Symp(2)\subset \Symp(3)$ and we have
   recovered the cohomogeneity one action of $\Symp(3)$ on $\CaP$.

If there are no orbits of codimension $8$,
\[
\dim(\G/\H)\le 2(4+5)=18.
\]
and hence $\dim(\G) -3\rank(\G)\le 15$. We now assume that there is
an orbit, say $\G/\Km$
  of codimension $5$.
Then $\Symp(2)\subset \Km$. The only groups, other than $\Symp(2)$,
  satisfying the above dimension estimate
and containing $\Symp(2)$ are $\Symp(3)$,  $\SU(5)$, $\Spin(7)$ and
$\Spin(6)$.

In the case of $\G=\Spin(6)$ or $\Spin(7)$, there is a unique
embedding of $\Sp(2)=\Spin(5)$ in $\G$  and hence $\H$ contains a
$4\times 4$-block and the Chain Theorem applies. In the case of
$\G=\Symp(3)$, $\H=\Symp(1)^2\subset \Symp(2)$ and  the isotropy
representation of $\G/\H$ contains a four dimensional
representation, which can only degenerate in an orbit of codimension
$8$, which we already dealt with. Thus we may assume $\G=\SU(5)$.
Because of $\rank(\Km)=3$, we have $\Kmo=\Symp(2)\cdot \S^1$ and
hence $\Ho=\SU(2)^2\cdot \S^1$. There is a four dimensional
subrepresentation of the isotropy representation of $\G/\Ho$ which
is not equivalent to a subrepresentation of $\Km/\Ho$. From the
Isotropy Lemma we deduce $\Kp=\SU(3)\cdot \SU(2)\cdot \S^1$ and all
groups are connected. We have recovered the linear action of
$\SU(5)$ on $\CP^{9}$.

We are left with the case that $l_\pm = 1,2,3$ or $5$. Hence
\[
\dim(\G/\H)\le 2(2+5)=14.
\]
and $\dim(\G)-3\rank(\G)\le 11 $. The only  simply connected compact
simple Lie
  groups satisfying this dimension
estimate are  $\S^3$, $\SU(3)$, $\Symp(2)$, $\G_2$, and $\SU(4)$.

This only leaves us to consider the case  $\G=\SU(4)$. If $\Ho$ is
abelian, we obtain a contradiction to the Lower Weyl Group Bound.
There are two six dimensional subgroups of $\SU(4)=\Spin(6)$, one
from $\SO(4)\subset\SO(6)$, where the Chain Theorem applies, and the
other from $\SO(3)\SO(3)\subset\SO(6)$ which contradicts the
Isotropy Lemma.

In the case of $\dim(\Ho)=4$, Table \ref{sph} implies that $\H$
contains an $\SU(2)$-block. Since the four dimensional
representation of $\SU(2)$ has to degenerate, $\Km=\U(3)$. Since
$\U(3)$ is maximal in $\SU(4)$, $\Km$ and hence also $\H$ are
connected. If $|W|=2$, the Lower Weyl Group Bound implies that
$l_+=7$ which is not possible since $\Sp(2)$ is maximal in $\SU(4)$.
Now the Upper Weyl Group Bound implies that $l_+$ is even. Thus
$l_+=2$,  $\Kp=\SU(2)^2, \H=\S^1\SU(2)$, and we have recovered the
action of $\SU(4)$ on $\HP^3$.
\end{proof}

\smallskip

For $\G$ simple it remains to consider the case where $\H$ has a higher rank normal subgroup.

\begin{prop}
Assume  $\G$ is simple with corank 1  and $\H$ contains a simple
subgroup of rank $\ge 2$. If the action is essential, the pair $(M,\G)$  is
one of the following: $(\CP^{n-1},\SO(n))$, $(\HP^{n-1},\SU(n))$,
$(\CP^{15},\Spin(10))$, $(\Sph^{14},\Spin(7))$, or $(\CP^7,\Spin(7))$ with the
actions given by Table \ref{even}.
\end{prop}

\begin{proof}

  By Lemma \ref{Hrank2}, there can be only one connected normal subgroup of $\H$
  which has rank larger than one, which we denote by $\H'$.

\smallskip

\begin{center}
  $\G=\Sp(k)$ or $\SU(k)$
\end{center}

If $\G=\Symp(k)$, Table \ref{sph} implies that $\H'$ is given by an
$h\times h$ block $h\ge 2$, and the Chain Theorem applies.

If $\G=\SU(k)$, Table \ref{sph} implies that either
  $\H'$ is given by an $h\times h$ block and the Chain Theorem applies, or
$\H'=\Symp(2)\subset \SU(4)\subset \SU(k)$.
  The latter case can be ruled out  as follows. If $k=4$, then
clearly $\G$ has a fixed point. If $k\ge 5$,  there is an eight
dimensional irreducible representation of $\H$ which can only
degenerate in an isotropy group $\Km$ containing $\Symp(3)$, and
furthermore $\rank(\Km)=\rank(\G)$. But this is impossible since
$\Symp(3)$ is not a regular subgroup of $\SU(k)$.

\smallskip

\begin{center}
  $\G=\Spin(k)$
\end{center}

If $\G=\Spin(k)$, then by Table \ref{sph}, either $\H'$ is given by
a block and
  the Chain Theorem applies,
or $\H'\cong \G_2,\Spin(7),\SU(4), \Symp(2)$, or $\SU(3)$.

 If
$\H'=\G_2$ or $\Spin(7)$, we first claim  that $\Ho=\H'$. Indeed, if
$\Ho\ne\H'$,   it follows from Table \ref{trans} that only one of
$\H'$ or $\Ho/\H'$ can act non-trivially on each irreducible
subrepresentation of $\G/\H$, which as in the proof of Lemma
\ref{Hrank2} contradicts the assumption that $\G$ is simple.
 If
$\Ho=\G_2$, the Rank Lemma implies $\G=\Spin(7)$ and $\G$ has a
fixed point. If $\Ho=\Spin(7)$, then $\G=\Spin(8)$ or $\G=\Spin(9)$
and  $\G/\Ho$ is a sphere. Then $\G$ either has a fixed point or the
action is not primitive.

If $\H'=\SU(4)$  or $\Symp(2)$, then $k\ge 8$. If $k=8$,  $\H'$ is,
up to an outer automorphism of $\Spin(8)$, given as a $6\times 6$ or
a $5\times 5$ block and the Chain Theorem applies. If $k\ge 9$, let
$\iota$ be the order 2 central element in $\H'$ so that
$\N(\iota)_{\subo}=\Spin(k-8)\cdot \Spin(8)$ acts with cohomogeneity one
on the reduction $M^{\iota}_c$ and the principal isotropy group
contains, up to outer automorphisms, a $5\times 5$ or a $6\times 6$
block. By induction it must be induced by
a tensor product action,
$\H'=\SU(4)$, $k=10$ and  $\Ho=\SU(4)\cdot\S^1$ by the Rank Lemma.
Hence $\Kmo=\SU(5)\cdot\S^1$ since the 8 dimensional
  representation has to degenerate.
  The isotropy representation of $\SO(10)/\U(5)$, when restricted to
  $\U(4)$,
  contains a six dimensional representation, which has to degenerate
  in
    $\Kp/\H$  and
hence $\Kp=\Spin(7)\cdot\S^1$. We have recovered the action of
$\Spin(10)$ on $\CP^{15}$.

Finally, we consider the case $\H'=\SU(3)\subset \Spin(6)\subset
\Spin(k) $. We first rule out $k\ge 8$. In this case $\rank(\H)\ge
3$, and hence there exists an irreducible summand in the isotropy
representation of $\G/\H$ on which both $\H'$ and $\H/\H'$ acts
nontrivially (cf.\eqref{Hrank2}). Thus not all the six dimensional
representations of $\SU(3)$ can degenerate in $\G_2$. Thus an
isotropy group, say $\Km$, contains $\SU(4)$ as a normal subgroup,
and we consider the fixed point set of the central involution $\iota
\in \SU(4)$. Since it is central in a $\Spin(8)$-block and acts as
$-\id$ on the slice, it has
  a homogeneous fixed point component $\Spin(k-8)\cdot
  \Spin(8)/( \Km\cap\Spin(k-8)\cdot
  \Spin(8))$ which cannot have positive curvature.

  In the case of $k=6$, either $\G$ has
a fixed point, or the action is not primitive. Thus we may assume $k
= 7$ and hence $\Ho=\SU(3)$. The isotropy representation of
$\Spin(7)/\SU(3)$ consists of the sum of a trivial representation, a
6 dimensional representation corresponding to the isotropy
representation in $\SU(4)=\Spin(6)$, and a second 6 dimensional
representation orthogonal to it. Thus the only connected subgroups
in between $\SU(3)$ and $\Spin(7)$ are $\U(3)$, $\Spin(6)$ and
$\G_2$, and the normalizer $\N(\Ho)/\Ho\cong\S^1$ acts transitively
on the possible embeddings of $\Spin(6)$ and $\G_2$.

If $\Kmo=\SU(4)\cong \Spin(6)$ occurs as isotropy group,
then $\Kpo=\S^1\cdot\SU(3)$ or $\Kpo=\Spin(6)$ and the action is not primitive.
Thus
$\H$ is connected, $\Kp\cong \G_2$, and we have recovered the action
on $\Sph^{14}$.

If $\SU(4)$ does not occur as isotropy group, primitivity implies
  that $\Kmo=\G_2$,
and  $\Kpo= \S^1\cdot \SU(3)$. As the center of $\Spin(7)$ is
contained in $\Kp$, it must be contained in $\H$ also, since
otherwise $\Spin(7)/\Kp$ would be totally geodesic. This also shows that
$\Spin(7)/\Km\cong \RP^6$, and $\Km$ and $\H$ have precisely two
components. We have recovered the linear action of $\Spin(7)$ on
$\CP^7$.

\smallskip

\begin{center}
  $\G=\F_4 \, , \, \E_6 \, , \, \E_7 \, , \, \E_8$
\end{center}

If $\G$ is one of $\F_4$, $\E_6$, $\E_7$, or $\E_8$,  Table
\ref{sph} implies that $\H'$ is one of the groups $\Spin(k)$, $k\le
8$, $\G_2$ or $\SU(3)$, where we again used the fact that
$\Ho=\Spin(9)$ is not possible. If $\H'\neq \Spin(7)$, we have
$\dim(\K_\pm/\H)\le 8$ and hence $\dim(\G/\H)\le 32$ by the Lower
Weyl Group Bound. This implies that $\dim \G - 3\rank \G \le 29$
which is clearly not possible.

  For $\H'=\Spin(7)$, it follows as before that $\H'=\Ho$
  and thereby
$\G=\F_4$. In one singular orbit the $8$-dimensional representation
of $\Spin(7)$ has to degenerate. This implies $\Kmo=\Spin(9)\subset
\F_4$ and since $\Spin(9)$ is maximal in $\F_4$, $\Km$ and $\H$ are
connected. Since $l_\pm$ are one of $1, 7,15$, the Upper Weyl Group
Bound implies that $|\W|=2$, which contradicts the Lower Weyl Group
Bound.
\end{proof}

\bigskip


\section{Appendix II: Group Diagrams for Compact Rank One Symmetric Spaces.}

In this Appendix we will collect various known information that will be
used throughout the proof of Theorem A.

To describe the representations, we use the notation $\rho_n\, , \,
\mu_n\, , \, \nu_n$ for the defining representations of $\SO(n)$,
 $\U(n)\, , \, \Sp(n)$ respectively. $\Delta_n$ denotes the spin
representation for $\SO(n)$ and $\Delta_{2n}^\pm$ the half spin
representation. Also $\phi$ denotes a 2 dimensional representation of
$\S^1$ and for all others we use $\psi_N$ for an N-dimensional
irreducible representation.

\smallskip

In Table \ref{sph} we reproduce the list of spherical simple
subgroups of the simple Lie groups from \cite[Proposition
7.2-7.4]{Wi: sym} since it will be an important tool in our
classification. All embeddings are standard embeddings among
classical groups, and $\Spin(7)\subset \SO(8)$ is the embedding via
the spin representation. We point out that the case of a rank one
group in the exceptional Lie groups was not included in \cite{Wi:
sym}. But in our proof, this case will only be needed for a rank one
group in $\G_2$, where one easily shows that  $\SU(2)\subset\SU(3)$
is the only spherical subgroup.

\renewcommand{\thetable}{\Alph{table}}

\renewcommand{\arraystretch}{1.4}
\begin{table}[!h]
      \begin{center}
          \begin{tabular}{|c|c|c|}
\hline

$\G$          &$\H$            & Inclusions \\
\hline \hline

$\SU(n)$ &$\SU(2)$ &$\SU(2)\subset\SU(n)$ given by $p\,(\mu_2)\oplus q\, id$ \\
\hline

$\SU(n)$ &$\Sp(2)$ &$\Sp(2)\subset\SU(4)\subset\SU(n)$\\
\hline

$\SU(n)$ &$\SU(k)$ & $k\times k$ block
       \\
\hline\hline

$\SO(n)$ &$\Sp(1)$ &$\Sp(1)\subset\SO(n)$ given by $p\,\nu_1\oplus q\, id$ \\
\hline

$\SO(n)$ &$\SU(3)$ &$\SU(3)\subset\SO(6)\subset\SO(n)$\\
\hline

$\SO(n)$ &$\Sp(2)$ &$\Sp(2)\subset\SU(4)\subset\SO(8)\subset\SO(n)$\\
\hline

$\SO(n)$ &$\G_2$ &$\G_2\subset\SO(7)\subset\SO(n)$\\
\hline

$\SO(n)$ &$\SU(4)$ &$\SU(4)\subset\SO(8)\subset\SO(n)$\\
\hline

$\SO(n)$ &$\Spin(7)$ &$\Spin(7)\subset\SO(8)\subset\SO(n)$\\
\hline

$\SO(n)$ &$\SO(k)$ & $k\times k$ block \\
\hline\hline

$\Sp(n)$ &$\Sp(1)$ &$\Sp(1)=\{\diag(q,q,\cdots ,q,1,\cdots ,1)\mid
q\in \S^3\}$ \\
\hline

$\Sp(n)$ &$\Sp(k)$ & $k\times k$ block\\
\hline\hline

$\G_2$ &$\SU(3)$ & maximal subgroup \\
\hline

$\F_4 , \E_6$ &$\Spin(k)$
&$\H\subset\Spin(9)\subset\F_4\subset\E_6\subset\E_7
\subset\E_8$  \\

$\E_7 , \E_8$& & $\H=\Spin(k)\, , k= 5 ,\cdots , 9$ standard embedding\\

\hline

$\F_4 , \E_6\cdots\E_8$ &$\SU(3)$
&$\SU(3)\subset\SU(4)\subset\Spin(8)\subset\F_4\subset\E_6\subset\E_7
\subset\E_8$  \\
\hline

$\F_4 , \E_6\cdots\E_8$ &$\G_2$
&$\G_2\subset\Spin(7)\subset\Spin(8)\subset\F_4\subset\E_6\subset\E_7
\subset\E_8$  \\
\hline

\hline

          \end{tabular}
      \end{center}
      \vspace{0.1cm}
      \caption{$\G/\H$ with spherical isotropy representations}\label{sph}
\end{table}

\smallskip

In Table \ref{trans} we list the transitive actions on spheres and
their isotropy representation.  Notice that $\nu_n\hat{\otimes}\nu_1
$  is the representation on ${\QH}^n=\R^{4n}$ given by left
multiplication of $\Sp(n)$ and right multiplication of $\Sp(1)$ on
quaternionic vectors and $\nu_n\hat{\otimes}\phi$ is the
subrepresentation under $\U(1)\subset\Sp(1)$. Notice also  that for
each irreducible subrepresentation $\fm$ in the isotropy
representation of $\K/\H$ the group $\H$ acts transitively on the
unit sphere in $\fm$, as long as $\dim \fm > 1$. This elementary but
important property is used in Isotropy \lref{isotropy}.

\smallskip

In Table \ref{homog} we list the  remaining simply connected homogeneous spaces with
positive curvature which will be used when one needs to check wether
a singular orbit can be totally geodesic.

\smallskip

  \stepcounter{equation}
\begin{table}[!h]
      \begin{center}
          \begin{tabular}{|c|c|c|c|}
\hline n &$\K$          &$\H$            & Isotropy representation \\
\hline \hline

n&$\SO(n+1)$ &$\SO(n)$  &$\rho_n$\\
\hline

$2n+1$&$\SU(n+1)$ & $\SU(n)$  &$\mu_n\oplus id$   \\
\hline

$2n+1$&$\U(n+1)$ & $\U(n)$  &$\mu_n\oplus id$   \\
\hline

$4n+3$&$\Sp(n+1)$ & $\Sp(n)$  &$\nu_n\oplus 3\, id$   \\
\hline

$4n+3$&$\Sp(n+1)\Sp(1)$ & $\Sp(n)\Delta\Sp(1)$
&$\nu_n\hat{\otimes}\nu_1 \oplus
id \,\hat{\otimes}\rho_3  $   \\
\hline

$4n+3$&$\Sp(n+1)\U(1)$ & $\Sp(n)\Delta\U(1)$
&$\nu_n\hat{\otimes}\phi \oplus
id\,\hat{\otimes}\phi \oplus id  $   \\
\hline

$15$&$\Spin(9)$   &$\Spin(7)$  &$\rho_7\oplus\Delta_8 $   \\
\hline

$7$&$\Spin(7)$   &$\G_2$  &$\phi_7 $   \\
\hline

$6$&$\G_2$   &$\SU(3)$  &$\mu_3 $   \\
\hline

\hline
          \end{tabular}
      \end{center}
      \vspace{0.1cm}
      \caption{Transitive actions on $\Sph^n$}\label{trans}
\end{table}

\stepcounter{equation}
\begin{table}[!h]
      \begin{center}
          \begin{tabular}{|c|c|c|}
\hline n &$\G$          &$\K$   \\
\hline \hline

2n&$\SU(n+1)$ &$\U(n)$  \\
\hline
4n&$\Sp(n+1)$ &$\Sp(n)\Sp(1)$  \\
\hline

4n&$\Sp(n+1)$ &$\Sp(n)\U(1)$  \\
\hline

16&$\F_4$ &$\Spin(9)$  \\
\hline

6&$\SU(3)$ &$\T^2$  \\
\hline

12&$\Sp(3)$ &$\Sp(1)^3$  \\
\hline

24&$\F_4$ &$\Spin(8)$  \\
\hline

7&$\SU(3)$ &$\S^1=\diag(z^p,z^q,\bar{z}^{p+q})$  \\
& & $(p,q)=1$, $pq(p+q)\ne 0$\\ \hline

7&$\U(3)$ & $\T^2$ \\
\hline

13&$\SU(5)$ &$\Sp(2)\cdot\S^1$  \\
\hline

\hline
          \end{tabular}
      \end{center}
      \vspace{0.1cm}
      \caption{Homogeneous spaces $M^n=\G/\K$ with positive
curvature,  which are not spheres}\label{homog}
\end{table}

Information about \co one actions on spheres is taken from
\cite{straume} where the group $\G$ and the principal isotropy group
$\H$ are given. The groups $\Kpm$, which can not all be found in the
literature, amusingly are obtained along the way in our proof. In
other words,  once an essential action arises in the induction proof
with  $\G$ and $\H$ from Straume's list, our obstructions leave only
one possibility for $\Kpm$. And all essential actions indeed arise
in the proof along the way.

\smallskip

In Table \ref{odd} we describe the essential group actions on odd
dimensional spheres and in Table \ref{even} the ones on even
dimensional rank one symmetric spaces. In these Tables $k$ is an
integer larger or equal to one. We also include the normal
extensions since these extensions will be used in the induction
proof. The \co one actions on $\CP^n$ and $\HP^n$ are obtained from
an action on an odd dimensional sphere when $\U(1)$ or $\Sp(1)$ is a
normal subgroup in $\G$ with induced action given by a Hopf action.
Conversely, an action on $\CP^n$ or $\HP^n$ lifts to such an action
on a sphere.

\smallskip

We will also use some knowledge about non-essential actions (apart
from the extensions of essential ones). One easily sees that if an
action of $\G$ on a sphere in $\R^n$ is not essential, and no
normal subaction is essential, then
$\G=\L_1\L_2$ and $\R^n=V_1\oplus V_2$ such that $\G$ preserves
$V_i$ and  $\L_i$ acts transitively on the unit sphere in $V_i$. The
most elementary ones are the sum actions by $\G(n)\G'(m)$ operating
on $V^n\oplus V^m$ as $f_n\hat{\otimes}id\oplus id\hat{\otimes}f_m $
where $\G(n)$ is any of the classical Lie groups with their defining
representations $f_n$. The property we sometimes use is that the
principal isotropy group is given by $\G(n-1)\G'(m-1)$.
This also
includes the case where one $\G(n)$ is absent, which corresponds to
actions with a fixed point.  Such sum actions can be further
modified by replacing the action of $\G(n)$ on $V^n$ by any one of
the other transitive actions on spheres. The corresponding isotropy
groups are given in Table \ref{trans}.

In these simplest kind of sum actions, each simple normal subgroup
of $\G$ acts nontrivially on only one of the subspaces $V^n , V^m$.
One can modify them  further by considering actions where some of
the factors (necessarily of rank one) operate on both. This is only
possible if the rank one factor commutes with both transitive
actions on $V_i$. Hence there are three such actions:
\begin{enumerate}
  \item $(\G,\H)$ is given by $(\U(1)\G(n)\G(m)\, , \,
\Delta\U(1)\G(n-1)\G'(m-1) )$
  acting via
$\phi^k\hat{\otimes}f_n\hat{\otimes}id\oplus
\phi^l\hat{\otimes}id\hat{\otimes}f_m$ with $(k,l)=1$ and $\G(n)$ is
one of $\SU(n)$ or $\Sp(n)$.   If one of the groups $\G(n)$ is
absent, the principal isotropy group is $\Delta\Z_k\G(n-1)$.

\item $(\G,\H)$ is given by
$(\Sp(1)\Sp(n)\Sp(m),\Delta\Sp(1)\Sp(n-1)\Sp(m-1))$
  acting via
$\nu_1\hat{\otimes}\nu_n\hat{\otimes}id\oplus
\nu_1\hat{\otimes}id\hat{\otimes}\nu_m$, including the case where
$\Sp(m)$ is absent.

\item  $(\G,\H)$ is given by $(\Sp(1)\Sp(n),\Delta\U(1)\Sp(n-1))$ acting via
$\nu_1\hat{\otimes}\nu_n\oplus \rho_3\hat{\otimes}id$, which is an
action on an even dimensional sphere.

\end{enumerate}

For the even dimensional rank one symmetric spaces one also has sum
actions  by $\Sp(n)\Sp(m)$ on $\HP^{n+m-1}$ and by $\SU(n)\SU(m)$ or
$\S(\U(n)\U(m))$  on $\CP^{n+m-1}$, where one or both unitary groups
can also be replaced by symplectic groups. If one of the groups is
absent, they are the actions with a fixed point.

\smallskip

Finally, in  Table \ref{symm}  we list the symmetric spaces $\G/\K$
where $\K$ and $\G$ have the same rank. They occur as normalizers of
elements $\iota$ whose square, but not $\iota$ itself, lies in the
center of $\G$. We point out that in this Table the group
$\Spin(4k)/\Z_2$ is the quotient of $\Spin(4k)$ which is not
isomorphic to $\SO(4k)$.

{\setlength{\tabcolsep}{0.08cm}
\renewcommand{\arraystretch}{1.6}
\stepcounter{equation}
\begin{table}[!h]
      \begin{center}
          \begin{tabular}{|c||c|c|c|c|c|c|c|}
\hline
    n     &$\G$ &$\chi$   &$\Km$    &$\Kp$     &$\H$& $(l_-,l_+)$ & $W$   \\
\hline \hline

$8k+7$ & $\Sp(2)\Sp(k+1)$ & $\nu_2\hat{\otimes}\nu_{k+1}$ &
$\triangle\Sp(2)\Sp(k-1)$ & $\Sp(1)^2\Sp(k)$ & $\Sp(1)^2\Sp(k-1)$ &
$( 4,4k+1 )$ & $\D_4$
\\

\hline

$4k+7$ & $\SU(2)\SU(k+2)$ & $\mu_2\hat{\otimes}\mu_{k+2}$ &
$\triangle\SU(2)\SU(k)$ & $\S^1\cdot\SU(k+1)$ & $\S^1\cdot\SU(k)$ &
$(2 ,2k+1 )$ & $\D_4$
  \\

 & $\U(2)\SU(k+2)$ & $\mu_2\hat{\otimes}\mu_{k+2}$ &
$\triangle\U(2)\SU(k)$ & $\T^2\cdot\SU(k+1)$ & $\T^2\cdot\SU(k)$ &  & \\

$7$ & $\U(2)\SU(2)$ & $\mu_2\hat{\otimes}\mu_2$ &
$\triangle\SU(2)$ & $\T^2$ & $\S^1$ & $(2,1)$ & \\

\hline

\hline

$2k+3$ & $\SO(2)\SO(k+2)$ & $\rho_2\hat{\otimes}\rho_{k+2}$ &
$\triangle\SO(2)\SO(k)$ & $\Z_2\cdot\SO(k+1)$ & $\Z_2\cdot\SO(k)$
& $( 1,k )$ & $\D_4$\\

\hline

$15$ & $\SO(2)\Spin(7)$ & $\rho_2\hat{\otimes}\Delta_7$ &
$\triangle\SO(2)\SU(3)$ & $\Z_2\cdot\Spin(6)$ & $\Z_2\cdot\SU(3)$ &
$(1 ,7 )$ & $\D_4$ \\
\hline

$13$ & $\SO(2)\G_2$ & $\rho_2\hat{\otimes}\phi_7$ &
$\triangle\SO(2)\SU(2)$ & $\Z_2\cdot\SU(3)$ & $\Z_2\cdot\SU(2)$ & $(
1,5 )$ & $\D_4$\\
\hline

$7$ & $\SO(4)$ & $\nu_1\hat\otimes\nu_3$ &
$\S(\O(2)\O(1))$ & $\S(\O(1)\O(2))$ & $\Z_2\oplus\Z_2$ & $(1 ,1 )$ & $\D_6$ \\
\hline

$15$ & $\Spin(8)$ & $\rho_8\oplus\Delta_8^\pm$ &
$\Spin(7)$ & $\Spin(7)$ & $\G_2$ & $(7 ,7 )$ & $\D_2$ \\
\hline

$13$ & $\SU(4)$ & $\mu_4\oplus\rho_6$ &
$\SU(3)$ & $\Sp(2)$ & $\SU(2)$& $( 5, 7)$ & $\D_2$ \\

  & $\U(4)$ & $\mu_4\oplus\rho_6$ &
$\S^1\cdot\SU(3)$ & $\S^1\cdot\Sp(2)$ & $\S^1\cdot\SU(2)$&  &  \\
\hline

$19$ & $\SU(5)$ & $\Lambda^2\mu_5$ &
$\Sp(2)$ & $\SU(2)\SU(3)$ & $\SU(2)^2$ & $( 4, 5 )$ & $\D_4$\\

& $\U(5)$ & $\Lambda^2\mu_5$ &
$\S^1\cdot\Sp(2)$ & $\S^1\cdot\SU(2)\SU(3)$ & $\S^1\cdot\SU(2)^2$ &  &  \\
\hline

$31$ & $\Spin(10)$ & $\Delta_{10}^\pm $ &
$\SU(5)$ & $\Spin(7)$ & $\SU(4)$& $( 9,6 )$ & $\D_4$ \\

  & $\S^1\cdot\Spin(10)$ & $\Delta_{10}^\pm $ &
$\Delta\S^1\cdot\SU(5)$ & $\S^1\cdot\Spin(7)$ & $\S^1\cdot\SU(4)$ & &\\
\hline\hline

$7$ & $\SU(3)$ & $ad$ &
$\S(\U(2)\U(1))$ & $\S(\U(1)\U(2))$ & $\T^2$& $( 2,2 )$ & $\D_3$ \\
\hline

$9$ & $\SO(5)$ & $ad$ &
$\U(2)$ & $\SO(3)\SO(2)$ & $\T^2$& $(2 ,2 )$ & $\D_4$ \\
\hline

$13$ & $\G_2$ & $ad$ &
$\U(2)$ & $\U(2)$ & $\T^2$ & $( 2,2 )$ & $\D_6$\\
\hline

$13$ & $\Sp(3)$ & $\psi_{14}$ &
$\Sp(2)\Sp(1)$ & $\Sp(1)\Sp(2)$ & $\Sp(1)^3$& $( 4,4 )$ & $\D_3$ \\
\hline

$25$ & $\F_4$ & $\psi_{26}$ &
$\Spin(9)$ & $\Spin(9)$ & $\Spin(8)$& $(8 ,8 )$ & $\D_3$ \\
\hline
  \hline

          \end{tabular}
      \end{center}
      \vspace{0.1cm}
      \caption{Essential cohomogeneity one actions and extensions on
      $\Sph^{2n+1}$}\label{odd}
\end{table} }

\vspace{200pt}
\newpage
\clearpage

{\setlength{\tabcolsep}{0.10cm}
\renewcommand{\arraystretch}{1.6}
\stepcounter{equation}
\begin{table}[!h]
      \begin{center}
          \begin{tabular}{|c||c|c|c|c|c|c|c|}
\hline
         &$\G$    &$\Km$    &$\Kp$     &$\H$ & $(l_-,l_+)$ & $W$   \\
\hline \hline

$\Sph^4$ & $\SO(3)$ & $\S(\O(2)\O(1))$ & $\S(\O(1)\O(2))$ &
$\Z_2\oplus\Z_2$ & $( 1,1 )$ & $ \D_3 $
\\ \hline

$\Sph^{14}$ & $\Spin(7)$  & $\Spin(6))$ & $\G_2$ & $\SU(3)$ & $( 1,6
)$ & $\D_2  $
\\ \hline\hline

$\CP^{k+1}$ & $\SO(k+2)$  & $\SO(2)\SO(k)$ & $\O(k+1)$ &
$\Z_2\cdot\SO(k)$ & $(
1,k)$ & $ \D_2 $  \\

 \hline

$\CP^{2k+1}$ & $\SU(2)\SU(k+1)$  & $ \Delta\SU(2)\S^1\SU(k-1)  $ & $
\S^1\U(k)$ & $\T^2\SU(k-1)$ &
$( 2 ,2k-1 )$ & $\D_2  $  \\

\hline

$\CP^6$ & $\G_2$ &
$\U(2)$ & $\Z_2\cdot\SU(3)$ & $\Z_2\cdot\SU(2)$ & $( 1,5 )$ & $ \D_2 $  \\
\hline

$\CP^7$ & $\Spin(7)$ & $\S^1\SU(3)$ & $\Z_2\cdot\Spin(6)$ &
$\Z_2\cdot\SU(3)$ & $(1 ,7 )$ &
$ \D_2 $  \\
\hline

$\CP^9$ & $\SU(5)$  &
$\S^1\cdot\Sp(2)$ & $\S(\U(2)\U(3))$ & $\S^1\cdot\SU(2)^2$ & $( 4,5
)$ & $\D_2  $  \\
\hline

$\CP^{15}$ & $\Spin(10)$  &
$\S^1\SU(5)$ & $\S^1\Spin(7)$ & $\S^1\SU(4)$ & $( 9,6 )$ & $ \D_2 $  \\
\hline\hline

$\HP^{k+1}$ & $\SU(k+2)$  & $\SU(2)\SU(k)$ & $\U(k+1)$ &
$\S^1\cdot\SU(k)$ & $(
2,2k+1 )$ & $\D_2  $  \\
& $\S^1\SU(k+2)$  & $\Delta\S^1\SU(2)\SU(k)$ & $\S^1\U(k+1)$ &
$\T^2\cdot\SU(k)$ &
&\\

\hline\hline

$\CaP$ & $\Sp(3)$   &
$\Sp(2)$ & $\Sp(1)\Sp(2)$ & $\Sp(1)^2$  & $( 11,8 )$ & $ \D_2 $ \\

& $\S^1\Sp(3)$   & $\Delta\S^1\Sp(2)$ & $\S^1\Sp(1)\Sp(2)$ &
$\S^1\Sp(1)^2$ &  & \\

 & $\Sp(1)\Sp(3)$   &
$\Delta\Sp(1)\Sp(2)$ & $\Sp(1)^2\Sp(2)$ & $\Sp(1)^3$ &  &  \\
\hline

  \hline

          \end{tabular}
      \end{center}
      \vspace{0.1cm}
      \caption{Essential cohomogeneity one actions and extensions
       in even dimensions}\label{even}
\end{table}}

\stepcounter{equation}
\begin{table}[!h]
      \begin{center}
           \begin{tabular}{|c|c|}
\hline $\G$          &$\K$   \\
\hline \hline

$\SO(2n)$ &$\SO(2k)\SO(2n-2k) \, , \, \U(n)$  \\
\hline

$\SO(2n+1)$ &$\SO(2k+1)\SO(2n-2k)$  \\
\hline

$\SU(n)$ &$\S(\U(k)\U(n-k))$  \\
\hline

$\Sp(n)$ &$\Sp(k)\Sp(n-k) \, , \, \U(n)$  \\
\hline

$\G_2$ &$\SO(4)$  \\
\hline

$\F_4$ &$\Spin(9)\, , \, \Sp(3)\Sp(1)$  \\
\hline

$\E_6$ &$\SU(6)\SU(2) \, , \, \Spin(10)\cdot \S^1$  \\
\hline

$\E_7$ &$\SU(8) \, , \,
\Spin(12)/\Z_2\cdot \S^3 \, ,  \, \E_6\cdot \S^1$  \\
\hline

$\E_8$ &$\Spin(16)/\Z_2 \, , \, \E_7\cdot\S^3 $  \\
\hline

\hline
          \end{tabular}
      \end{center}
      \vspace{0.1cm}
      \caption{Equal rank symmetric subgroups}\label{symm}
\end{table}

\bigskip

\bigskip

\bigskip

\bigskip

\allowdisplaybreaks

\pagebreak

\providecommand{\bysame}{\leavevmode\hbox
to3em{\hrulefill}\thinspace}


\begin{thebibliography}{99999}

\bibitem[AA]{AA} A.V.~Alekseevsy and D.V.~Alekseevsy,
{\em $G$- manifolds with one dimensional orbit space,} Ad. in Sov.
Math. {\bf 8} (1992), 1--31.

\bibitem[AW]{AW} S.~Aloff and N.~Wallach,
{\em An infinite family of 7--manifolds admitting positively curved
Riemannian structures,} Bull. Amer. Math. Soc. {\bf 81}(1975),
93--97.

\bibitem[BH]{BH} A.~Back and W.Y.~Hsiang,
{\em Equivariant geometry and Kervaire spheres,} Trans. Amer. Math.
Soc. {\bf 304} (1987), no. 1, 207--227.

\bibitem[Ba]{Ba} Y.V.~Bazaikin, {\em On a certain family of closed
13--dimensional Riemannian manifolds of positive curvature,} Sib.
Math. J. {\bf 37}, No. 6 (1996), 1219-1237.


\bibitem[BB]{Bb}
      L.~B\'erard Bergery,
      {\em  Les vari\'et\'es riemanniennes homog\`enes simplement connexes
            de dimension impaire \`a courbure strictement positive},
            J.\ Math.\ pure et appl.\ {\bf 55}\,(1976), 47--68.

\bibitem[Be]{Be} M.~Berger, {\em Les varietes riemanniennes
     homogenes normales simplement connexes a
Courbure strictment positive}, Ann. Scuola Norm. Sup. Pisa {\bf 15}
(1961), 191-240.

\bibitem[Bo]{borel} A.~Borel,
{\em Sous-groupes commutatifs et torsion des groupes de Lie compacts
connex\'es}, Tohuko Math. J. {\bf 3}
(1961), 216--240.

\bibitem[BG]{BG}
C. P.~Boyer and K.~Galicki \emph{ 3-Sasakian manifolds}, Surveys in
differential geometry: Essays on Einstein manifolds, Surv. Differ.
Geom., VI, Int. Press, Boston, MA, (1999), 123--184.

\bibitem[BGM]{BGM}
C. P.~Boyer, K.~Galicki and B. M.~Mann {\em The Geometry and
Topology of 3-Sasakian Manifolds}, J. f\"{u}r die reine und angew.
Math. {\bf 455} (1994), 183-220.

\bibitem[BGMR]{BGMR}
C. P.~Boyer, K.~Galicki, B. M.~Mann and E.G.~Rees {\em Compact
3-Sasakian 7-Manifolds with Arbitrary Second
 Betti Number}, Inv. Math. {\bf 131} (1998), 321-344.


\bibitem[Br]{Br}
G.E.~Bredon, \emph{Introduction to compact transformation groups},
Academic
  Press, New York, 1972, Pure and Applied Mathematics, Vol. 46.



\bibitem[De1]{De1}
O.~Dearricott, \emph{Positive sectional curvature on 3-Sasakian manifolds},
Ann. Global Anal. Geom. \textbf{25} (2004), 59--72.

\bibitem[De2]{De2}
O.~Dearricott, \emph{Positively curved self-dual Einstein metrics on
weighted projective spaces}, Ann. Global Ana. Geom. \textbf{27}
(2005), 79--86.

\bibitem[DE]{DE}
O.~Dearricott and J.~Eschenburg, \emph{Totally geodesic embeddings
of 7-manifolds in positively curved 13-manifolds}, Manuscripta Math.
\textbf{114} (2004), 447--456.



\bibitem[E1]{Es1} J.H.~Eschenburg, {\em New examples of manifolds with
strictly positive curvature,} Inv. Math {\bf 66} (1982), 469-480.

\bibitem[E2]{Es2} J.H.~Eschenburg,
     {\em Freie isometrische Aktionen auf kompakten Lie-Gruppen
     mit positiv gekr\"ummten Orbitr\"aumen,}
     Schriftenr.~Math.~Inst.~Univ.~M\"unster {\bf 32} (1984).




\bibitem[FR1]{fang-rong:finiteness}
F.~Fang and X.~Rong, \emph{Positive pinching, volume and second
Betti number}, Geom. Funct. Anal. \textbf{9} (1999), 641--674.


\bibitem[FR2]{fang-rong:sym}
F.~Fang and X.~Rong, \emph{Positively curved manifolds with maximal discrete
symmetry rank}, Amer. Math. J.
\textbf{126} (2004), 227--245.

\bibitem[FR3]{fang-rong:dim8and9}
F.~Fang and X.~Rong, \emph{Homeomorphism classification of
positively curved manifolds with almost maximal symmetry rank,}
Math. Ann.  \textbf{332} (2005), 81--101.

\bibitem[Fu]{fukaya}
K.~Fukaya, \emph{Hausdorff convergence of Riemannian manifolds and its
applications.}  Recent topics in differential and analytic geometry,
 Adv. Stud. Pure Math., \textbf{18-{\rm I}}, Academic Press, Boston, MA
(1990), 143--238.

\bibitem[GKS]{GKS}
S.~Goette, N.~Kitchloo and K.~Shankar, \emph{Diffeomorphism type of
the Berger space ${\rm SO}(5)/{\rm SO}(3)$}, Amer. Math. J.
\textbf{126} (2004), 395--416.


\bibitem[Gr]{grove:survey} K.~Grove,
{\em Geometry of, and via, Symmetries}, Amer. Math. Soc. Univ.
Lecture Series {\bf 27} (2002), 31--53.

\bibitem[GK]{grove-kim}
K.~Grove and C.W.~Kim, \emph{Positively curved manifolds with low
fixed point cohomogeneity}, J. Diff. Geom. \textbf{67} (2004),
1--33.


\bibitem[GS1]{grove-searle:rank}
K.~Grove and C.~Searle,
\newblock Positively curved manifolds with maximal symmetry-rank,
\newblock {\em J. Pure Appl. Algebra} \textbf{91} (1994), 137--142.





\bibitem[GS2]{grove-searle:rep} K.~Grove and C.~Searle,
{\em Differential topological restrictions by curvature and
symmetry,}
       J. Diff. Geom. {\bf 47} (1997), no. 3, 530--559.




\bibitem[GS3]{grove-searle:core} K.~Grove and C.~Searle,
{\em  Global G-manifold Reductions and Resolutions}, Ann. Global
Anal. Geom. {\bf 18} (2000), 437--446.


\bibitem[GSZ]{GSZ}
K. ~Grove, K.~Shankar and W.~Ziller, \emph{Symmetries of Eschenburg
spaces and the Chern Problem}, in preparation.




\bibitem[GVWZ]{kervaire}
K.~Grove, L.~Verdiani, B.~Wilking and W.~Ziller, \emph{Obstructions
to nonnegative curvature in cohomogeneity one}, in preparation.




\bibitem[GZ]{GZ}
K.~Grove and W.~Ziller, \emph{Curvature and symmetry of {M}ilnor
spheres}, Ann. of Math. \textbf{152} (2000), 331--367.





\bibitem[Hi1]{Hi1}
N.~Hitchin, {\emph A new family of Einstein metrics, }Manifolds and
geometry (Pisa, 1993), 190--222, Sympos. Math., XXXVI, Cambridge
Univ. Press, Cambridge, 1996.

\bibitem[Hi2]{Hi2}
N.~Hitchin, {\emph Poncelet Polygons and the Painleve equations},
Proc of TATA Institue Conference 1991.



\bibitem[HK]{hsiang-kleiner}
W.Y.~Hsiang and B.~Kleiner, \emph{On the topology of positively
curved $4$-manifolds with symmetry},
\newblock J. Diff. Geom. \textbf{29} (1989), 615--621.


\bibitem[HL]{HL}
W.Y.~Hsiang and B.~Lawson, \emph{Minimal submanifolds of low
  cohomogeneity}, J. Diff. Geom. \textbf{5} (1971), 1--38.


\bibitem[Iw1]{iwata:1}
K.~Iwata, \emph{Classification of compact transformation groups on
cohomology quaternion projective spaces with codimension one
orbits},  Osaka J. Math. \textbf{15} (1978), 475--508.

\bibitem[Iw1]{iwata:2}
K.~Iwata, \emph{Compact transformation groups on rational cohomology
Cayley projective planes}, Tohoku Math. J. \textbf{33} (1981),
429--442.

\bibitem[Ne]{neumann}
W.~Neumann, \emph{3-dimensional G manifolds with two dimensional
orbits}, in: P.S.Mostert (ed.) Proceedings of Conference on
Transformation Groups, Springer Verlag (1968), 220--222.

\bibitem[Pe]{Pe} A.~Petrunin, {\em Parallel transportation for
     Alexandrov spaces with curvature bounded below,} Geom. Funct. Anal.
     {\bf 8} (1998), 123-148.


\bibitem[PT]{PT}
A.~Petrunin and W.~Tuschmann, \emph{Diffeomorphism finiteness,
positive pinching, and second homotopy}, Geom. Funct. Anal.
\textbf{9} (1999),
       736--774.

\bibitem[PV1]{PV1}
F.~Podesta and L.~Verdiani, \emph{Totally geodesic orbits of
isometries}, Ann. Glob. Anal. Appl. \textbf{16} (1998), 399--412.


\bibitem[PV2]{PV2}
F.~Podesta and L.~Verdiani, \emph{Positively curved 7-dimensional
manifolds}, Quat. Math. Oxford \textbf{50} (1999), 497--504.

\bibitem[P\"{u}]{putmann}
T.~P\"{u}ttmann, \emph{Homogeneity Rank and Atoms of Actions}, Ann.
Glob. Ana. and Geom. \textbf{22} (2002), 375--399.


\bibitem[Ro]{rong:s5}
X.~Rong, \emph{Positively curved manifolds with almost maximal
symmetry rank,} Geom. Dedicata \textbf{95} (2002), 157--182.




\bibitem[Se]{searle} C.~Searle,
{\em Cohomogeneity and positive curvature in low dimensions,} Math.
Z. {\bf 214} (1993), 491--498: Err. ibet. {\bf 226} (1997),
165--167.



\bibitem[Sh]{shankar:chern}
K.~Shankar, \emph{On the fundamental groups of positively curved
manifolds,}  J. Diff. Geom. \textbf{49} (1998), 179--182.


\bibitem[St]{straume} E.~Straume,
     {\em Compact connected Lie transformation groups on spheres with
     low cohomogeneity I,} Mem. Amer. Math. Soc. {\bf 119} (1996), no.
569, vi+93 pp.


\bibitem[Ta]{Ta}
I.A.~Taimanov, \emph{A remark on positively curved manifolds of
dimensions 7 and 13,} Tensor and vector analysis, Gordon and Breach,
Amsterdam (1998), 282--295.


    \bibitem[Uc]{uchida} F.~Uchida,
{\em Classification of compact transformation groups on cohomology
complex projective spaces with codimension one orbits}, Japan J.
Math. {\bf 3} (1977), 141--189.


\bibitem[V1]{verdiani:1} L.~Verdiani, {\em Cohomogeneity one
Riemannian manifolds
     of even dimension with strictly positive sectional curvature, I,}
     Math. Z. {\bf 241} (2002), 329--339.




\bibitem[V2]{verdiani:2} L.~Verdiani, {\em Cohomogeneity one manifolds
     of even dimension with strictly positive sectional curvature,}
     J. Diff. Geom. {\bf 68} (2004), 31--72.



     \bibitem[Wa]{Wa} N.~Wallach, {\em Compact homogeneous
     Riemannian manifolds
with strictly positive curvature}, Ann. of Math. {\bf 96} (1972),
277-295.

\bibitem[WZ]{WZ}
M.Y.~Wang and W.~Ziller, \emph{Einstein metrics on principal torus
      bundles}, J. Diff. Geom. \textbf{31} (1990), 215--248.


     \bibitem[Wi1]{Wi} B.~Wilking, {\em The normal homogeneous space
$\SU(3)\times \SO(3)/\U(2)$ has positive sectional curvature,} Proc.
     of  Amer. Math. Soc.
     {\bf 127} (1999), 1191-1994.

\bibitem[Wi2]{Wi: torus} B.~Wilking,
       {\em Torus actions on manifolds of positive sectional curvature,}
       Acta Mathematica, {\bf 191} (2003), 259-297.

\bibitem[Wi3]{Wi: sym} B.~Wilking,
       {\em Positively curved manifolds with symmetry,} Ann. of Math., to
appear.

\bibitem[Wi4]{wilking:dual} B.~Wilking,
       {\em A duality theorem for singular Riemannian foliations in
nonnegative curvature}, Preprint 2005.

\bibitem[Wo]{wolf} J.~Wolf, Spaces of Constant Curvature,
 4th edition, Publish or Perish Inc., (1977).

\bibitem[Zi]{Z}
W.~Ziller, \emph{Homogeneous spaces, biquotients, and manifolds with
positive curvature}, Lecture Notes 1998, unpublished.

\end{thebibliography}
\end{document}